\documentclass[12pt]{amsart}

\usepackage{amssymb,amsbsy}
\usepackage{latexsym}
\usepackage{verbatim}
\usepackage{fullpage}

\input {cyracc.def}
\font\ququ=cmr10 scaled \magstep1
 \font\tencyr=wncyr10 %scaled \magstephalf
\font\tencyi=wncyi10 %scaled \magstephalf
\font\tencysc=wncysc10 %scaled \magstephalf
\def\rus{\tencyr\cyracc}
\def\rusi{\tencyi\cyracc}
\def\rusc{\tencysc\cyracc}

\newcommand{\re}[1]{\textrm  (\ref{#1})}

\renewenvironment{proof}
{\noindent {\sl Proof.}\quad }{\hfill
$\square$ \vskip1.1ex\noindent }

\newenvironment{proof*}
{\noindent {\sl Proof.}\quad }{\hfill
$\square$}

\renewcommand{\theequation}{\thesection .\arabic{equation}}
\renewcommand{\thesubsubsection}{\theequation .\arabic{subsubsection}}

\catcode`\@=\active
\catcode`\@=11
\def\@eqnnum{\hbox to
.01pt{}\rlap{\hskip-\displaywidth(\mathbf{\theequation})}}
\catcode`\@=12
\newenvironment{s}[1]
{ \vskip1.2ex \refstepcounter{equation}
\noindent {\bf \theequation\enspace #1.} \begin{sl}}{\end{sl}
\vskip1.1ex\noindent }

\newenvironment{rem}[1]
{ \vskip1.2ex \refstepcounter{equation}
\noindent {\bf \theequation\enspace {#1}.} }{ \vskip1.1ex\noindent }

\newcommand {\ah}{{\frak a}}
\newcommand {\be}{{\frak b}}
\newcommand {\ce}{{\frak c}}
\newcommand {\de}{{\frak d}}

\newcommand {\g}{{\frak g}}
\newcommand {\h}{{\frak h}}
\newcommand {\ka}{{\frak k}}

\newcommand {\ma}{{\frak m}}
\newcommand {\n}{{\frak n}}

\newcommand {\p}{{\frak p}}
\newcommand {\q}{{\frak q}}
\newcommand {\rr}{{\frak r}}
\newcommand {\es}{{\frak s}}
\newcommand {\te}{{\frak t}}

\newcommand {\fX}{{\frak X}}
\newcommand {\fY}{{\frak Y}}
\newcommand {\z}{{\frak z}}
\newcommand {\gln}{{\frak gl}_n}
\newcommand {\sln}{{\frak sl}_n}

\newcommand {\glv}{{\frak gl}(V)}

\newcommand {\spv}{{\frak sp}(V)}
\newcommand {\spn}{{\frak sp}_{2n}}

\newcommand {\sono}{{\frak so}_{2n+1}}
\newcommand {\sone}{{\frak so}_{2n}}
\newcommand {\son}{{\frak so}_{n}}

\newcommand {\esi}{\varepsilon}
\newcommand {\eps}{\epsilon}
\newcommand {\ap}{\alpha}

\newcommand {\lb}{\lambda}
\newcommand {\vp}{\varphi}

\newcommand {\wF}{\widehat F}

\newcommand {\ck}{{\mathcal K}}

\newcommand {\N}{{\mathcal N}}
\newcommand {\co}{{\mathcal O}}

\newcommand {\BV}{{\mathbb V}}

\newcommand {\BZ}{{\mathbb Z}}
\newcommand {\BN}{{\mathbb N}}

\newcommand {\md}{/\!\!/}
\newcommand {\isom}{\stackrel{\sim}{\longrightarrow}}

\newcommand {\ad}{{\mathrm{ad\,}}}
\newcommand {\ads}{{\mathrm{ad}^*}}
\newcommand {\Ad}{{\mathrm{Ad\,}}}

\newcommand {\codim}{{\mathrm{codim\,}}}

\newcommand {\ed}{{\mathrm{edim\,}}}

\newcommand {\hot}{{\mathrm{ht\,}}}
\newcommand {\ind}{{\mathrm{ind\,}}}
\newcommand {\Lie}{{\mathrm{Lie\,}}}

\newcommand {\Ima}{{\mathrm{Im\,}}}
\newcommand {\Mor}{\operatorname{Mor}}

\newcommand {\rk}{{\mathrm{rk\,}}}
\newcommand {\spe}{{\mathrm{Spec\,}}}

\newcommand {\trdeg}{{\mathrm{trdeg\,}}}

\newcommand {\tri}{{\frak sl}_2}
\newcommand {\GR}[2]{{\textrm{{\bf #1}}}_{#2}}

\newcommand {\ov}{\overline}
\newcommand {\un}{\underline}

\newcommand {\nv}{{\goth N}(V)}
\newcommand {\ntg}{{\goth N}(\hat\g)}
\newcommand {\ntv}{{\goth N}(\hat V)}

\newcommand {\qus}{\hfill $\square$ \vskip1.1ex}

\newcommand {\beq}{\begin{equation}}
\newcommand {\eeq}{\end{equation}}

\renewcommand{\le}{\leqslant}
\renewcommand{\ge}{\geqslant}

\newfam\Bbbfam\newfam\eufam\newfam\eusfam%
\font\Bbbfont=msbm10 scaled 1200%
\font\bbbfont=msbm10 scaled 1000%
\font\olala=msam10 scaled 1200%
\font\frak=eufm10 scaled 1400%
\font\Bbbsmallfont=msbm8%
\textfont\Bbbfam=\Bbbfont\scriptfont\Bbbfam=\Bbbsmallfont%%
\font\euzw=eufm10 scaled 1200%
\font\euac=eufm7 scaled 1200%
\font\euacc=eufm7 scaled 1000%
\textfont\eufam=\euzw\scriptfont\eufam=\euac%
\scriptscriptfont\eufam=\euacc%
\font\euszw=eusm10 scaled 1200%
\font\eusac=eusm7 scaled 1200%
\font\eusacc=eusm7 scaled 1000%
\textfont\eusfam=\euszw\scriptfont\eusfam=\eusac%
\scriptscriptfont\eusfam=\eusacc%
\def\frak{\fam\eufam}%
\def\goth{\fam\eusfam}%

\def\varnothing{\hbox {\Bbbfont\char'077}}
\def\square{\hbox {\olala\char"03}}
\def\bbk{\hbox {\Bbbfont\char'174}}
\def\bkk{\hbox {\bbbfont\char'174}}

\begin{document}
\setlength{\parskip}{2pt plus 4pt minus 0pt}
%%{\scriptsize Preliminary version!} 
\hfill {\scriptsize October 17, 2007}%% 
\vskip1ex
\vskip1ex

\title[Semi-direct products of Lie algebras]{Semi-direct products of Lie algebras 
and their invariants}
\author{Dmitri I. Panyushev}
\thanks{This research was supported in part by 
RFBI Grants 05-01-00988 and 06-01-72550}
\subjclass[2000]{14L30, 17B20, 22E46}
\maketitle
\begin{center}
{\footnotesize
{\it Independent University of Moscow,
Bol'shoi Vlasevskii per. 11 \\
119002 Moscow, \quad Russia \\ e-mail}: {\tt panyush@mccme.ru }\\
}
\end{center}

\tableofcontents  
\section*{Introduction}

\noindent
The ground field $\bbk$ is algebraically closed and of characteristic zero.
The goal of this paper is to extend the standard invariant-theoretic design,
well-developed in the reductive case, to the setting of non-reductive group 
representations. This concerns the following notions and results:
the existence of generic stabilisers and generic isotropy groups for
(finite-dimensional rational) representations; structure of the fields and
algebras of invariants; quotient morphisms and structure of their fibres.
One of the main tools for obtaining non-reductive Lie algebras is the semi-direct 
product construction.
There is a number of articles devoted to the study of the coadjoint representations
of non-reductive Lie algebras; in particular, semi-direct products, see e.g. 
\cite{aif-rais,cambr,rais,rt,tayu,ksana1}. 
In this article, we consider 
such algebras from a broader point of view. 
In particular, we found that the adjoint representation is an interesting object, too.
Our main references for Invariant Theory are \cite{brion} and \cite{VP}.
All algebraic groups are assumed to be linear.

If an algebraic group $A$ acts on an affine variety $X$, then $\bbk[X]^A$ 
stands for the algebra of $A$-invariant regular functions on $X$. If $\bbk[X]^A$
%%this algebra 
is finitely generated, then $X\md A:=\spe \bbk[X]^A$, and
the {\it quotient morphism\/} $\pi_A: X\to X\md A$ is the mapping associated with
the embedding $\bbk[X]^A \hookrightarrow \bbk[X]$. If $\bbk[X]^A$ is polynomial, then
the elements of any set of algebraically independent homogeneous generators 
will be referred to as {\it basic invariants\/}.
\\[.6ex]
Let $G$ be a connected reductive algebraic group with Lie algebra $\g$.
Choose a Cartan subalgebra $\te\subset \g$ with the corresponding 
Weyl group $\mathsf W$. 
The adjoint representation $(G:\g)$ has a number of good properties,
some of which are listed below:

$\bullet$ \ The adjoint representation is self-dual, and $\te$ is a generic stabiliser 
for it;

$\bullet$ \ The algebra of invariants $\bbk[\g]^G$ is polynomial;

$\bullet$ \ the restriction homomorphism $\bbk[\g]\to\bbk[\te]$ induces the isomorphism
$\bbk[\g]^G \simeq \bbk[\te]^{\mathsf W}$ (Chevalley's theorem);

$\bullet$ \ The quotient morphism $\pi_G:\g\to \g\md G$ is equidimensional and the 
fibre of the
origin, $\N:=\pi_G^{-1}(\pi_G(0))$, is an irreducible complete intersection.
The ideal of $\N$ in $\bbk[\g]$ is generated by the basic invariants;

$\bullet$ \ $\N$ is the union of finitely many $G$-orbits.
\\[.6ex]
Each of these properties may fail if $\g$ is replaced with an arbitrary
algebraic Lie algebra $\q$. In particular, one have to distinguish the adjoint
and coadjoint representations of $\q$. As usual, $\ad$ (resp. $\ads$) stands
for the adjoint (resp. coadjoint) representation.
Write $Q$ for a connected group with Lie algebra $\q$. 

First, we consider the problem of existence of generic
stabilisers for $\ad$ and $\ads$.
(See Section~\ref{prelim} for precise definitions).
It turns out that 
if $(\q,\ad)$ has a generic stabiliser, say $\h$, 
then $\h$ is commutative and $\n_\q(\h)=\h$.
This yields a Chevalley-type theorem for the fields of invariants:
$\bbk(\q)^Q\simeq \bbk(\h)^{W}$, where $W=N_Q(\h)/Z_Q(\h)$ is finite. 
We also notice that 
$(\q,\ad)$ has a generic stabiliser if and only if the Cartan subalgebras of
$\q$ are commutative.
If $(\q,\ads)$ has a generic stabiliser, say $\h$, then $\h$ is commutative,
$\dim N_Q(\h)=\dim (\q^*)^\h$,
and $\bbk(\q^*)^Q\simeq \bbk((\q^*)^\h)^{N_Q(\h)}$. 
But unlike the adjoint case, the action $(N_Q(\h):(\q^*)^\h)$
does not necessarily reduce to a finite group action.
We prove that under a natural constraint 
the representation of the identity component of $N_Q(\h)$ on $(\q^*)^\h$ is
the coadjoint representation. 

Our main efforts are connected with the following situation. Suppose that
$(\q,\ad)$ or $(\q,\ads)$ has some of the above good properties and $V$ is a 
(finite-dimensional rational) $Q$-module. Form the Lie algebra $\q\ltimes V$. 
It is the semi-direct product of $\q$ and $V$, $V$ being a commutative ideal
in it. The corresponding connected algebraic group is $Q\ltimes V$.
(See section~\ref{covariants} for the details.)
Then we want to realise to which extent those good properties are preserved under this
procedure. This surely depends on $V$, 
and we are essentially interested in two cases:

(a) \ $\q$ is arbitrary and $V=\q$ or $\q^*$
(the adjoint or coadjoint $\q$-module);

(b) \ $\q=\g$ is reductive and $V$ is an arbitrary $G$-module.
\\[.6ex]
For (a), we prove that if $(\q,\ad)$ has a generic stabiliser, then so do
$(\q\ltimes \q,\ad)$ and $(\q\ltimes \q^*,\ad)$.
Furthermore, the passages $\q\leadsto \q\ltimes \q$ and
$\q\leadsto \q\ltimes \q^*$ does not affect the generalised
Weyl group $W$, and both fields $\bbk(\q\ltimes \q)^{Q\ltimes\q}$ and
$\bbk(\q\ltimes \q^*)^{Q\ltimes\q^*}$ are purely transcendental extensions
of $\bbk(\q)^Q$. It is also true that if $(\q,\ads)$ has a generic stabiliser, then
so does $(\q\ltimes\q,\ads)$.

For (b), we prove that $(\g\ltimes V,\ad)$ always has a generic stabiliser.
But this is not the case for $\ads$. Recall that any $\g$-module $V$
has a generic stabiliser. The following result seems to be quite unexpected.
Suppose generic $G$-orbits in $V$ are closed (i.e., the action
$(G:V)$ is {\it stable}), then $(\g\ltimes V,\ads)$  has a generic stabiliser
if and only if $V$ is a {\it polar\/} $G$-module in  the sense of \cite{polar}.
The assumption of stability is relatively harmless, since there are only
finitely many $G$-modules without that property.
On the other hand, the hypothesis of being polar is quite restrictive,
because for any $G$ there are only finitely many polar representations.

One of our main observations is that there are surprisingly many 
nonreductive Lie algebras $\ah$ and $\ah$-modules $M$
such that $\bbk[M]^{A}$ is a polynomial algebra. Furthermore, the 
basic invariants of $\bbk[M]^{A}$ can explicitly be constructed using certain 
modules of covariants. This concerns the following cases:

-- If $\g$ is reductive and $V$ is an arbitrary $\g$-module, then one takes
$\ah=M=\g\ltimes V$;

-- If the action $(Q:V)$ satisfies some good properties, then one
takes $\ah=\q\ltimes \q$ and $M=V\ltimes V$. 
Furthermore, the passage $(\q,V)\mapsto (\hat\q=\q\ltimes\q, \hat V=V\ltimes V)$ 
can be iterated.
\\[.7ex]
The precise statements are given below.

\begin{s}{Theorem}  \label{thm-intro-2}  
Let\/ $V$ be an arbitrary $G$-module. Set\/ $\q=\g\ltimes V$, $Q=G\ltimes V$,
and\/ $m=\dim V^\te$. Notice that $1\ltimes V$ is a commutative normal subgroup 
of $Q$ (in fact, the unipotent radical of\/ $Q$). Then
\begin{itemize}
\item[\sf (i)] \  $\bbk[\q]^{1\ltimes V}$ is a polynomial algebra of Krull
dimension  $\dim\g+m$.
It is freely generated by the coordinates on $\g$ and the functions
$\wF_i$, $i=1,\ldots,m$, associated with covariants of type $V^*$.
\item[\sf (ii)] \  $\bbk[\q]^{Q}$ is a polynomial algebra of
Krull dimension $\dim\te+m$.
It is freely generated by the basic invariants of\/ $\bbk[\g]^G$ and the same
functions $\wF_i$, $i=1,\ldots,m$. 
\item[\sf (iii)] \  $\max\dim_{x\in \q} Q{\cdot}x=\dim\q-\dim \q\md Q$;
\item[\sf (iv)] \ If\/ $\pi:\q\to \q\md Q$ is the quotient morphism and
$\Omega:=\{x\in\q \mid \text{d}\pi_x\ \text{ is onto }\}$,
then $\q\setminus \Omega$ contains no divisors.
\end{itemize}
\end{s}%
Given a $\q$-module $V$, the space $V\times V$ can be regarded as $\q\ltimes\q$-module
in a very natural way. Write $\hat V$ or $V\ltimes V$ for this module.

\begin{s}{Theorem}  \label{thm-intro-1}  
Suppose the action\/ $(Q:V)$ satisfies the following conditions:

(1) \ $\bbk[V]^Q$ is a polynomial algebra;

(2) \ $\max\dim_{v\in V} Q{\cdot}v=\dim V-\dim V\md Q$;

(3) \ If\/ $\pi_Q: V\to V\md Q$ is the quotient morphism and
$\Omega:=\{v\in V \mid (\text{d}\pi_Q)_v\ \text{ is onto }\}$, then\/ 
$V\setminus\Omega$ contains no divisors.
\\
Set\/ $\hat\q=\q\ltimes\q$ and $\hat Q=Q\ltimes\q$. 
Then
\begin{itemize}
\item[\sf (i)] \  $\bbk[\hat V]^{1\ltimes\q}$ is a polynomial algebra of Krull
dimension\/ $\dim V+\dim V\md Q$, which is generated by the coordinates on the first
factor of $\hat V$ and the polynomials $\wF_1,\ldots,\wF_m$ associated with
the differentials of basic invariants in $\bbk[V]^Q$;
\item[\sf (ii)] \ $\bbk[\hat V]^{\hat Q}$ is a polynomial algebra of Krull
dimension\/ 2$\dim V\md Q$, which is freely generated by the basic invariants 
of\/ $\bbk[V]^Q$ and the same functions $\wF_i$, $i=1,\ldots,m$.
\item[\sf (iii)] \ The $\hat Q$-module $\hat V$ satisfies conditions (1)--(3), too.
\end{itemize}
\end{s}%
Since the adjoint representation of a reductive Lie algebra $\g$ satisfies the above 
properties (1)-(3),
one may begin with $\q=\g=V$, and iterate the procedure ad infinitum.
For the adjoint representation of a semisimple Lie algebra, 
the assertion in part (ii) is due to Takiff \cite{takiff}.
For this reason Lie algebras of the form $\q\ltimes\q$ are called {\it Takiff (Lie)
algebras.} We will also say that the $\hat\q$-module $\hat V$ is the 
{\it Takiffisation\/} of the $\q$-module $V$.
But $(\g, \ad)$ is not the only possible point of departure for the infinite 
iteration process.
In view of Theorem~\ref{thm-intro-2},
the algebras $\q=\g\ltimes V$ and their adjoint representations
can also be used as initial bricks in the 
Takiffisation procedure.
%% of  Theorem~\ref{thm-intro-1}.

If $\bbk[V]^{Q}$ is polynomial, then it is natural to study the fibres of
the quotient morphism $\pi_{Q}$. The {\it null-cone}, ${\goth N}(V)=\pi_{Q}^{-1}
(\pi_{Q}(0))$, is the most 
important fibre. For instance,  $\bbk[V]$ is a free 
$\bbk[V]^{Q}$-module if and only $\dim {\goth N}(V)=
\dim V-\dim V\md {Q}$, i.e., $\pi_{Q}$ is equidimensional.
We consider properties of null-cones arising in the context of semi-direct
products and their representations.
%% of Theorems~\ref{thm-intro-1} and \ref{thm-intro-2}.

For $\q=\g\ltimes V$, as in Theorem~\ref{thm-intro-2}, a necessary and sufficient
condition for the equidimensionality of $\pi_Q$
is stated in terms of a stratification
of $\N$ determined by the covariants on $\g$ of type $V^*$. 
Using this stratification and some technique from \cite{jac} and \cite{ls},
we prove the following:

If ${\goth N}(\q)$ is irreducible, then 
(i) \ $\pi_Q$ is equidimensional; (ii) \ the morphism $\kappa: \q\to \q$ defined by
$\kappa(x,v)=(x,x{\cdot}v)$, $x\in\g,v\in V$, has the property that the closure of
$\Ima(\kappa)$ is a factorial complete intersection and its ideal in $\bbk[\q]$
is generated
by the polynomials $\wF_i$, $i=1,\ldots,m$, mentioned in Theorem~\ref{thm-intro-2}.
This is a generalisation of \cite[Prop.\,2.4]{ls}. Similar results hold for 
the Takiffisation of $G$-modules $V$ having good properties, as in 
Theorem~\ref{thm-intro-1}.
In this case, conditions of equidimensionality for 
$\pi_{\hat G}: \hat V\to \hat V\md \hat G$
are stated in terms of a stratification of ${\goth N}(V)$ 
determined by the covariants on $V$ of type $V^*$.
See Section~\ref{null-cone} for the details.

In general, it is difficult to deal with the stratifications of $\N$ and 
${\goth N}(V)$, 
but, for isotropy contractions and $\BZ_2$-contractions of reductive
Lie algebras, explicit results can be obtained.
Let $\h$ be a reductive subalgebra of $\g$ and $\g=\h\oplus\ma$  a direct sum of $\h$-modules.
Then $\h\ltimes\ma$ is called an {\it isotropy contraction\/} of $\g$.
If $\g=\h\oplus\ma$ is a $\BZ_2$-grading, then we say about a $\BZ_2$-{\it contraction}.
(The word ``contraction'' can be understood in the usual sense of deformation theory
of Lie algebras.)
Semi-direct products occurring in this way have some interesting properties. As a sample, 
we mention the following useful fact: $\ind (\h\ltimes\ma)= \ind\g + 2c(G/H)$, where 
$\ind(.)$ is the index of a Lie algebra and 
$c(.)$ is the complexity of a homogeneous space. 
In particular, $\ind(\h\ltimes\ma)=\ind\g$ 
if and only if $H$ is a spherical subgroup of $G$.

Our main results on the equidimensionality of quotient morphisms and
irreducibility of null-cones are related to the $\BZ_2$-contractions of simple
Lie algebras.
Given a $\BZ_2$-grading $\g=\g_0\oplus\g_1$, Theorem~\ref{thm-intro-2} applies
to the semi-direct product $\ka=\g_0\ltimes \g_1$, so that $\bbk[\ka]^K$ is a polynomial
algebra of Krull dimension $\rk\g$. Using the classification of $\BZ_2$-gradings,
we prove that ${\goth N}(\ka)$ is irreducible.
Therefore the good properties discussed in a preceding paragraph
hold for the morphism $\kappa:\ka\to \ka$, $\kappa(x_0,x_1)=(x_0, [x_0,x_1])$. 
Our proof of irreducibility of ${\goth N}(\ka)$ basically reduces to
the verification of certain inequality for the nilpotent $G_0$-orbits in $\g_0$.
Actually, we notice that one may prove a stronger constraint (cf. 
inequalities~\re{old-brilliant} and \re{brilliant}). This leads to the following
curious result: Consider $\tilde \ka=\g_0\ltimes(\g_1\oplus \g_1)$. 
(In view of Theorem~\ref{thm-intro-2}, $\bbk[\tilde \ka]^{\tilde K}$ is polynomial.) 
Then $\pi_{\tilde K}$ is still equidimensional, although ${\goth N}(\tilde\ka)$ can
already be reducible.

To discuss similar results for the Takiffisation of $\q$-modules,
i.e., $\hat \q$-modules $\hat V$, 
one has to impose more constraints on $V$. We also assume below
that $\q=\g$ is reductive.

\begin{s}{Theorem}  \label{thm-intro-3}
Suppose the $G$-module $V$ satisfies conditions (1)--(3) of Theorem~\ref{thm-intro-1}
and also the following two conditions:

(4) \ $\nv:=\pi_G^{-1}(\pi_G(0))$ consists of finitely many $G$-orbits; 

(5) \  ${\nv}$ is irreducible and has only rational
singularities. 
\\
For $\pi_{\hat G}: \hat V\to \hat V\md \hat G$ and 
$\ntv=\pi_{\hat G}^{-1}(\pi_{\hat G}(0))$, we then have, in addition to the conclusions
of Theorem~\ref{thm-intro-1},
\begin{itemize}
\item[\sf (i)] \ $\ntv$ is an irreducible complete intersection 
and the ideal
of\/ $\ntv$ in $\bbk[\hat V]$ is generated by the basic invariants
in\/ $\bbk[\hat V]^{\hat G}$;
\item[\sf (ii)] \  $\pi_{\hat G}$ is equidimensional
and\/ $\bbk[\hat V]$ is a free $\bbk[\hat V]^{\hat G}$-module.
\end{itemize}
\end{s}%
For $G$ semisimple, conditions (2) and (3) are satisfied for all $V$, 
therefore the most essential conditions are (4) and (5). The main point here is to prove
the irreducibility.
The crucial step in proving this theorem is the use of the Goto-Watanabe inequality
\cite[Theorem\,2']{nw}
which relates the dimension and embedding dimension of the local rings that are 
complete intersections with only rational singularities, see Section~\ref{tt}.
(We refer to \cite{hess2} for the definition of rational singularities.)
For $V=\g$, the idea of using that inequality is due to M.\,Brion.
The irreducibility of $\ntg$  was first proved by 
F.\,Geoffriau \cite{geof1} via case-by-case checking. Then, 
applying  the Goto-Watanabe inequality, Brion found 
a conceptual proof of Geoffriau's result~\cite{brion2}.
Our observation is that Brion's idea applies in a slightly more general 
setting of the Takiffisation of representations $(G:V)$ satisfying conditions (1)--(5).

The irreducibility of $\ntg$ is equivalent to that a certain inequality holds
for all non-regular nilpotent elements (orbits). Here is it:
\[
    \dim\z_\g(x)+\rk(\textsl{d}\pi_G)_x > 2\rk\g
\ \ \text{ if } \ x\in \N\setminus \N^{reg} \ .
\] 
Using case-by-case checking, we prove a stronger inequality
\[
  \dim\z_\g(x)+2\rk(\text{d}\pi_G)_x- 3\rk\g\ge 0 
\ \ \text{ for all } \ x\in \N\ .
\] 
It seems that the last inequality is more fundamental, because it is stated more
uniformly, can be written in different equivalent forms, and has geometric 
applications. For instance, if $\g=\g_0\oplus\g_1$ is a $\BZ_2$-grading of maximal rank
and $\hat\g_1=\g_1\ltimes\g_1$, then the equidimensionality of 
$\pi_{\hat G_0}: \hat\g_1\to \hat\g_1\md \hat G_0$ is essentially 
equivalent to the last inequality. This result cannot be deduced from
Theorem~\ref{thm-intro-3}, because ${\goth N}(\g_1)$ is not normal.
Furthermore, ${\goth N}(\hat\g_1)$ can be reducible.

Our methods also work for generalised Takiff algebras
introduced in \cite{rt}.
The vector space $\q_\infty:=\q\otimes \bbk[\mathsf T]$ has a natural Lie algebra 
structure such that $[x\otimes \mathsf T^l, y\otimes \mathsf T^k]=
[x,y]\otimes \mathsf T^{l+k}$.
Then $\q_{\ge (n+1)}=\displaystyle\bigoplus_{j\ge n+1} \q\otimes \mathsf T^j$ is an 
ideal of $\q_\infty$, and the respective quotient
is a {\it generalised Takiff Lie algebra}, denoted $\q\langle n\rangle$.
Write $Q\langle n\rangle$ for the corresponding connected group.
Clearly, $\dim\q\langle n\rangle=(n+1)\dim\q$ 
and  $\q\langle 1\rangle\simeq \q\ltimes \q$.
We prove that
if $(Q:\q)$ satisfies conditions (1)--(3) of Theorem~\ref{thm-intro-1},
then the similar conclusions hold for the adjoint action
$(Q\langle n\rangle: \q\langle n\rangle)$. In particular,
$\bbk[\q\langle n\rangle]^{Q\langle n\rangle}$ is a polynomial algebra 
of Krull dimension $(n+1)\dim \q\md Q$. 

For $\q=\g$ semisimple, our methods enable us to deduce the 
equidimensionality of
$\pi_{G\langle 2\rangle} : \g\langle 2\rangle \to
\g\langle 2\rangle\md G\langle 2\rangle$ 
from the same fact related to the semi-direct product $\g\ltimes(\g\oplus\g)$.
%%is equidimensional.
%%(See Section~\ref{generalised} for the details).
%%After completing the first version of this article,
%%I learned that 
However, it was shown by Eisenbud and Frenkel that 
$\pi_{G\langle n\rangle} : \g\langle n\rangle \to
\g\langle n\rangle\md G\langle n\rangle$ is equidimensional for any $n$,
see  \cite[Appendix]{must}. Their proof exploits the interpretation of
${\goth N}(\g\langle n\rangle)$ as a jet scheme and uses
the deep result of Musta\c t\u a concerning the irreducibility of jet schemes
\cite[Theorem\,3.3]{must}. 
%%Our approach in case $n=2$ is much more elementary.

%%${\goth N}(\g\langle n\rangle)$.

\noindent
{\small
{\bf Acknowledgements.} Work on this article commenced during my visits to 
the Universit\'e de Poitiers (France) in 1996--98. I would like to thank 
Thierry Levasseur for arranging those visits, inspiring conversations, and
drawing my attention to work of Geoffriau. Thanks are also due to Michel Brion
for sharing some important insights on  Takiff algebras.
I am grateful to Sasha Premet for drawing my attention to results
of Eisenbud and Frenkel.
}

                %%%%
%%%%%%%%%%%%%%%%
%%%%%%%%%%%%%%%%   Section 1
%%%%%%%%%%%%%%%%
                %%%%

\section{Preliminaries}
\label{prelim}
\setcounter{equation}{0}

\noindent
Algebraic groups are denoted by capital Latin letters and their Lie
algebras are denoted by the corresponding lower-case Gothic letters.
The identity component of an algebraic group $Q$ is denoted by $Q^o$.
\\[.6ex]
Let $Q$ be an affine algebraic group acting regularly on an 
irreducible  variety $X$. Then $Q_x$ stands for the isotropy group
of $x\in X$. Likewise, the stabiliser of
$x$ in $\q=\Lie Q$ is denoted by $\q_x$.
We write $\bbk[X]^Q$ (resp. $\bbk(X)^Q$) for the algebra
of regular (resp. field
of rational) $Q$-invariants on $X$. A celebrated theorem of {M.\,Rosenlicht}
says that there is a dense open $Q$-stable subset $\tilde\Omega\subset X$
such that $\bbk(X)^H$ separates the $Q$-orbits in $\tilde\Omega$, see e.g.
\cite[1.6]{brion},\,\cite[2.3]{VP}.
In particular, $\trdeg \bbk(X)^Q=\dim X-\max\dim_{x\in X} Q{\cdot}x$.
We will use Rosenlicht's theorem in the following equivalent form:

\begin{s}{Theorem}  \label{rosenlicht}
Let\/ $\mathbb F$ be a subfield of $\bbk(X)^Q$. Then 
$\mathbb F=\bbk(X)^Q$ if and only if\/ $\mathbb F$ separates the 
$Q$-orbits in a dense open subset of $X$.
\end{s}%
We say that the action $(Q:X)$ {\it has a generic stabiliser\/}, 
if there exists
a dense open subset $\Omega\subset X$ such that all stabilisers $\q_\xi$, $\xi\in \Omega$,
are $Q$-conjugate. Then each of the subalgebras $\q_\xi$, $\xi\in\Omega$, 
is called a generic stabiliser.
The points of such an $\Omega$ are said to be {\it generic}.
Likewise, one defines  a {\it generic isotropy group}, 
which is a subgroup of $Q$. Clearly, the existence of 
a generic isotropy group implies that of a generic stabiliser.
That the converse is also true is proved by Richardson \cite[\S\,4]{Ri}.
The reader is also referred to \cite[\S 7]{VP} for a 
thorough discussion of generic stabilisers.
If $Y\subset X$ is irreducible, then
$Y^{reg}:=\{y\in Y\mid \dim Q{\cdot}y=\max_{z\in Y}\dim Q{\cdot}z\}$.
It is a dense open subset of $Y$. 
The points of $Y^{reg}$ are said to be {\it regular}.
Of course, these notions depend on $\q$. If we wish to make this dependence
explicit, we speak about $\q$-generic or $\q$-regular points.
Since $X^{reg}$ is dense in $X$, all generic points 
(if they do exist) are regular.
The converse is however not true.

If $Q$ is reductive and $X$
is smooth, then $(Q:X)$ always has a generic stabiliser \cite{Ri}. 
One of our goals is to 
study existence of generic stabilisers in case of non-reductive $Q$.
Specifically, we consider the adjoint and coadjoint representations of $Q$.
To this end, we recall some standard invariant-theoretic techniques and
a criterion for the existence of generic stabilisers.

Let $\rho: Q \to GL(V)$ be a finite-dimensional rational representation of $Q$ and
$\bar{\rho}:\q\to {\frak gl}(V)$  the corresponding representation of $\q$.
For $s\in Q$ and $v\in V$, we usually write $s{\cdot}v$ in place of $\rho(s)v$.
Similarly, $x{\cdot}v$ is a substitute for $\bar{\rho}(x)v$, $x\in \q$. 
(But for the adjoint representation, the standard bracket notation is used.)
It should be clear from the context which meaning of `$\cdot$' is meant.
%%The same notation is employed for the $\g$-module structure on $V$.
Given $v\in V$, consider 
\[
U=V^{\q_v}=\{y\in V\mid \q_v{\cdot}y=0\} \ , 
\]
the fixed point space of $\q_v$. 
Associated to $U\subset V$, there are two subgroups of $Q$:
\[
 {\mathsf N}(U)=\{s\in Q\mid s{\cdot}U \subset U\}, \quad
 {\mathsf Z}(U)=\{s\in Q\mid s{\cdot}u=u \ \text{ for all } \ u\in U\}.
\]
The following is well known and easy.
\begin{s}{Lemma}   \label{nz}
\begin{itemize}   
\item[\sf (i)] \ $\Lie {\mathsf Z}(U)=\q_v$ and 
${\mathsf Z}(U)$ is a normal subgroup of ${\mathsf N}(U)$; 
\item[\sf (ii)] \ ${\mathsf N}(U)=N_Q({\mathsf Z}(U))=N_Q(\q_v)$.
\end{itemize}
\end{s}%
It is not necessarily the case that ${\mathsf Z}(U)$ is connected; however,
${\mathsf Z}(U)$ and ${\mathsf Z}(U)^o$ have the same normaliser in $Q$.

\begin{s}{Lemma}    \label{inters}
If $y\in U^{reg}$ (i.e., $\q_y=\q_v$), then
$Q{\cdot}y\cap U={\mathsf N}(U){\cdot}y$ and $\q{\cdot}y\cap U=\n_\q(\q_v){\cdot}y$.
\end{s}\begin{proof}
1. Suppose $s{\cdot}y\in U$ for some $s\in Q$. Then $\q_{s{\cdot}y}=\q_v=\q_y$. 
Hence $s\in N_Q(\q_v)$, and we refer to Lemma~\ref{nz}.

2. Suppose $s{\cdot}y\in U$ for some $s\in \q$. Then 
$0=\q_v(s{\cdot}y)=[\q_v,s]{\cdot} y$. Hence $[\q_v,s]\subset \q_y=\q_v$.
\end{proof}%
Set $Y=\ov{Q{\cdot}U}$. It is a $Q$-stable irreducible subvariety of $V$.

\begin{s}{Proposition}  \label{field_isom}
The restriction homomorphism $(f\in \bbk(Y)) \mapsto f\vert_U$
yields an isomorphism
 $\bbk(Y)^Q\isom \bbk(U)^{{\mathsf N}(U)}=
\bbk(U)^{{\mathsf N}(U)/{\mathsf Z}(U)}$.
\end{s}\begin{proof*}
This follows from the first equality in Lemma~\ref{inters} and
Rosenlicht's theorem.
\end{proof*}%
\begin{rem}{Example}
Let $G$ be a semisimple algebraic group with Lie algebra $\g$, 
and $v=e\in\g$ a nilpotent element.
Then $\g_e=\z_\g(e)$ is the centraliser of $e$ and
$U=\{x\in\g\mid [x,\z_\g(e)]=0\}=:\de_\g(e)$ is the centre of $\z_\g(e)$.
Here ${\mathsf N}(U)=N_G(\z_\g(e))$ is the normaliser of $\z_\g(e)$ in $G$.
Letting $Y=\ov{G{\cdot}\de_\g(e)}$, we obtain an isomorphism
\[
     \bbk(Y)^G\simeq \bbk(\de_\g(e))^{N_G(\z_\g(e))} \ .
\]
It is known that $\de_\g(e)$ contains no semisimple elements \cite{rk1}, so that
$Y$ is the closure of a nilpotent orbit and hence $\bbk(Y)^G=\bbk$.
It follows that $N_G(\z_\g(e))$ has a dense orbit 
in $\de_\g(e)$. This fact was already noticed in \cite[\S\,4]{cambr}.
Actually, the dense $G$-orbit in $Y$ is just $G{\cdot}e$.
\end{rem}%
%
%We will especially be interested in the case $Y=V$.
Clearly, if $\ov{Q{\cdot}U}=V$, then $(Q:V)$ has a generic stabiliser and $v$ is a generic point.
A general criterion for this to happen is  proved in \cite[\S\,1]{alela}.
For future reference, we recall it here.

\begin{s}{Lemma {\ququ (Elashvili)}}   \label{ela_lem}
Let $v\in V$ be an arbitrary point. Then
$Q{\cdot}V^{\q_v}$ is dense in $V$ if and only if\/ 
$V=\q{\cdot}v + V^{\q_v}$. 
\end{s}%
The  existence of a non-trivial generic stabiliser  yields a Chevalley-type theorem for 
the field of invariants.
Indeed, it follows from Proposition~\ref{field_isom} that if $(Q:V)$ has a generic stabiliser, 
$v\in V$ is a generic point, and $U=V^{\q_v}$, 
then 
\begin{equation}   \label{sgp-field}
\bbk(V)^Q\simeq \bbk(U)^{{\mathsf N}(U)}=\bbk(U)^{{\mathsf N}(U)/{\mathsf Z}(U)} \ .
\end{equation}
In this context, the group $W:={\mathsf N}(U)/{\mathsf Z}(U)$ is called the 
{\it Weyl group\/}
of the action $(Q:V)$. Notice that this $W$ is not necessarily finite.

The corresponding question for the algebras of invariants is much more
subtle. The restriction homomorphism 
$f\mapsto f\vert_U$ certainly induces an embedding 
$\bbk[V]^Q \hookrightarrow \bbk[U]^{{\mathsf N}(U)/{\mathsf Z}(U)}$.
However, if $Q$ is non-reductive, then it is usually not onto.

                %%%%
%%%%%%%%%%%%%%%%
%%%%%%                                  Section 2.
%%%%%%%%%%%%%%%%
                %%%%

\section{Generic stabilisers (centralisers) for the adjoint representation}  
\label{sgp-adj}
\setcounter{equation}{0}

\noindent
In what follows, $Q$ is a connected algebraic group.
In this section, we elaborate on the existence of generic stabilisers 
and its consequences for the adjoint
representations $\Ad : Q\to GL(\q)$ and  $\ad: \q\to \mathfrak{gl}(\q)$.

For $x\in \q$, the stabiliser $\q_x$ is nothing but the {\it centraliser of $x$}
in $\q$,
so that we write $\z_\q(x)$ in place of $\q_x$. The centraliser of $x$ in $Q$ 
is denoted by $Z_Q(x)$. If $(\q,\ad)$ has a generic
stabiliser, then we also say that $\q$ {\it has a generic centraliser\/}.
By Lemma~\ref{ela_lem},  a point $x\in \q$ is generic if and only if
\[
    [\q,x] + \q^{\z_\q(x)}=\q \ .
\]
Since $\q^{\z_\q(x)}$ is the centre of the Lie algebra $\z_\q(x)$ and 
$\dim[\q,x]=\dim\q- \dim\z_\q(x)$, 
one immediately derives

\begin{s}{Proposition}   \label{sgp-ad}
An algebraic Lie algebra $\q$ has a generic centraliser if and only if there is an
$x\in\q$ such that
\begin{gather}  \label{sop-com}
\text{ $\z_\q(x)$ \ is commutative and } \\ \label{eq-summa}
\text{ $[\q,x]\oplus \z_\q(x)=\q$. }
\end{gather}
\vskip-1ex
\end{s}%
Equality~\re{eq-summa} implies that $\Ima(\ad x)=\Ima(\ad x)^2$. 
The latter is never satisfied if $\ad x$ is nilpotent and $\Ima(\ad x)\ne 0$.
That is, if $\q$ is nilpotent and $[\q,\q]\ne 0$, then $\q$ has no 
generic centralisers.
It also may happen that neither of the centralisers $\z_\q(x)$ 
is commutative.
(Consider the Heisenberg Lie algebra $H_n$ of dimension $2n+1$ for $n\ge 2$.)
On the other hand, 
if there is a semisimple $x\in \q$ such that $\z_\q(x)$ is commutative, then 
the conditions of Proposition~\ref{sgp-ad} are satisfied, so that a 
generic centraliser exists. [{\sl Warning}: this does not imply that the semisimple
elements are dense in $\q$.]
\begin{s}{Lemma}  \label{n=z}
Let $x\in \q$ be a generic point. Then $\n_\q(\z_\q(x))=\z_\q(x)$.
\end{s}\begin{proof}
Assume that $\n_\q(\z_\q(x))\ne\z_\q(x)$. In view of Eq.~\re{eq-summa}, there is 
then a nonzero 
$y\in \n_\q(\z_\q(x))\cap [\q,x] $. That is, $y=[s,x]$ for some $s\in\q$. Then
\[
[y,\z_\q(x)]=[[s,\z_\q(x)],x]\subset [\q,x]
\]
and hence $[y,\z_\q(x)]=0$. Thus, $y\in \z_\q(x)\cap [\q,x]=0$, and we are done. 
\end{proof}%
Recall that a subalgebra $\h$ of $\q$ is called a {\it Cartan subalgebra\/}
if $\h$ is nilpotent and $\n_\q(\h)=\h$. Every Lie algebra has a Cartan subalgebra,
and all Cartan subalgebras of $\q$ are conjugate under $Q$, see \cite[Ch.\,III]{serre}.

\begin{s}{Proposition}  \label{cartan}
An algebraic Lie algebra $\q$ has a generic centraliser if and only if the 
Cartan subalgebras of $\q$ are commutative. 
\end{s}\begin{proof}
If $\q$ has a generic centraliser, then, by Lemma~\ref{n=z}, such a 
centraliser is a (commutative) Cartan subalgebra. Conversely, any 
Cartan subalgebra of $\q$ is of the form 
$\h=\{y\in\q\mid (\ad x)^n y=0 \ \text{ for }\ n\gg 0\}$ 
for some $x\in\q$ \cite[Ch.\,III.4, Cor.\,2]{serre}.
Therefore, the commutativity of $\h$ implies that 
$\h=\z_\q(x)$ and $\ad x$ is invertible on $[\q,x]$.
\end{proof}%
As is already mentioned, the existence of a generic centraliser implies
that of a generic isotropy group. For this reason, we always assume 
that a generic point $x$ has the property that $Z_Q(x)$ is a generic isotropy
group. (This is only needed if a generic isotropy group
is disconnected.)

\begin{s}{Theorem}   \label{field-ad}
Suppose $\q$ has a generic centraliser. Let $x\in \q$ be a generic point
such that $Z_Q(x)$ is a generic isotropy group.
Then (i) \ ${\mathsf Z}(\z_\q(x))=Z_Q(x)$ and (ii) \ 
$\bbk(\q)^Q\simeq \bbk (\z_\q(x))^W$, where
$W=N_Q(\z_\q(x))/Z_Q(x)$ is a finite group. 
\end{s}\begin{proof}
(i) \ Since $x\in\z_\q(x)$, we have ${\mathsf Z}(\z_\q(x))
\subset Z_Q(x)$. Hence one has to prove that $Z_Q(x)$ acts trivially on
$\z_\q(x)$. Assume that the fixed point space of $Z_Q(x)$
is a proper subspace of $\z_\q(x)$, say $M$. Since
$\dim Q{\cdot}M\le \dim[\q,x]+\dim M < \dim\q$, $Q{\cdot}M$ cannot be dense
in $\q$, which contradicts the fact that $Z_Q(x)$ is a generic isotropy group.

(ii) \ This follows from Eq.~\re{sgp-field} and Lemma~\ref{n=z}.
\end{proof}%
Below, we state a property of generic points related to 
the dual space $\q^*$.

\begin{s}{Proposition} \label{2.7}
Let $x\in\q$ be a generic point, as in Theorem~\ref{field-ad}. Then

{\sf (i)} \  
$\q^*=x{\cdot}\q^*\oplus(\q^*)^x=x{\cdot}\q^*\oplus(\q^*)^{\z_\q(x)}$ and \ 
{\sf (ii)} \  
$(\q^*)^{Z_Q(x)}=(\q^*)^{\z_\q(x)}$.
%%\end{itemize}
\end{s}\begin{proof}
(i) \ We have $[\q,x]^\perp=(\q^*)^x$ and $\z_\q(x)^\perp= x{\cdot}\q^*$.
Hence the first equality follows from Eq.~\re{eq-summa}.
\\[.6ex]
The second equality means that $(\q^*)^x = (\q^*)^{\z_\q(x)}$.
Clearly, $(\q^*)^x \supset (\q^*)^{\z_\q(x)}$. Taking the annihilators provides 
the inclusion $[\q,x]\subset [\q,\z_\q(x)]$. 
Then using Eq.~\re{sop-com} and \re{eq-summa} yields

$[\q,\z_\q(x)]\subset [\z_\q(x)+[\q,x],\z_\q(x)]=[[\q,x],\z_\q(x)]=
[[\q,\z_\q(x)],x]\subset [\q,x]$.
\\[.7ex]
(ii) \ In view of (i),
%%these decompositions, 
$(\q^*)^{\z_\q(x)}$ is identified with
$(\z_\q(x))^*$. Hence the assertion stems from 
Theorem~\ref{field-ad}(i).
\end{proof}%
Thus, the very existence of a generic centraliser 
implies that $\q$ has some properties in common with
reductive Lie algebras. For instance, the Weyl group of $(Q:\q)$ is finite,
and the decomposition of $\q^*$  with respect to a generic element $x\in\q$
is very similar to that of $\q$.
It will be shown below that there is a vast stock of such Lie algebras.

                %%%%
%%%%%%%%%%%%%%%%
%%%%%%%%%%%%%%%%   Section 3
%%%%%%%%%%%%%%%%
                %%%%

\section{Generic stabilisers for the coadjoint representation}  
\label{sgp-coadj}
\setcounter{equation}{0}

\noindent
In this section, we work with the coadjoint representations of $Q$ and $\q$. 
Usually, we use lowercase Latin (resp. Greek) letters to denote elements of 
$\q$ (resp. $\q^*$). 
%%and Greek letters to denote elements of $\q^*$. 
By Lemma~\ref{ela_lem},  a point $\xi\in \q^*$ is generic if and only if
\[
    \q{\cdot}\xi + (\q^*)^{\q_\xi}=\q^* \ .
\]
As was noticed by Tauvel and Yu \cite{tayu}, taking the annihilators yields a 
simple condition, entirely in terms of $\q$.
Namely, $\xi$ is generic if and only if
\begin{equation}  \label{coad-sgp}
       \q_\xi \cap [\q,\q_\xi] =\{0\} \ .
\end{equation}
Below, we assume that $(\q, \ads)$ has a generic stabiliser and thereby
Eq.~\re{coad-sgp} is satisfied for some $\xi$.
This readily implies that $\q_\xi$ is commutative and
$\n_\q(\q_\xi)=\z_\q(\q_\xi)$. 
However, unlike the adjoint representation case, 
$\q_\xi$ can be a proper subalgebra of 
$\z_\q(\q_\xi)$. In other words, the Weyl group of
$(Q:\q^*)$ is not necessarily finite. Our goal is to understand 
what isomorphism \re{sgp-field} means in this situation.
Set $\h=\q_\xi$ and $U=(\q^*)^{\q_\xi}$. 
Then we can write
\[
   \bbk(\q^*)^Q\simeq (\bbk(U)^{Z_Q(\h)^o})^{N_Q(\h)/Z_Q(\h)^o} \ .
\]
That is, one first takes the invariants of the {\sl connected\/} group $Z_Q(\h)^o$, and
then the invariants of the {\sl finite\/} group $N_Q(\h)/Z_Q(\h)^o$.
\begin{s}{Lemma}   \label{dim-U}
$\dim U=\dim \z_\q(\h)$.
\end{s}\begin{proof}
By Lemma~\ref{inters} and Eq.~\re{coad-sgp}, we have $\q{\cdot}\xi\cap U=
\z_\q(\h){\cdot}\xi$. Equating the dimensions of these spaces 
yields the assertion.
\end{proof}%
In view of this equality, it is tempting to interpret $U$ as the space of the
coadjoint representation of $\z_\q(\h)=\Lie Z_Q(\h)^o$.
However it seems to only be possible under an additional assumption on
$\h$. 

\begin{rem}{Definition}   \label{def-near}
We say that a subalgebra $\h$ is {\it near-toral\/} if $[\q,\h]\cap\z_\q(\h)=\{0\}$.
\end{rem}%
This condition is stronger than \re{coad-sgp}. It is obviously satisfied if 
$\h$ is a toral Lie algebra (=\,Lie algebra of a torus). 
\\
Recall that the {\it index\/} of (a Lie algebra) $\q$, $\ind\q$,
is the minimal codimension of $Q$-orbits in $\q^*$. Equivalently,
$\ind\q=\trdeg \bbk(\q^*)^Q$. If $\ind\q=0$, then $\q$ is called {\it Frobenius\/}.

\begin{s}{Theorem}   \label{thm-yakoby}
Suppose the generic stabiliser $\h$ is  near-toral. Then 
\begin{itemize}
\item[\sf (i)] \ $[\q,\h]\oplus\z_\q(\h)=\q$ and $U\simeq \z_\q(\h)^*$;
\item[\sf (ii)] \ $\ind\q=\ind\z_\q(\h)=\dim\h$ and
$\h$ is the centre of $\z_\q(\h)$
\end{itemize}
\end{s}\begin{proof*} 
(i) \ It is easily seen that $[\q,\h]^\perp=(\q^*)^\h=U$. Therefore
Definition~\ref{def-near} says that $\z_\q(\h)^\perp+ U=\q^*$. From Lemma~\ref{dim-U},
it then follows  that this sum (of $\z_\q(\h)$-modules)
is direct. Hence $U\simeq \q^*/\z_\q(\h)^\perp\simeq
\z_\q(\h)^*$.

(ii) \ Since $\xi$ is generic and hence regular in $\q^*$, we have
$\ind\q=\dim\h$. 
\\
For $\nu\in U^{reg}$, we have $U\cap \h^\perp=
U\cap \q{\cdot}\nu=\z_\q(\h){\cdot}\nu$. In particular, $\dim\z_\q(\h){\cdot}\nu=
\dim U-\dim\h$. Hence almost all $Z_Q(\h)$-orbits in $U$ are of
codimension $\dim\h$. This also means that the centre of $\z_\q(\h)$ cannot
be larger than $\h$.
\end{proof*}%
\begin{s}{Corollary}    \label{cor-field}
If the generic stabiliser $\h$ is near-toral, then 
$\bbk(\q^*)^Q\simeq (\bbk(\z_\q(\h)^*)^{Z_Q(\h)^o})^{F}$,
where $F=N_Q(\h)/Z_Q(\h)^o$ is finite. That is, one first takes the 
invariants of the coadjoint representation for a smaller Lie algebra and then the 
invariants of a finite group. 
\end{s}%
Under the assumption that $\h$ is near-toral, 
$\es:=\z_\q(\h)$ has the property that $\ind\es=\dim\z(\es)$.
The following results present some properties of such algebras.

\begin{s}{Proposition}   \label{prop-l2}
Suppose $\ind\es=\dim\z(\es)$. Then \par
1. \ The closure of any regular $S$-orbit in $\es^*$ is an affine space.
\par 
2. \ If\/ $\z(\es)$ is toral, then $\es/\z(\es)$ is Frobenius.
\end{s}\begin{proof}
1. If $y\in (\es^*)^{reg}$, then $\es_y=\z(\es)$ and hence
$\es{\cdot}y=\z(\es)^\perp$. Hence all points of the orbit
$S{\cdot}y$ have one and the same tangent space.
Therefore $S{\cdot}y$ is open and dense in the affine space
$y+\z(\es)^\perp$.

2. Since $\z(\es)$ is reductive, one has a direct sum of Lie algebras
$\es=\rr\dotplus\z(\es)$, and $\ind\rr=\ind\es-\ind\z(\es)=0$.
\end{proof}%
It is not, however, always true that $\es/\z(\es)$ is Frobenius. 
For instance, the Heisenberg Lie algebra $H_n$ has one-dimensional centre
and $\ind H_n=1$. But $H_n/\z(H_n)$ is commutative, so that
$\ind (H_n/\z(H_n))=2n$.

\begin{rem}{Examples}
{\sf 1.}  Let $\be$ be a Borel subalgebra of a simple Lie algebra $\g$.
Then $(\be,\ads)$ has a generic stabiliser, which is always a toral Lie algebra,
see e.g. \cite{tayu}. If $\h$ is such a stabiliser, then by
Proposition~\ref{prop-l2}, $\z_\be(\h)/\h$ is a Frobenius Lie algebra.
It is not hard to compute this quotient for all cases in which $\h\ne 0$.

$\bullet$ \ If $\g=\sln$, then $\dim\h=\left[\frac{n-1}{2}\right]$ and 
$\z_\be(\h)/\h\simeq \be(\tri)^{[n/2]}$.

$\bullet$ \ If $\g={\frak so}_{4n+2}$, then then $\dim\h=1$ and
$\z_\be(\h)/\h\simeq \be({\frak so}_{4n})$.

$\bullet$ \ If $\g=\GR{E}{6}$, then $\dim\h=2$ and
$\z_\be(\h)/\h\simeq \be({\frak so}_{8})$.

{\sf 2.}  If $\g=\sln$ or $\spn$ and $\es$ is a seaweed subalgebra of
$\g$, then a generic stabiliser for $(\es,\ads)$
always exists, and it is a toral subalgebra~\cite{aif-rais}. 
For instance, let $\p\subset{\frak gl}_{2n}$ be a maximal parabolic 
subalgebra whose Levi part is $\gln\dotplus\gln$. Then a generic stabiliser
for $(\p,\ads)$ is $n$-dimensional and toral, and
$\z_\p(\h)/\h\simeq \be(\tri)^{n}$.

{\sf 3.} There are  non-trivial examples of Lie algebras such that
a generic stabiliser for $\ads$ exists, is near-toral, 
and equals its own centraliser,
but it is not toral. Let $e$ be a nilpotent element in $\g=\sln$ 
and $\q=\z_\g(e)$. Then a generic stabiliser for
the coadjoint representation of $\q$ exists, see \cite{ksana1}.
If $\h$ is such a stabiliser, then the description of $\h$ given in 
\cite[Theorems\,1\,\&\,5]{ksana1}
shows that $\z_\q(\h)=\h$. Hence, by Corollary~\ref{cor-field},
$\bbk(\q^*)^Q$ is the field of invariants of a finite group.   
\vskip-1.5ex
\end{rem}

                %%%%
%%%%%%%%%%%%%%%%
%%%%%%%%                           Section 4
%%%%%%%%%%%%%%%%
                %%%%

\section{Semi-direct products of Lie algebras and modules of covariants}  
\label{covariants}
\setcounter{equation}{0}

\noindent
In this section, we review some notions and results that will play the
principal role in the following exposition.

{\bf (I)} \ Recall a semi-direct product construction for Lie groups and algebras.
\\[.6ex]
Let $V$ be a $Q$-module, and hence a $\q$-module.
%%If $\q$ is a Lie algebra and $V$ a $\q$-module, 
Then $\q\times V$ has a natural structure of Lie algebra,
$V$ being an Abelian ideal in it. Explicitly, if $x,x'\in \q$ and $v,v'\in V$, then
\[
   [(x,v), (x',v')]=([x,x'], x{\cdot}v'-x'{\cdot}v) \ .
\]
This Lie algebra is denoted by $\q\ltimes V$ or $\q\oplus\eps V$.
Accordingly, an element of this algebra is denoted
by either $(x,v)$ or $x+\eps v$. Here $\eps$ is regarded as a formal symbol.
Sometimes, e.g. if $V=\q$, it is convenient to think of $\eps$
as element of the ring 
of dual numbers $\bbk[\eps]=\bbk\oplus\bbk\eps$, $\eps^2=0$.  
A connected algebraic group with Lie algebra $\q\ltimes V$ 
is identified set-theoretically with $Q\times V$, and 
we write $Q\ltimes V$ for it. The product in $Q\ltimes V$ is given by
\[
    (s,v)(s',v')= (ss', (s')^{-1}{\cdot}v+v') \ .
\]
In particular,  $(s,v)^{-1}=(s^{-1}, -s{\cdot}v)$.
The adjoint representation of $Q\ltimes V$ is given by the formula
\begin{equation}    \label{adj-QV}
(\Ad(s,v))(x',v')=(\Ad(s)x', s{\cdot}v'-x'{\cdot}v) \ ,
\end{equation}
where $v,v'\in V$, $x\in\q$, and $s\in Q$.
\\[.7ex]
Note that $V$ can be regarded as either a commutative
unipotent subgroup of $Q\ltimes V$ or a commutative nilpotent subalgebra of
$\q\ltimes V$. Referring to $V$ as subgroup of $Q\ltimes V$, 
%%(resp. subalgebra of $\q\ltimes V$), 
we write $1\ltimes V$.
%% (resp. $\eps V$). 
A semi-direct product $\q\ltimes V$ is said to be {\it reductive\/} if
$\q$ is a reductive (algebraic) Lie algebra.

{\bf (II)} \ Our second important ingredient is the notion of modules of covariants.
\\[.6ex]
Let $A$ be an algebraic group, acting on an affine variety $X$, and
$V$ an $A$-module.  
The set of all $A$-equivariant morphisms from $X$ to $V$, denoted
$\Mor_A(X,V)$, has a natural structure of $\bbk[X]^A$-module.
This $\bbk[X]^A$-module is said to be the {\it module of covariants\/} (of type $V$).
It is easily seen that $\Mor_A(X,V)$ can be identified with $(\bbk[X]\otimes V)^A$. 
For any $x\in X$, we denote by $\esi_x$ the evaluation homomorphism
$\Mor_A(X,V)\to V$, which takes $F$ to $F(x)$. Obviously, $\Ima(\esi_x)\subset
V^{A_x}$.

Assume for a while that $A=G$ is reductive. Then the algebra $\bbk[X]^G$ is finitely 
generated and
$\Mor_G(X,V)$ is a finitely generated $\bbk[X]^G$-module, see e.g. \cite[2.5]{brion},
\cite[3.12]{VP}. A review of recent results on modules of covariants in the reductive
case can be found in \cite{vdB}.
The following result is proved in \cite[Theorem\,1]{cov}.

\begin{s}{Theorem}   \label{esi-onto}
If\/ $\ov{G{\cdot}x}$ is normal and
$\codim_{\ov{G{\cdot}x}}(\ov{G{\cdot}x}\setminus
G{\cdot}x)\ge 2$, then $\Ima(\esi_x)=V^{G_x}$.
\end{s}% 
Let $\g^{reg}$ be the set of regular elements of $\g$ and $T$ a maximal torus of $G$. 
The following fundamental result is due to Kostant \cite[p.\,385]{ko63}.

\begin{s}{Theorem}   \label{kost-free}
%Let $T$ denote a maximal torus of $G$.
Let $V$ be a $G$-module. Then $\dim V^{G_x}=\dim V^T$ for
any $x\in\g^{reg}$ and $\Mor_G(\g,V)$ is a free $\bbk[\g]^G$-module of 
rank\/ $\dim V^T$.
\end{s}%
In particular,  if $V^T=0$, then there is no non-trivial
$G$-equivariant mappings from $\g$ to $V$.
These modules of covariants are graded, and the degrees of minimal generating systems 
are uniquely determined. These degrees are called the {\it generalised exponents\/} of $V$.
The multiset of generalised exponents of a $\g$-module $V$ is denoted by
$\textsl{g-exp}_\g(V)$.
Similar results hold if $\g$ is replaced 
with a "sufficiently good" $G$-module, see \cite[Ch.\,III, \S\,1]{vust}
and \cite[Prop. 4.3,\, 4.6]{gerry1}.
Namely,

\begin{s}{Theorem}   \label{cov-free}
Let $\tilde V$ be a $G$-module such that $\bbk[\tilde V]^G$ is a polynomial 
algebra and the quotient morphism $\pi: \tilde V\to \tilde V\md G$ is equidimensional.
Then $\Mor_G(\tilde V,V)$ is a free $\bbk[\tilde V]^G$-module for any $G$-module $V$.
Furthermore, if $(G:\tilde V)$ is stable, then the rank of\/ $\Mor_G(\tilde V,V)$ equals 
$\dim V^H$, where $H$ is a generic isotropy group for $(G:\tilde V)$. 
%%(which is reductive in this setting).
\end{s}%
An action $(G:V)$ is said to be {\it stable\/}, if the union of closed $G$-orbits
is dense in $V$ (see \cite[7.5]{VP} and \cite{stab} about stable actions). 
If $(G:V)$ is stable, then
a generic stabiliser is reductive and $\bbk(V)^G$ is the quotient field of $\bbk[V]^G$.

In some cases, a basis for free modules
of covariants can explicitly be indicated. For any $f\in \bbk[V]$,
the differential of $f$ can be regarded as a covector field on $V$: \ 
$v\mapsto \textsl{d}f_v\in V^*$.
Starting with $f\in \bbk[V]^G$, one obtains in this way a covariant
$\textsl{d}f \in \Mor_G(V,V^*)$. 
The following result of Thierry Vust appears in \cite[Ch.\,III, \S\,2]{vust}.
%%, see also \cite[\S\,11]{lift}.

\begin{s}{Theorem}   \label{gen-free} 
Let a $G$-module $\tilde V$ satisfy all the assumptions of 
Theorem~\ref{cov-free}. Suppose also that $N_G(H)/H$ is finite.
Let $f_1,\ldots, f_m$ be a set of  basic invariants in\/ 
$\bbk[\tilde V]^G$. Then $\Mor_G(\tilde V,\tilde V^*)$ is freely generated by 
$\text{d}f_i$, $i=1,\ldots,m$.
\end{s}%

{\bf (III)} \ Here we point out a connection between modules of covariants and
invariants of semi-direct products.
\\
%%There is an obvious method of constructing invariant 
%%polynomials using covariants. 
For $F\in\Mor_A(X,V)$, define
the polynomial  $\hat F\in \bbk[X\times V^*]^A$ by the rule
$\wF(x,\xi)=\langle F(x),\xi\rangle$, where $\langle\ ,\ \rangle: V\times V^*\to \bbk$
is the natural pairing.

%This fits in the above setting as follows.

\begin{s}{Lemma}   \label{marvel}
Consider the Lie algebra $\q\ltimes V$ and the $\bbk[\q]^Q$-module
$\Mor_Q(\q, V^*)$. Then for any $F\in \Mor_Q(\q, V^*)$, we have $\wF\in 
\bbk[\q\ltimes V]^{Q\ltimes V}$. 
\end{s}\begin{proof}
Clearly, $\wF$ is $Q$-invariant. The invariance with respect to $1\ltimes V$-action
means that
\[
   \langle F(x), v\rangle = \langle F(x), v+ x{\cdot}v'\rangle
\]
holds for any $x\in\q$ and $v,v'\in V$. To this end, we notice that
$\langle F(x), x{\cdot}v'\rangle=\langle x{\cdot}F(x), v'\rangle$, and 
$x{\cdot}F(x)=0$,  since $F:\q\to V^*$ is a $Q$-equivariant morphism.
\end{proof}%
The point is that $\wF$ turns out to be invariant
with respect to the action of the unipotent group $1\ltimes V$.

                %%%%
%%%%%%%%%%%%%%%%
%%%%%%%%%%%%%%%%   Section 5
%%%%%%%%%%%%%%%%
                %%%%

\section{Generic stabilisers and rational invariants for semi-direct products}  
\label{sgp-semidir}
\setcounter{equation}{0}

\noindent
Given $Q$ and $V$, one may ask the following questions:

\begin{itemize}
\item[\bf (Q1)] \ When does a generic centraliser for $\q\ltimes V$ exist? What are 
invariant-theoretic consequences of this?
\end{itemize}
%%\\[.7ex]
It is easily seen that the existence of a generic centraliser for
$\q$ is a necessary condition. We will therefore assume that
this is the case.

\begin{s}{Theorem}    \label{sgp-qV}
Let $x\in \q$ be a generic point. Suppose $V^x=V^{Z_Q(x)}$ and 
$V^x \oplus x{\cdot}V=V$. Then 
\begin{itemize}
\item[\sf (i)] \ each point of the form
$x+\eps v$, $v\in V^{\z_\q(x)}$, is generic and $\z_\q(x)\oplus \eps V^{\z_\q(x)}$
is a generic centraliser for $\q\ltimes V$.
\item[\sf (ii)] \ 
$\trdeg \bbk(\q\ltimes V)^{Q\ltimes V}=\trdeg\bbk(\q)^Q+\dim V^{\z_\q(x)}$;
\item[\sf (iii)] \ The Weyl groups of $(\q,\ad)$ and $(\q\ltimes V,\ad)$ are isomorphic;
\item[\sf (iv)] \ 
$\bbk(\q\ltimes V)^{Q\ltimes V}$ is a purely
transcendental extension of\/ $\bbk(\q)^Q$.
\end{itemize}
\end{s}\begin{proof}
Set $\h=\z_\q(x)$, $R=Q\ltimes V$, and $\rr=\q\ltimes V$. It follows from the 
assumptions that $ V^x=V^\h$.
\\
(i) \ Let $v\in V^\h$ be arbitrary. Let us verify that 
Proposition~\ref{sgp-ad} applies here. A direct calculation shows 
that $\z_\rr(x+\eps v)=\h\oplus\eps V^\h$ and this algebra is commutative.
Next, 
\[
[\rr, x+\eps v]=\{ [z,x]+\eps(z{\cdot}v) \mid z\in\q\}+ \eps(x{\cdot}V) \ .
\]
Notice that $\q{\cdot}v=([\q,x]\oplus\h){\cdot}v=[\q,x]{\cdot}v=x{\cdot}(\q{\cdot}v)
\subset x{\cdot}V$. Hence $z{\cdot}v\subset x{\cdot}V$ for any $z\in\q$ and
$[\rr, x+\eps v]=[\q,x]\oplus\eps(x{\cdot}V)$.
Therefore the equality $[\rr, x+\eps v]\oplus \z_\rr(x+\eps v)=\rr$ is equivalent to that
$V^x \oplus x{\cdot}V=V$.

(ii) \ By part (i), $\tilde\h:=\h\oplus \eps V^\h$ is a generic centraliser
for $\rr$. Since $\trdeg\bbk(\q)^Q=\dim\h$, the claim follows. 

(iii) \ Using formula~\re{adj-QV}, one easily verifies that 
$N_R(\tilde\h)=N_Q(\h)\ltimes V^\h$ and 
$Z_R(\tilde\h)=Z_Q(\h)\ltimes V^\h$. Hence using Theorem~\ref{field-ad}, we obtain
\[
   \tilde W=N_R(\tilde\h)/Z_R(\tilde\h)\simeq N_Q(\h)/Z_Q(\h)=W \ .
\]
\indent (iv) \ Here we may work entirely with invariants of  $W$. 
In view of (iii) and Theorem~\ref{field-ad},
it suffices to prove that $\bbk(\h)^W \hookrightarrow \bbk(\h\oplus V^\h)^W$
is a purely transcendental extension. Actually, a transcendence basis of
$\bbk(\h\oplus V^\h)^W$ over $\bbk(\h)^W$ can explicitly be constructed.
This follows from Theorem~\ref{red-inv} below, since the representation 
of $W$ on $\h$ is faithful.
\end{proof}%
The following result concerns fields of invariants of reductive algebraic groups.

Recall from Section~\ref{covariants}(III) that 
one may associate the invariant $\wF\in \bbk[V_1\times V_2]^G$
to any $F\in\Mor_G(V_1,V_2^*)$. If $D$ is a domain, then we write $D_{(0)}$ for the field
of fractions. 

\begin{s}{Theorem}   \label{red-inv}
Let\/ $\rho_i: G\to GL(V_i)$, $i=1,2$, be representations of a reductive group
$G$. Set\/ $m=\dim V_2$ and $J=\bbk[ V_1]^G$.
Suppose that %%$J_{(0)}=\bbk(V_1)^G$, 
a generic isotropy subgroup
for $(G:V_1)$ is trivial, and $(G:V_1)$ is stable. 
%%if $x\in V_1$ is generic point, then $\codim_{\ov{G{\cdot}x}}(\ov{G{\cdot}x}
%%\setminus G{\cdot}x)\ge 2$.
Then 
\begin{itemize}
\item[\sf (i)] \ $\dim_{J_{(0)}} \Mor_G( V_1, V_2^*)\otimes_J J_{(0)}=m$;
\item[\sf (ii)] \ Let $F_1,\ldots, F_m\in \Mor_G( V_1, V_2^*)$ be covariants such that 
$\{F_i\otimes 1\mid i=1,\ldots,m\}$ form a basis for the 
$J_{(0)}$-vector space in (i). Then 
$\bbk(V_1\oplus V_2)^G=\bbk(V_1)^G(\wF_1,\ldots, \wF_m)$. In other words,
any such basis for $\Mor_G(V_1, V_2^*)\otimes_J J_{(0)}$ gives rise to a
transcendence basis for the field $\bbk( V_1\oplus V_2)^G$ over $\bbk( V_1)^G$.
\end{itemize}
\end{s}\begin{proof}
(i) \ Because $(G:V_1)$ is stable, $J_{(0)}=\bbk(V_1)^G$.
Since $\Mor_G( V_1, V_2^*)$ is a finitely-generated $J$-module,
$\un{M}= \Mor_G( V_1, V_2^*)\otimes_J J_{(0)}$ is a finite-dimensional
$J_{(0)}$-vector space. By the assumptions, there is an $x\in V_1$ such
that the isotropy group $G_x$ is trivial and $G{\cdot}x=\ov{G{\cdot}x}$.
%%the boundary of $G{\cdot}x$ is "sufficiently small". 
Then by Theorem~\ref{esi-onto},
\\[.6ex]
\hbox to \textwidth{\ $(\diamond)$\hfil
the evaluation map $\esi_x: \Mor_G( V_1, V_2^*)\to  V_2^*=( V_2^*)^{G_x}$ is onto.
\hfil}
\vskip.6ex\noindent
Hence $\dim \un{M}\ge m$. On the other hand, it cannot be greater than $m$.

(ii) \ In view of Theorem~\ref{rosenlicht}, it suffices to prove that
$\bbk( V_1)^G(\wF_1,\ldots, \wF_m)$ separates the generic $G$-orbits in
$V_1\oplus V_2$. First, the field $\bbk(V_1)^G$ separates the generic $G$-orbits
in $ V_1$. Therefore, 
for generic points $(x_1,x_2)$, $x_i\in V_i$, the first coordinate is determined
uniquely up to $G$-conjugation by the values $f(x_1)$, where $f$ runs over $\bbk(V_1)^G$.
By condition $(\diamond)$, $F_1(x_1),\ldots,F_m(x_1)$ form a basis for $V_2^*$
if $x_1$ is generic.
Hence given a generic $x_1$ and arbitrary values of the invariants
$\wF_i$, the second coordinate (i.e., $x_2$) is uniquely determined.
\end{proof}%
{\bf Remarks.}
1. Most of the assumptions of Theorem~\ref{red-inv} are always satisfied if
$G$ is either finite or semisimple. For $G$ finite, it suffices to only require that
$\rho_1$ is faithful. For $G$ semisimple, it suffices to require that
a generic isotropy group of $(G:V_1)$ is trivial.

2. The assertion that the field extension in (ii) is purely transcendental
is known, see e.g. \cite[p.\,6]{dolg}. But the explicit construction of a transcendence 
basis via modules of covariants seems to be new. 
\\[.7ex]
The following assertion demonstrates important instances, where Theorem~\ref{sgp-qV}
applies.

\begin{s}{Proposition}   \label{imp-cases}
Theorem~\ref{sgp-qV} applies
to the following $\q$-modules $V$:
  
1. \ $\q$ is an arbitrary Lie algebra having a generic centraliser and
$V$ is either $\q$ or $\q^*$. 

2. \ $\q=\g$ is reductive and $V$ is an arbitrary $\g$-module.
\end{s}\begin{proof}
1. For $\q\ltimes\q$, 
the conditions of Theorem~\ref{sgp-qV}
are satisfied in view of Proposition~\ref{sgp-ad} and Theorem~\ref{field-ad}.
%% to the existence of a generic centraliser for $\q$.
For $\q\ltimes\q^*$, these conditions
%% of Theorem~\ref{sgp-qV} 
are satisfied in view of Proposition~\ref{2.7}.

2. Here $x\in\g$ is a regular semisimple element and $Z_G(x)$ is a 
maximal torus. Therefore $V^x$ is the zero weight space of $V$ and
$x{\cdot}V$ is the sum of all other weight spaces.
\end{proof}%
{\bf Remark.} For the semi-direct products as in Proposition~\ref{imp-cases}(2),
we are able to describe the polynomial invariants, see Section~\ref{sect:red}.
\\[.6ex]
A Lie algebra is said to be {\it quadratic\/} whenever
its adjoint and coadjoint representations are equivalent. 
It is easily seen that $\q\ltimes\q^*$ is quadratic for any Lie algebra
$\q$. For, if $\langle\ ,\ \rangle$ is the pairing of $\q$ and $\q^*$, then
the formula 
$(x_1+\eps \xi_1, x_2+\eps \xi_2)=\langle x_1 ,\xi_2 \rangle+
\langle x_2 ,\xi_1 \rangle$ determines a non-degenerate symmetric
$\q\ltimes\q^*$-invariant form. 
For $\q\ltimes\q^*$, there is no difference between the adjoint and
coadjoint representations.
So, previous results of this section describe some properties of the 
coadjoint representation of $\q\ltimes \q^*$ as well. 
However, for an arbitrary $V$ the adjoint and coadjoint representation of 
$\q\ltimes V$ are very different. Hence our second problem is:

\begin{itemize}
\item[\bf (Q2)] \ When does a generic stabiliser for $(\q\ltimes V,\ads)$ exist? What are 
invariant-theoretic consequences of this?
\end{itemize}
This problem is quite different from {\bf (Q1)}.
It seems to be more involved and restrictive. 

Set $\rr=\q\ltimes V$ and $R=Q\ltimes V$.
The dual space $\rr^*$ is identified with $\q^*\oplus V^*$, and
a typical element of it is denoted by $\eta=(\ap,\xi)$. For $(s,v)\in\rr$, 
the coadjoint representation is given by
\begin{equation}  \label{action*}
   (\ad_\rr^*(s,v))(\ap,\xi)=(\ad_\q^*(s)(\ap)-v{\ast}\xi,
    s{\cdot}\xi) \ . 
\end{equation}
Here the mapping $((s,\xi)\in \q\times V^*)\mapsto (s{\cdot}\xi\in V^*)$ is the 
natural $\q$-module structure on $V^*$, and 
$((v,\xi)\in V\times V^*)\mapsto (v{\ast}\xi\in \q^*)$ is the moment mapping
with respect to the symplectic structure on $V\times V^*$.
\\[.6ex]
To describe the stabiliser of any point in $\rr^*$, we need
some notation. For $\ap\in\q^*$, let
$\ck_\ap$ denote the Kirillov form on $\q$ associated with $\ap$, i.e.,
$\ck_\ap(s_1,s_2)=\langle\ap,[s_1,s_2]\rangle$. Then $\ker(\ck_\ap)=\q_\ap$,
the stabiliser of $\ap$.
If $\h$ is a subalgebra of $\q$, then $\ck_\ap\vert_\h$ can also be
regarded as the Kirillov form
associated with $\ap\vert_\h\in \h^*$.

\begin{s}{Proposition} \label{stab*}
For any $\eta=(\ap,\xi)\in\rr^*$, we have
\[
\rr_\eta= \{(s,v)\mid s\in\ker(\ck_\ap\vert_{\q_\xi}) \ \ \& \ \ 
\ad_\q^*(s)\ap=v{\ast}\xi\}\ .
\]  \vskip-1ex
\end{s}\begin{proof}
Straightforward. The first condition imposed on $s$ guarantees us the equality
$s{\cdot}\xi=0$ and that
the equation $\ad_\q^*(s)\ap=v{\ast}\xi$ has a solution $v$ for any such $s$.
\end{proof}%
It follows that $\rr_\eta$ is a direct sum of the
space $\{w\in V\mid w{\ast}\xi=0\}=(\q{\cdot}\xi)^\perp$, sitting
in $V$, and a space of dimension $\dim\ker(\ck_\ap\vert_{\q_\xi})$,
which is embedded in $\q\ltimes V$ somehow diagonally.
(We will see below that under additional constraints this second space
lies entirely in $\q$.)

A result of Ra\"\i s on semi-direct products \cite{rais} describes $\rr$-regular 
points in $\rr^*$ and
gives the value of $\ind\rr$, that is, the dimension of the 
stabilizer of the $\rr$-regular points in $\rr^*$. Namely, if $\xi\in V^*$ is 
$\q$-regular, then $(\ap,\xi)$ is $\rr$-regular if and only if
$\ap$ is $\q_{\xi}$-regular as element of $\q_{\xi}^*$ 
(with respect to the coadjoint representation
of $\q_\xi$). By a theorem of Duflo-Vergne \cite{dv}, 
the stabiliser of any regular point in the coadjoint representation is commutative, 
see also \cite[1.8]{cambr} for an invariant-theoretic proof.

It seems to be difficult to find out a general condition ensuring
that Eq.~\re{coad-sgp} holds for
some regular point in $\rr^*$. For this reason, we only look at  
the three cases occurring already in Proposition~\ref{imp-cases} in
connection with generic centralisers.

$\bullet$ \ If $(\q, \ads)$ has a generic stabiliser, then 
$(\q\ltimes\q, \ads)$ has. \\
Indeed, if $\q_\xi$ is a generic stabiliser ($\xi\in\q^*$),
then $\q_\xi\ltimes\q_\xi$ is the stabiliser of 
$\eta=(0,\xi)\in (\q\ltimes\q)^*$
and $[\q\ltimes\q, \q_\xi\ltimes\q_\xi]\cap (\q_\xi\ltimes\q_\xi)=\{0\}$.
%%On can take $\eta=(0,\xi)\in (\q\ltimes\q)^*$. Thus, $\eta$ is generic.

$\bullet$ \ If $(\q, \ads)$ has a generic stabiliser, then 
$(\q\ltimes\q^*, \ads)$ may have no generic stabilisers. \\
{\it Example}. Let $\q$ be the 3-dimensional Heisenberg algebra $H_1$.
{\it The\/} generic stabiliser for $(\q,\ads)$ exists and equals the centre of $\q$.
But $\hat\q=\q\ltimes \q^*$ is nilpotent and quadratic. Therefore $(\hat\q,\ads)
\simeq (\hat\q,\ad)$ has no generic stabiliser.

$\bullet$ \ Suppose $\q=\g$ is reductive.  Then $\g\simeq \g^*$.
%%The rest of the section is devoted to this case.  
\\[.6ex]
By \cite{Ri}, $(\g:V^*)$ always has a generic stabiliser. Assume that this stabiliser
is reductive. There is no much harm in it, since  there are finitely
many $\g$-modules whose generic stabiliser is not reductive.
Then our goal is to prove that the existence of a generic stabiliser for 
$(\rr,\ads)$ imposes a very strong constraint on the action $(G:V)$.

Let $\Omega_{V^*}$ be the open subset of $\g$-generic points in $V^*$.
Fix a generic stabiliser $\h\subset\g$ and a Cartan subalgebra $\te_\h\subset\h$.

\begin{s}{Lemma}   \label{Xi}
There is an open $R$-stable subset\/ 
$\Xi\subset (\rr^*)^{reg}\cap (\g\times \Omega_{V^*})$ 
such that if $\eta=(\ap,\xi)\in \Xi$, then $\rr_\eta$
is a direct sum of two spaces, one lying in $\g^*$ and another lying in
$V^*$. Furthermore, eventually replacing $\eta$ with an $H$-conjugate point, one can 
achieve that\/  $\rr_\eta=\te_\h\ltimes (\g{\cdot}\xi)^\perp$.
\end{s}\begin{proof}
Since $\h$ is reductive, the $\h$-modules $\h$ and $\h^*$ can be identified
using the restriction to $\h$ of a non-degenerate $\g$-invariant symmetric
bilinear form on $\g$. 
Suppose $\eta=(\ap,\xi)\in (\rr^*)^{reg}\cap (\g\times \Omega_{V^*})$. Without loss of
generality, assume that $\g_\xi=\h$.
As was explained above, the $\rr$-regularity of $\eta$ means that
$\ap$ is $\h$-regular as an element of $\h^*$. Having identified $\h^*$ and $\h$, we may
assume that $\ap$ is regular semisimple. This last condition distinguishes the 
required subset $\Xi$. Then 
$\ker(\ck_\ap\vert_\h)$ is a Cartan subalgebra of $\h$, and
if $s\in \ker(\ck_\ap\vert_\h)$, then $\ad_\g^*(s)\ap=0$.
Comparing this with Proposition~\ref{stab*}, we see that 
$\rr_\eta=\ker(\ck_\ap\vert_\h)\ltimes (\g{\cdot}\xi)^\perp$. 
Taking an $H$-conjugate, which does not affect $\xi$, we may
achieve that $\ker(\ck_\ap\vert_\h)=\te_\h$.
\end{proof}%
Thus, for almost all $\rr$-regular points in $\rr^*$, their stabilisers are
conjugate to subalgebras of the form $\tilde\h=\te_\h\ltimes (\g{\cdot}\xi)^\perp$.
Set $U=(\g{\cdot}\xi)^\perp$. By the very construction, $U$ is $\h$-stable.
Since $\eta$ is regular
and therefore $\rr_\eta$ is commutative, 
$\te_\h$ acts trivially on $U$, i.e., $\te_\h{\cdot}U=0$.

\begin{s}{Proposition}   \label{zth-polar} \\
{\sf 1}. \ Suppose $(\rr,\ads)$ has a generic 
stabiliser. Then\/  $\g{\cdot}U\cap U=\{0\}$. 
\\
{\sf 2.} \ If\/ $\h=0$, then the converse is also true.
\end{s}\begin{proof}
1. \ By Lemma~\ref{Xi} and Eq.~\re{coad-sgp}, $(\rr,\ads)$ has a generic 
stabiliser if and only if 
$[\rr,\tilde\h]\cap\tilde\h =\{0\}$. We have 
\[
[\rr,\tilde\h]=[\g\ltimes V,\te_\h\ltimes U]=[\g,\te_\h\ltimes U]+
\te_\h{\cdot}V \ .
\]
Clearly, $\g{\cdot}U$ is a subspace of $[\g,\te_\h\ltimes U]$. 
Hence we get the condition that $\g{\cdot}U\cap U=\{0\}$.
\\[.6ex]
2. Let $V=U\oplus V'$ be an $\h$-stable decomposition. Then 
$\te_\h{\cdot}V=\te_\h{\cdot}V'\subset V'$. Hence this summand causes no harm.
If $\h=0$, then $[\g,\te_\h\ltimes U]=\g{\cdot}U$. Therefore the condition 
$\g{\cdot}U\cap U=\{0\}$ appears to be necessary and sufficient for the existence of
a generic stabiliser.
\end{proof}%
Recall from \cite{polar} the notion of a polar representation of 
a reductive group.
Let $v\in V$ be semisimple, i.e., $G{\cdot}v$ is closed.
Define  $\ce_v=\{x\in V\mid \g{\cdot}x\subset \g{\cdot}v\}$. Then $(G:V)$
is said to be {\it polar\/} if there is a semisimple $v\in V$ such
that $\dim\ce_v=\dim V\md G$. Such $\ce$ is called a {\it Cartan subspace}.
%%It was shown in \cite{polar} that 
Polar representations have a number of nice (and hence restrictive) properties.
For instance, all points of $\ce$ are semisimple,
all Cartan subspaces are $G$-conjugate,
the group $W_\ce:={\mathsf N}(\ce)/{\mathsf Z}(\ce)$ is finite,
and $\bbk[V]^G\simeq \bbk[\ce]^{W_\ce}$ \cite{polar}. 
The latter implies that $\bbk[V]^G$ is polynomial and the morphism
$\pi_G: V\to V\md G$ is equidimensional \cite{isom}.

%%The action $(G:V)$ is said to be {\it stable\/} if the union of 
%%closed $G$-orbits contains a dense open subset of $V$ (see \cite[7.5]{VP} about
%%stable actions).

Our main result related to Question {\bf (Q2)} is:

\begin{s}{Theorem}   \label{main-Q2}
Suppose the action $(G{:}V)$ is stable.
Then $(\rr=\g\ltimes V, \ads)$ has a generic stabiliser if and only if
$(G{:}V)$ is a polar representation. 
\end{s}\begin{proof}
1. Suppose $(\g\ltimes V,\ads)$ has a generic stabiliser. 
\\[.6ex]
Choose $\eta=(\ap,\xi)\in \Xi$ as prescribed by Lemma~\ref{Xi}, so that
$\rr_\eta=\te_\h\ltimes U$ is a generic stabiliser and hence $\g{\cdot}U\cap U=\{0\}$
(Proposition~\ref{zth-polar}). In view of stability, 
we may also assume that $\xi$ is ($\g$-regular and) semisimple. Let us 
prove that $U$ is a Cartan subspace of $V$.

As is well known, $\dim V\md G=\dim V^*\md G$ and $(G:V)$ is stable
if and only if $(G:V^*)$ is, see e.g. \cite{stab}. By the stability hypothesis, 
\[
\max_{\nu\in V^*}\dim G{\cdot}\nu=\max_{v\in V}\dim G{\cdot}v=\dim V-\dim V\md G \ .
\]
Hence $G{\cdot}\xi=\dim V-\dim V\md G$ and
$\dim U=\dim V\md G$.
\\[1.1ex]
{\bf Claim.} \  %%\label{U-reg}
{\sl There is a closed $G$-orbit of maximal dimension meeting\/ $U$.}
%%\end{rems}%
\\[.8ex]
{\sl Proof of the Claim}. \quad 
The proof of main results in \cite{polar} is based on transcendental methods
(compact real forms of $G$, Kempf--Ness theory). This is an excuse for our
using similar methods below. In the next paragraph, we assume that $\bbk=\mathbb C$.

 Let $G_c$ be a maximal compact subgroup
of $G$ with Lie algebra $\g_c$. Fix a $G_c$-invariant Hermitian form $<\ ,\ >$
on $V^*$. Without loss of generality, we may assume that $\xi$ is of minimal length
in $G{\cdot}\xi$ and hence $<\g{\cdot}\xi,\xi >=0$, see \cite[Sect.\,1]{polar}.
Upon the identification the $\g_c$-modules $V$ and $V^*$ via
$<\ ,\ >$, $\xi$ appears to be a point of $U$. If $\tilde v\in U$  
corresponds to $\xi$ under this identification, then
we still have $<\g_c{\cdot}\tilde v,\tilde v>=0$, and therefore
$<\g{\cdot}\tilde v,\tilde v>=0$. Hence $G{\cdot}\tilde v$ is closed
\cite[Theorem\,1.1]{polar}. Since $({\g_c})_\xi=({\g_c})_{\tilde v}$ and
$({\g_c})_{\tilde v}$ is a compact real form of $\g_{\tilde v}$
\cite[Prop.\,1.3]{polar}, we conclude that
$\dim\g_{\tilde v}=\dim\g_\xi=\dim\h$. \hfill $\Box$
\\[.6ex]
%%\noindent 
For $\tilde v\in U$, we have $\dim\g{\cdot}\tilde v=\dim V-\dim U$.
Hence $\g{\cdot}\tilde v=\g{\cdot}U$ for dimension reason. In particular,
$\g{\cdot}y\subset \g{\cdot}\tilde v$ for any $y\in U$.
Thus, $U$ satisfies all conditions in the definition of a Cartan subspace.

2. Suppose $(G:V)$ is stable and polar.
\\[.6ex]
Let $v\in V$ be a regular semisimple element and
$\ce=\ce_v$ the corresponding Cartan subspace.
Then $V=\g{\cdot}\ce\oplus \ce$ and $\g{\cdot}\ce=\g{\cdot}v$ 
\cite[Section\,2]{polar}. 
Set $\h=\g_v$.
The Lie algebra $\es:=\te_\h\ltimes \ce$ is commutative, and a direct 
verification shows that it satisfies Eq.~\re{coad-sgp}. Indeed,
\[
  [\g\ltimes V,\te_\h\ltimes \ce]=[\g,\te_\h\ltimes \ce]+
\te_\h{\cdot}V \ .
\]
Using the $\te_\h$-stable decomposition $V=\g{\cdot}\ce\oplus \ce$,
we see that $\te_\h{\cdot}V\subset \g{\cdot}\ce$.
As for the first summand,  its $\g$-component does not belong to
$\te_\h$ and its $V$-component  belongs to $\g{\cdot}\ce$. Hence
$[\g\ltimes V,\te_\h\ltimes \ce]\cap (\te_\h\ltimes \ce)=\{0\}$.
It remains to find an $\eta\in\rr^*$ such that $\rr_\eta=\es$.

The dual version of the previous Claim shows that 
$(\g{\cdot}\ce)^\perp$ is a Cartan subspace of $V^*$ and that, for
sufficiently general $\xi\in (\g{\cdot}\ce)^\perp$, we have $\g_\xi=\h$
and $\g{\cdot}\xi=\ce^\perp$. Now, take an $\ap\in\g^*$ such that
under the identification $\g^*\simeq \g$ it becomes a regular element
of $\te_\h$ (i.e., $\ap\in (\te_\h)^{reg}$).
Then $\gamma=(\ap,\xi)\in (\rr^*)^{reg}$ and $\rr_\gamma=\es$.
\end{proof}%
We mention without proof the following consequence of Theorem~\ref{main-Q2}.

\begin{s}{Corollary}  \label{near}
If a generic stabiliser $\h$ for $(\rr,\ads)$ is near-toral, then 
$\rk\h=\rk\g$ and $U=V^\h$. In case of $\g$ simple, this implies that
$V$ is either the adjoint or "little adjoint" $\g$-module.
{\ququ (The latter means that the highest weight is the short dominant
root, in case $\g$ has roots of different length.)}
\end{s}%
{\bf Remark.}  It may happen that a generic stabiliser for $(G:V^*)$ is not 
reductive, but $(\g\ltimes V,\ads)$ still has a generic stabiliser. Indeed, there are 
$G$-modules $V$ such that $\rr=\g\ltimes V$ is Frobenius, i.e., $\rr^*$ has a dense 
$R$-orbit, which certainly ensures the existence of a generic stabiliser.
For $G$ simple, the list of such $V$ is obtained in \cite{ela2}.
\\[.7ex]
Finally, we consider the field of rational invariants for the coadjoint
representation of $\rr=\g\ltimes V$.
By \cite{rais}, 
\[
\trdeg\bbk(\rr^*)^R=\trdeg\bbk(V^*)^G+\ind\h \ ,
\]
where $\h$ is a generic stabiliser for $(G:V^*)$.
It follows from Eq.~\re{action*} 
that $\bbk(V^*)^G$ can be regarded as a subfield
of $\bbk(\rr^*)^{R}$. 
%%If $\ind\h=0$, then the previous equality seems to
%%suggest that these two fields can be equal. This is really the case.

\begin{s}{Theorem}  \label{ind=0}
If\/ $\ind\h=0$, then $\bbk(\rr^*)^R\simeq \bbk(V^*)^G$.
\end{s}\begin{proof}
It suffices to verify that $\bbk(V^*)^G$
separates $R$-orbits in a dense open subset of $\rr^*$. 
Let $p:\rr^*=V^*\oplus \g^*\to V^*$ denote the projection.
If $\co\subset V^*$ is a generic $G$-orbit, then we will prove
that $p^{-1}(\co)$ contains a dense $R$-orbit. The latter is equivalent
to that, for any $\xi\in \co$,
 $G_\xi\ltimes V$ has a dense orbit in $p^{-1}(\xi)=\{\xi\}\times \g^*$.
Since $1\ltimes V$
is a normal subgroup of $G_\xi\ltimes V$, we first look at its orbits.
For any $(\xi,\ap)\in p^{-1}(\xi)$, we have
$(1\ltimes V){\cdot}(\xi,\ap)=(\xi, \ap+ V\ast\xi)$. Hence all orbits are parallel
affine space of dimension $\dim(V\ast\xi)$. Therefore, it will be sufficient to prove
that $G_\xi$ has a dense orbit in the (geometric) quotient
$p^{-1}(\xi)/(1\ltimes V)$.
Because $V\ast\xi=(\g_\xi)^\perp$, that quotient is isomorphic to 
$\g^*/(\g_\xi)^\perp\simeq (\g_\xi)^*$ as $G_\xi$-variety.
Now, the presence of a dense $G_\xi$-orbit in $(\g_\xi)^*$
exactly means that $\ind\g_\xi=0$, which is true as $\xi$ is generic.
\end{proof}%
{\bf Remarks.} 1. \ In Theorem~\ref{ind=0}, the reductivity of $G$ is not needed. It suffices 
to assume that $(G:V^*)$ has a generic stabiliser.

2.\ A related result for $\bbk(\rr^*)^R$ is obtained in \cite[Corollary\,2.9]{aif-rais}
under the assumption that $\trdeg\bbk(V^*)^G=0$, but without assuming that
$G$ is reductive.

                %%%%
%%%%%%%%%%%%%%%%
%%%%%%%%%%%%%%%%   Section 6
%%%%%%%%%%%%%%%%
                %%%%

\section{Reductive semi-direct products and their polynomial invariants}  
\label{sect:red}
\setcounter{equation}{0}
\setcounter{subsubsection}{0}

\noindent
In this section, we study polynomial invariants of semi-direct products
$\q=\g\ltimes V$, where $\g$ is reductive. 

Our main technical tool is the following result of Igusa
(see \cite[Lemma\,4]{ig},\,\cite[Theorem\,4.12]{VP}). 
For reader's convenience, we provide a proof.
Given an irreducible variety $Y$, we say that an open subset
$\Omega\subset Y$ is {\it big\/} if $Y\setminus\Omega$ contains no divisors.

\begin{s}{Lemma {\ququ (Igusa)}}   \label{igusa}
Let $A$ be an algebraic group acting regularly on an irreducible affine variety 
$X$. Suppose $S$ is an integrally closed finitely generated
subalgebra of\/ $\bbk[X]^A$ and the morphism $\pi: X\to \spe S=:Y$ has the properties:

(i) \ the fibres of $\pi$ over a dense open subset of\/ $Y$ 
contain a dense $A$-orbit;

(ii)  \ $\Ima\pi$ contains a big open subset of\/ $Y$.
\\
Then $S=\bbk[X]^A$. In particular, the algebra of $A$-invariants is finitely generated.
\end{s}\begin{proof}
From (i) and Rosenlicht's theorem, it follows that $\bbk(Y)=\bbk(X)^A$.
In particular, $\bbk(X)^A$ is the quotient field of $\bbk[X]^A$.
Assume that $S\ne \bbk[X]^A$. Then one can find a finitely generated
intermediate subalgebra: \ \ $S\subset \tilde S\subset \bbk[X]^A$ such that
$S\ne \tilde S$. The natural morphism $\tilde\pi:\spe\tilde S\to Y$ is
birational and its image contains a big open subset of $Y$ 
(because $\pi$ does). Since $Y$ is normal, the Richardson lemma \cite[3.2 Lemme\,1]{brion}
implies that $\tilde \pi$ is an isomorphism. This contradiction
shows that $S =\bbk[X]^A$.
\end{proof}%
Recall that $Q:=G\ltimes V$ is a connected group with Lie algebra $\q$.
Here $1\ltimes V$ is exactly the unipotent radical of $Q$, which is also denoted
$Q^{u}$. Let $T$ be a maximal torus of $G$ with the corresponding Cartan subalgebra
$\te$.
 
First, we consider the adjoint representation of $\g\ltimes V$.

\begin{s}{Theorem}   \label{main-reduct}
Let\/ $V$ be an arbitrary $G$-module, $\q=\g\ltimes V$, and $m=\dim V^\te$. Then
\begin{itemize}
\item[\sf (i)] \  $\bbk[\q]^{Q^u}$ is a polynomial algebra of Krull
dimension  $\dim\g+m$.
It is freely generated by the coordinates on $\g$ and the functions
$\wF_i$, $i=1,\ldots,m$, associated with covariants of type $V^*$ (see below).
\item[\sf (ii)] \  $\bbk[\q]^{Q}$ is a polynomial algebra of
Krull dimension $\dim\te+m$.
It is freely generated by the basic invariants of\/ $\bbk[\g]^G$ and the same
functions $\wF_i$, $i=1,\ldots,m$. 
\item[\sf (iii)] \  $\max\dim_{x\in \q} Q{\cdot}x=\dim\q-\dim \q\md Q$;
\item[\sf (iv)] \ If\/ $\pi:\q\to \q\md Q$ is the quotient morphism, then
$\Omega:=\{x\in\q \mid \text{d}\pi_x\ \text{ is onto }\}$ is a big open subset of $\q$. 
\end{itemize}
\end{s}\begin{proof*}
(i) \ By Theorem~\ref{kost-free}, $\Mor_G(\g, V^*)$ is a free
$\bbk[\g]^G$-module of rank $m$. Let $F_1,\ldots,F_m$ be a basis for
this module and $\wF_1,\ldots,\wF_m$ the corresponding $Q$-invariants on $\q$, i.e.,
${\wF_i(x+\eps v)}=\langle F_i(x), v\rangle$. 
To prove that $\bbk[\q]^{Q^u}$ is freely generated by the coordinate 
functions on $\g$ and the polynomials
$\wF_i$, $i=1,\ldots,m$, we wish to apply Lemma~\ref{igusa}. 

Set $\fX_m=\{x\in\g\mid \dim \text{span}\{F_1(x),\ldots,F_m(x)\}=m\}$. 
That is, $\fX_m$
is the set of those $x$, where the vectors $F_i(x)\in V^*$,
$i=1,\ldots,m$, are linearly independent.
\\[1.1ex]
{\bf Claim}. \   %% \label{omega-big}
{\sl $\fX_m$ is a big open subset of $\g$. More precisely,  
$\codim_\g (\g\setminus\fX_m)\ge 3$.}
\\[.8ex]
{\sl Proof of the claim.} \quad
The set of regular elements of $\g$, $\g^{reg}$, has the
property that $\codim (\g\setminus \g^{reg})\ge 3$ and
$\ov{G{\cdot}x}$ is normal for any $x\in \g^{reg}$ \cite{ko63}.
The condition that 
$\codim_{\ov{G{\cdot}x}} (\ov{G{\cdot}x}\setminus G{\cdot}x)\ge 2$
is satisfied for every $x\in\g$, since any $G$-orbit is even-dimensional.
By Theorems~\ref{esi-onto} and \ref{kost-free}, we conclude that $\fX_m\supset \g^{reg}$,
and the claim follows.
\\[.6ex]
Let $x_1,\ldots, x_n$ be the coordinates on $\g$, where $n=\dim\g$.
Then $x_1,\ldots, x_n,\wF_1,\ldots,\wF_m$ are algebraically independent, because their
differentials are linearly independent on $\fX_m\ltimes V$. 
Consider the mapping 
\[
   \tau: \q\to \spe\bbk[x_1,\ldots, x_n,\wF_1,\ldots,\wF_m]= \bbk^{n+m},
\]
where $\tau(x+\eps v)=(x,\wF_1(x+\eps v),\ldots,\wF_m(x+\eps v))$.
We identify $\bbk^{n+m}$ with $\g\times \bbk^{m}$.
If $x=(x_1,\ldots,x_n)\in\fX_m$, then the $F_i(x)$'s are linearly independent, so that
the system 
\[
    \wF_i(x+\eps v)=\langle F_i(x), v\rangle=\ap_i,\quad i=1,\ldots,m
\] 
has a solution $v$ for any $m$-tuple
$\ap=(\ap_1,\ldots,\ap_m)$. Hence $\Ima\tau\supset \fX_m\times \bbk^{m}$,
which means that $\Ima\tau$ contains a big open subset of $\bbk^{n+m}$.

It follows from the above Claim that 
$\g^{reg}=\fX_m\cap\{y\in \g\mid \dim G{\cdot}y=n-m\}$. 
%%It is still a non-empty open  $G$-stable subset of $\g$. Actually. $\Psi=$.
Take $x\in\g^{reg}$, and let $v_\ap$ be a solution to 
the system
$\wF_i(x+\eps v)= \ap_i$. Then 
$\tau^{-1}(x,\ap)\ni x+\eps v_\ap$ and
\[
   \tau^{-1}(x,\ap)\supset Q^u{\cdot}(x+\eps v_\ap)=\{ x+\eps (v_\ap+x\ast V)\}  \ .
\]
Since $x\in\fX_m$, we have $\dim \tau^{-1}(x,\ap)=n-m$. On the other hand,
$\dim[\g,x]=n-m$, by the definition of $\g^{reg}$.
Hence 
$\tau^{-1}(x,\ap)= Q^u{\cdot}(x+\eps v_\ap)$ for dimension reason.
Thus, a generic fibre of $\tau$ is a $Q^u$-orbit, 
and Lemma~\ref{igusa} applies here.

(ii) \ Clearly, 
\[
\bbk[\q]^{Q}=(\bbk[\q]^{Q^u})^G=\bbk[x_1,\ldots, x_n,\wF_1,\ldots,\wF_m]^G \ .
\]
Since the $\wF_i$'s are already $G$-invariant, the algebra in question is equal to
\[
    \bbk[\g]^G[\wF_1,\ldots,\wF_m] \ .
\]
\indent (iii) \ The dimension of a $Q$-orbit cannot be greater than
$\dim\q-\dim\q\md Q$, and 
if $x\in\te$ is regular, then $\dim Q{\cdot}(x+\eps 0)=
\dim Q- \dim\te-m$.

(iv) \ It follows from the previous discussion that $\Omega\supset\g^{reg}\ltimes V$.
\end{proof*}%
\begin{rem}{Remarks} 1. If $V^T=\{0\}$, then the module of covariants
of type $V^*$ is trivial, so that we obtain a natural isomorphism
$\bbk[\q]^Q\simeq\bbk[\g]^G$.
\\[.6ex]
2. From Theorem~\ref{sgp-qV} and Proposition~\ref{stab*}, it follows that 
$\tilde\te:=\te\ltimes V^\te$ is a generic centraliser in $\q$ and 
$\tilde W=N_Q(\tilde\te)/Z_Q(\tilde\te)$ is isomorphic to
$W=N_G(\te)/Z_G(\te)$, the usual Weyl group of $\g$. Therefore
\[
  \bbk(\q)^Q \simeq \bbk(\te\ltimes V^\te)^W= \bbk(\te\times V^\te)^W \ .
\]
Since $\bbk(\te)^W$ is a rational field, Theorem~\ref{sgp-qV}(iv) implies that
$\bbk(\q)^Q$ is rational, too. For $\g$ semisimple, the rationality of
$\bbk(\q)^Q$ also follows from Theorem~\ref{main-reduct}, because in this situation
$\bbk(\q)^Q$ is the quotient field of $\bbk[\q]^Q$. However, if
$V^\te\ne 0$, then the restriction homomorphism
\[
  res: \bbk[\q]^Q \to \bbk[\te\times V^\te]^W
\]
is not onto. For, the description of the generators of 
$\bbk[\q]^Q$ shows that $\bbk[V^\te]^W$ does not belong to the image
of $res$.
\end{rem}%
Now, we look at polynomial invariants of the coadjoint representation of 
$\q=\g\ltimes V$. 
As we know from Section~\ref{sgp-semidir}, the existence of a generic
stabiliser for $(\q,\ads)$ is a rare phenomenon; but this existence is
not always needed for describing  invariants.
It follows from Eq.~\re{action*} 
that $\bbk[V^*]^G$ can be regarded as a subalgebra
of $\bbk[\q^*]^{Q}$.
Recall that 
$\trdeg\bbk(\q^*)^Q=\trdeg\bbk(V^*)^G+\ind\h$,
where $\h$ is a generic stabiliser for $(G:V^*)$.
In particular, if $\g$ is semisimple and $\h$ is reductive, then
$\trdeg\bbk(\q^*)^Q=\trdeg\bbk(V^*)^G+\rk\h$.
Since the roles of $V$ and $\g$ are interchanged in the
dual space, one might hope that 
$\bbk[\q^*]^{Q}$ could be generated
by $\bbk[V^*]^G$ and certain invariants arising from $\Mor_G(V,\g^*)$.
This is however false, because it can happen that $\rk\h >0$,
but $\Mor_G(V,\g^*)=0$. In general, it is not clear  how to discover
"missing" invariants associated with the summand $\ind\h$ (or $\rk\h$). 
\\
The simplest case is that in which $\h=0$.
Then we are in a position to state an analogue of Theorem~\ref{main-reduct}.

\begin{s}{Theorem}   \label{main-dual}
As above, let $\q=\g\times V$ and $Q^u=1\ltimes V$.
Suppose a generic stabiliser for $(G:V^*)$ is trivial. Then
$\bbk[\q^*]^{Q^u}=\bbk[V^*]$ and $\bbk[\q^*]^Q=\bbk[V^*]^G$.
\end{s}\begin{proof}
The second equality stems from the first. To prove the first equality,
we use the same method as in Theorem~\ref{main-reduct}.
The natural projection $\q^*\to \q^*/\g^*\simeq V^*$ is
$Q^u$-equivariant and satisfies all the requirements of Lemma~\ref{igusa}.
The details are left to the reader.
\end{proof}%
{\bf Remark.} In Theorem~\ref{main-dual}, the reductivity of $G$ is not needed.

                %%%%
%%%%%%%%%%%%%%%%
%%%%%%%%%%%%%%%%   Section 7
%%%%%%%%%%%%%%%%
                %%%%

\section{Takiff Lie algebras and their invariants}  
\label{takif-2}
\setcounter{equation}{0}

\noindent
For $\g$ semisimple, some interesting results on the invariants of $(\g\ltimes\g, \ad)$ are 
obtained by Takiff in \cite{takiff}. For this reason, Lie algebras of the form 
$\q\ltimes\q$ are sometimes called {\it Takiff (Lie) algebras\/}, see \cite{rt},\cite{geof1}. 
We will follow this terminology.
%%Takiff Lie algebras can be obtained via the following procedure.
%%Consider $\q[T]=\q\otimes \bbk[T]$ with the natural structure of Lie algebra, where
%%$T$ is an indeterminate. Then $\q\ltimes\q$ is the quotient of $\q[T]$ by the ideal
%%generated by $1\otimes T^2$. Letting $\eps=T \mod (1\otimes T^2)$,
%%%%Below, we recall and generalise Takiff's results.

In this section, we consider orbits and invariants of certain
representations of a Takiff group $\hat Q=Q\ltimes \q$.
%% acting on $\q\ltimes\q$.
Some results on rational invariants have already appeared in Section~\ref{sgp-semidir}.
Our main object here is the polynomial (regular) invariants.
We obtain a generalisation of the main result in
\cite{takiff}, which concerns several aspects. First, in place of semisimple
Lie algebras, we consider a wider class. Second, the initial representation of $Q$
is not necessarily adjoint. Third, we also describe the invariants of
the unipotent group $1\ltimes\q\subset \hat Q$. 
%%acting on $\q\ltimes\q$. 
Fourth, our proof does not exploit complex numbers and complex topology.

If $V$ is a $\q$-module, then $V\times V$ can regarded as
$\q\ltimes\q$-module in a very natural way. For $(x_1,x_2)\in \q\ltimes\q$
and $(v_1,v_2)\in V\times V$, we define
\[
     (x_1,x_2){\cdot}(v_1,v_2):=(x_1{\cdot}v_1, x_1{\cdot}v_2-x_2{\cdot}v_1) \ .
\]
This $\q$-module will be denoted by $\hat V=V\ltimes V$.
We also write $v_1{+}\eps v_2$ for $(v_1,v_2)$.
%%As above, $N:=1\ltimes\q$ is a commutative unipotent normal subgroup of $Q$.
If $f\in \bbk[V]^Q$, then $\textsl{d}f\in \Mor_Q(V,V^*)$, and 
we define $\wF_f\in \bbk[\hat V]$ by the rule:
$\wF_f(x+\eps y)=\langle \textsl{d}f_{x}, y\rangle$.
Similarly to Lemma~\ref{marvel}, one proves 
\[
   \wF_f\in \bbk[\hat V]^{\hat Q} \ .
\]
Here one needs the fact that $\textsl{d}f_v$ annihilates the tangent space of
$Q{\cdot}v$ at $v\in V$.

\begin{s}{Theorem}   \label{main-twice}  
Let $V$ be a $Q$-module.
Suppose the action\/ $(Q:V)$ satisfies the following conditions:

(1) \ $\bbk[V]^Q$ is a polynomial algebra;

(2) \ $\max\dim_{v\in V} Q{\cdot}v=\dim V-\dim V\md Q$;

(3) \ If\/ $\pi_Q: V\to V\md Q$ is the quotient morphism and\/
$\Omega:=\{v\in V \mid (\text{d}\pi_Q)_v\ \text{ is onto }\}$, then\/ 
$V\setminus\Omega$ contains no divisors.
\\
%%Set\/ $\hat\q=\q\ltimes\q$ and $\hat Q=Q\ltimes\q$. Notice that $1\ltimes\q$
%%is a commutative unipotent normal subgroup of\/ $\hat Q$.
Then
\begin{itemize}
\item[\sf (i)] \  $\bbk[\hat V]^{1\ltimes\q}$ is a polynomial algebra of Krull
dimension\/ $\dim V+\dim V\md Q$, which is generated by the coordinates on the first
factor of $\hat V$ and the polynomials $\wF_1,\ldots,\wF_m$ associated with
the differentials of basic invariants in $\bbk[V]^Q$;
\item[\sf (ii)] \ $\bbk[\hat V]^{\hat Q}$ is a polynomial algebra of Krull
dimension\/ 2$\dim V\md Q$, which is freely generated by the basic invariants 
of\/ $\bbk[V]^Q$ and the same functions $\wF_i$, $i=1,\ldots,m$.
\item[\sf (iii)] \ The $\hat Q$-module $\hat V$ satisfies conditions (1)--(3), too.
\end{itemize}
\end{s}\begin{proof*}
The proof is very close in the spirit to the proof of Theorem~\ref{main-reduct},
though some technical details are different.

Set $N=1\ltimes \q$. 
%%It is a commutative unipotent subgroup of $\hat Q$.
Let $f_1,\ldots,f_m$, $m=\dim V\md Q$, be algebraically independent generators of
$\bbk[ V]^Q$. 
%%Then each $\textsl{d}f_i\in\Mor_Q( V, V^*)$, and we define the
As was noticed above, to each $f_i$ one may associate
the polynomial $\wF_i=\wF_{f_i}\in \bbk[ V\ltimes V]^{\hat Q}$.
%% by the rule $\wF_i(x+\eps y)
%%=\langle (\textsl{d}f_i)_{x}, y\rangle$. By Lemma~\ref{marvel}, 
%%$\wF_i\in\bbk[\hat\q]^{\hat Q}$.

(i) \ We are going to prove, using Lemma~\ref{igusa}, that $\bbk[\hat\q]^N$ is freely 
generated by the coordinate functions
on $V$ (which is the first component of $\hat V$) and the polynomials
$\wF_i$, $i=1,\ldots,m$.  
Let $x_1,\ldots, x_n$ be the coordinate functions on $V$.  
Then $x_1,\ldots, x_n,\wF_1,\ldots,\wF_m$ are algebraically independent, because their
differentials are linearly independent on $\Omega\ltimes V$. 
Consider the mapping 
\[
   \hat\tau: \hat V\to \spe\bbk[x_1,\ldots, x_n,\wF_1,\ldots,\wF_m]= \bbk^{n+m}\ .
\]
We identify $\bbk^{n+m}$ with $V\times \bbk^{m}$.
If $x=(x_1,\ldots,x_n)\in\Omega$, then $(\textsl{d}f_i)_x$ are linearly independent, so that
the system $\wF_i(x+\eps y)=\ap_i$, $i=1,\ldots,m$, 
has a solution $y$ for any $m$-tuple
$\ap=(\ap_1,\ldots,\ap_m)$. Hence $\Ima\hat\tau\supset \Omega\times \bbk^{m}$,
which means that $\Ima\hat\tau$ contains a big open subset of $\bbk^{n+m}$.

Next, consider $\Omega'=\Omega\cap\{y\in V\mid \dim Q{\cdot}y=n-m\}$. In view of
condition~(2), it is still a non-empty open 
$Q$-stable subset of $V$. Take $x\in\Omega'$,  and let $y_\ap$ be a solution to the system
$\wF_i(x+\eps y)=\langle(\textsl{d}f_i)_x,y\rangle= \ap_i$, $i=1,\ldots,m$. 
Then 
$\hat\tau^{-1}(x,\ap)\ni x+\eps y_\ap$ and
\[
   \hat\tau^{-1}(x,\ap)\supset N{\cdot}(x+\eps y_\ap)=\{ x+\eps (y_\ap+\q{\cdot}x)\}  \ .
\]
Since $x\in\Omega$, we have $\dim \hat\tau^{-1}(x,\ap)=n-m$. On the other hand,
$\dim Q{\cdot}x=n-m$, because of the definition of $\Omega'$.
Hence 
$\hat\tau^{-1}(x,\ap)= N{\cdot}(x+\eps y_\ap)$ for dimension reason.
Thus, a generic fibre of $\hat\tau$ is an $N$-orbit, and Lemma~\ref{igusa} applies here.

(ii) \ Clearly, 
\[
\bbk[\hat V]^{\hat Q}=(\bbk[\hat V]^N)^Q=\bbk[x_1,\ldots, x_n,\wF_1,\ldots,\wF_m]^Q \ .
\]
Since the $\wF_i$'s are already $Q$-invariant, the algebra in question is equal to
\[
    \bbk[V]^Q[\wF_1,\ldots,\wF_m]=\bbk[f_1,\ldots,f_m,\wF_1,\ldots,\wF_m] \ .
\]
(iii) \ We have to check that the $\hat Q$-module $\hat V$ satisfies 
properties (1)--(3). 

$\bullet$ \quad Property~(1) is verified in (ii).

$\bullet$ \quad If $x\in\Omega$, then $\dim \hat Q{\cdot}(x+\eps 0)=2n-2m$, which gives
property~(2) for $\hat Q$.

$\bullet$ \quad Set $\hat\Omega=\Omega\times V$. It is a big open subset of
$\hat V$. Explicit expressions for algebraically independent generators of
$\bbk[\hat V]^{\hat Q}$ show that their differentials are linearly independent
on $\hat\Omega$, which is exactly Property~(3) for $\hat V$.
\end{proof*}%
\begin{rem}{Remarks}
1. \ If the pair $(\q, V)$ satisfies properties~(1)-(3) of Theorem~\ref{main-twice},
then the passage $(\q, V) \mapsto (\hat\q,\hat V)$ can be iterated
{\sf ad infinitum\/} without losing those properties. 
%%In particular, we always obtain a Lie algebra
%%with a polynomial ring of invariants for the adjoint representations.

2. \ Since the adjoint representation of a semisimple Lie algebra $\g$
has properties (1)--(3),
iterating the Takiffisation procedure $\g\mapsto \g\ltimes \g$ always yields algebras 
with a polynomial ring of invariants for the adjoint representation. 
This is the main result of \cite{takiff}.
Explicit form of the basic invariants for $(\g\ltimes\g,\ad)$
is also pointed out there.  Notice also that Takiff's results follow from either Theorem~\ref{main-reduct} with $V=\g$ or
Theorem~\ref{main-twice} with $\q=V=\g$.

3. \  By Theorem~\ref{main-reduct}, the adjoint
representation of $\q=\g\ltimes V$
satisfies all the conditions of Theorem~\ref{main-twice}. Therefore these $\q$ can be used
as building blocks for Takiffisation procedure, which yields more and
more complicated Lie algebras having polynomial algebras of invariants.
\end{rem}%
Let us make some comments on the conditions of Theorem~\ref{main-twice}.  
If $Q=G$ is semisimple, then conditions (2) and (3) are always satisfied, 
regardless of the fact whether $\bbk[V]^G$ is polynomial.
For condition~(3) we refer to \cite[Satz\,2]{knop}, while (2)
follows since $G$ has no rational characters and
therefore $\bbk(V)^G$ is the quotient field of $\bbk[V]^G$. 
Thus, we have

\begin{s}{Corollary}
If $\rho:G\to GL(V)$ is a representation of a semisimple group such
that $\bbk[V]^G$ is polynomial, then Theorem~\ref{main-twice} applies to
the $G\ltimes \g$-module $V\ltimes V$.
\end{s}%
%
\begin{comment}
%%%%%%%%%%%%%%%%%%%%%%%%%%
Comparing the proofs of Theorems~\ref{main-reduct}
and \ref{main-twice}, one may notice that a basis for
the $\bbk[\g]^G$-module $\Mor_G(\g,V^*)$ was used for producing missing
invariants in the former; while the differentials
of the basic invariants in $\bbk[\q]^Q$ were used in the latter.
for the same purpose.
Since $\textsl{d}f_i\in \Mor_Q(\q,\q^*)$, 
one may ask 

{\it Is it true that $\Mor_Q(\q,\q^*)$ is a free $\bbk[\q]^Q$-module
and $(\textsl{d}f_1,\ldots,\textsl{d}f_m)$ is a basis for it?}
%%%%%%%%%%%%%%%%%%%%%%%%%%%%%%%
\end{comment}

                %%%%
%%%%%%%%%%%%%%%%
%%%%%%%%%%%%%%%%   Section 8
%%%%%%%%%%%%%%%%
                %%%%

\section{The null-cone and its irreducibility}  
\label{null-cone}
\setcounter{equation}{0}

\noindent
In previous sections, we described several instances of representations of
nonreductive Lie algebras having a polynomial algebra of invariants. 
If $  Q\subset GL(\tilde V)$ and $\bbk[\tilde V]^{Q}$ is polynomial,
then it is natural to inquire of
whether it is true that $\bbk[\tilde V]$ is a free $\bbk[\tilde  V]^{Q}$-module. 
As is well known, the freeness is equivalent to that 
the quotient morphism $\pi: \tilde V\to \tilde V\md Q$ is
{\it equidimensional\/}, i.e., 
has the property that $\dim\pi^{-1}(\pi(0))=\dim \tilde V-\dim \tilde V\md Q$.
As in the case of reductive group actions, we say that $\pi^{-1}(\pi(0))$
is the {\it null-cone}, denoted ${\goth N}^{Q}(\tilde V)$
or ${\goth N}(\tilde V)$.

In this section, we only deal with reductive semi-direct products and their
representations. Our goal is to describe necessary and
sufficient conditions for equidimensionality
of $\pi$ and point out some consequences of it.
We consider two types of representations:

{\bf A)} \ $ \q=\g\ltimes V$,  where $V$ is a $\g$-module, and  
$\tilde V=\q$, i.e., we consider the adjoint representation of
$\q$.

{\bf B)} \ $\q=\g\ltimes\g$ is a reductive Takiff algebra and 
$\tilde V=V\ltimes V$, where $V$ is a $\g$-module.

We begin with case {\bf A)}. 
Recall that $m=\dim V^T$ and 
$F_1,\ldots, F_m$ is a basis for the $\bbk[\g]^G$-module 
$\Mor_G(\g,V^*)$. The null-cone for $(\g,\ad)$ is denoted by $\N(\g)$ or merely by $\N$. 
In other words, $\N$ is the set of nilpotent
elements of $\g$. Recall that $\N$ is irreducible and
$\dim\N=\dim\g-\dim\g\md G=\dim\g-\dim\te$.
\\[.6ex]
Theorem~\ref{main-reduct} says that
if $V^T=0$, then $\bbk[\q]^Q=\bbk[\g]^G$ and therefore
${\goth N}(\q)\simeq \N\times V$. In this trivial case, $\pi_Q$ is equidimensional,
since it is so for $\pi_G: \g\to \g\md G$. Therefore we assume below that $V^T\ne 0$.

Define a stratification of $\g$ in the following way:
\[
   \fX_{i,V}={\fX}_i =\{ x\in\g\mid \dim \text{span}\{F_1(x),\ldots,
F_m(x)\}=i\} \ .
\]
Then $\ov{{\fX}}_i\subset \ov{{\fX}}_{i+1}$ and $\ov{{\fX}}_m=\g$.
The induced stratification on the null-cone is
$\fX_i(\N):=\fX_i\cap\N$. 
As is shown in the proof of Theorem~\ref{main-reduct}, 
$\fX_m$ is a big open subset of $\g$
containing $\g^{reg}$. Therefore $\fX_m(\N)$ is a big open subset
of $\N$ containing the regular nilpotent orbit.

\begin{s}{Theorem}  \label{cond-EQ}
\begin{itemize}
\item[\sf 1.] \ The quotient morphism $\pi_{  Q}:  \q \to  \q\md   Q$
is equidimensional if and only if \\
 $\codim_\N \fX_i(\N) \ge m-i$.
\item[\sf 2.] \ If\/ ${\goth N}(  \q)$ is irreducible, then $\pi_{  Q}$
is equidimensional;
\item[\sf 3.] \  ${\goth N}(  \q)$ is irreducible if and only if\/ 
$\codim_\N \fX_i(\N)\ge m-i+1$ for $i<m$.
\end{itemize}
\end{s}\begin{proof}
1. Since $\pi_Q$ is dominant, all irreducible components of ${\goth N}(\q)$ are of
dimension $\ge \dim\q-\dim\q\md Q$.
By Theorem~\ref{main-reduct}, ${\goth N}(\q)=\{(x,v) \mid
x\in\N \ \ \& \ \ \wF_i(x,v)=\langle F_i(x),v\rangle=0 \ \ \forall i\}$. 
Let $p: {\goth N}(  \q)
\to \N$ be the projection onto the first factor. Then 
${\goth N}( \q)=\displaystyle \bigsqcup_{i=0}^m p^{-1}(\fX_i(\N))$ and
$\dim p^{-1}(\fX_i(\N))=\dim \fX_i(\N)+\dim V-i$.

2. By Theorem~\ref{main-reduct}, if $e\in\N^{reg}$, then $(e,0)\in \q^{reg}$ and
$(\textsl{d}\pi_{Q})_{(e,0)}$ is onto . Therefore,
$(e,0)$ is a smooth point of ${\goth N}(\q)$, and the
unique irreducible component of ${\goth N}(\q)$
to which $(e,0)$ belongs is of dimension $\dim \q-\dim \q\md Q$.
On the other hand, $\dim p^{-1}(\fX_m(\N))=\dim \N +\dim V-m=\dim \q-
\dim  \q\md   Q$. Hence $\ov{p^{-1}(\fX_m(\N))}$ is {\sl the\/} irreducible 
component of ${\goth N}(\q)$ containing $(e,0)$.

3. The proof of part~2 shows that $\ov{p^{-1}(\fX_m(\N))}$ is an irreducible component
of ${\goth N}(\q)$  of expected dimension. To ensure the irreducibility, we have to require that
$\ov{p^{-1}(\fX_i(\N))}$ cannot be an irreducible component for $i<m$.
Since all irreducible components of ${\goth N}(\q)$ are of dimension
$\ge \dim \q-\dim  \q\md   Q$, the condition that
$\dim p^{-1}(\fX_i(\N)) < \dim \q-\dim  \q\md   Q$ for $i< m$
is equivalent to the irreducibility.
\end{proof}%
The following is now immediate.

\begin{s}{Corollary}  \label{m12}
If\/ $m=1$, then ${\goth N}(\q)$ is irreducible;
if\/ $m=2$, then $\pi_Q$ is equidimensional.
\end{s}%
\vskip-1ex
\begin{rem}{Remarks}  \label{remarks-8}
1. Since $\N$ consists of finitely many $G$-orbits, condition~\ref{cond-EQ}(1)
is equivalent to the following: if $G{\cdot}x \subset \fX_i(\N)$, then
$\dim G{\cdot}x\le \dim\N-(m-i)$, or
\[
\dim\z_\g(x)-\rk\g\ge m- \dim \bigl(\text{span}\{F_1(x),\ldots,F_m(x)\}\bigr) \ .
\]
Furthermore, a more careful look at the projection ${\goth N}(\q)\to \N$ shows that
if last condition is satisfied, then the number of the irreducible components of
${\goth N}(\q)$ equals the number of the $G$-orbits $G{\cdot}x\subset\N$ such that 
$\dim\z_\g(x)-\rk\g = m- \dim (\text{span}\{F_1(x),\ldots,F_m(x)\})$.
\\   \indent
2. The condition in Theorem~\ref{cond-EQ}(1) for $i=0$ reads
$\dim\N-\dim (\fX_0(\N))\ge m$, or $\dim V^T\le \dim\N-\dim \fX_0(\N)\le
\dim\N$. This is a rough necessary condition for $\pi_Q$ to be
equidimensional.
Let $G$ be simple and $V_\lb$ a simple $G$-module with highest weight $\lb$. 
Then $(V_\lb)^T\ne 0$ if and only if $\lb$ lies in the root lattice,
$\mathcal R$. The function
$n\mapsto \dim(V_{n\lb})^T$, $\lb\in\mathcal R$, has a polynomial growth.
The only case in which this function is constant is that of
$G=SL_p$, $\lb=p\vp_1$ or $p\vp_{p-1}$. Here $\vp_i$'s are fundamental weights,
and $\dim (V_{n\lb})^T=1$ for any $n\in \BN$. Thus, modulo this exception, 
there are finitely many simple $G$-modules $V$ such that $V^T\ne 0$ and
$\pi_Q$ is equidimensional.
\end{rem}%
For future use, we record a relationship between the stratifications of
$\N$ and $\g$.

\begin{s}{Proposition}      \label{F2}
If\/ $\codim_\N \fX_i(\N)\ge m-i+1$ for $i<m$, then 
$\codim_\g \fX_i \ge m-i+2$. 
\end{s}\begin{proof}
It follows from the definitions that $\ov{\fX_i}\subset\underset{j\le i}{\cup}\fX_j$
and $\ov{\fX_i(\N)}\subset\underset{j\le i}{\cup}\fX_j(\N)$. 
For $i< m$, we have $\fX_i \cap\g^{reg}=\varnothing$ and hence 
$\pi_G(\ov{\fX_i})=\ov{\fX_i}\md G$ is a proper subvariety of $\g\md G$.
Therefore
\\[.7ex]
  $\dim \ov{\fX_i}\le \dim\g\md G-1+\dim \ov{\fX_i(\N)}\le \dim\g-(m-i+2)$.
\end{proof}%
There is another interesting cone related
to $\q=\g\ltimes V$. 
Consider the morphism $\bar\pi:\q\to \bbk^m$,
$(x,v) \mapsto (\wF_1(x,v),\ldots,\wF_m(x,v))$.
The zero-fibre of $\bar\pi$
is denoted by ${\goth N}^u(\q)$. Thus,
\[
  {\goth N}^u(\q)=\{(x,v)\in\q \mid \langle F_i(x),v\rangle=0 \ \ i=1,\ldots,m\} \ .
\] 
The proof of the following result is entirely similar to that of 
Theorem~\ref{cond-EQ}. One should only consider the projection
${\goth N}^u(\q)\to \g$.

\begin{s}{Theorem}      \label{cond-EQu}
\begin{itemize}
\item[\sf 1.] \ The morphism $\bar\pi: \q\to \bbk^m$
is equidimensional if and only if 
 $\codim_\g \fX_i \ge m-i$.
\item[\sf 2.] \ If\/ ${\goth N}^u(\q)$ is irreducible, then $\bar\pi$
is equidimensional;
\item[\sf 3.] \  ${\goth N}^u(\q)$ is irreducible if and only 
if\/ $\codim_\g \fX_i\ge m-i+1$ for $i<m$.
\end{itemize}
\end{s}%
Now, comparing Theorem~\ref{cond-EQ}(iii), Proposition~\ref{F2},
and Theorem~\ref{cond-EQu}(iii), one concludes that
%%\centerline{\sl 
if  ${\goth N}(\q)$ is irreducible, then so is \ ${\goth N}^u(\q)$.

\vskip.7ex\noindent
But one can derive a much stronger assertion on ${\goth N}^u(\q)$
from the irreducibility of ${\goth N}(\q)$.  This is related to properties of symmetric 
algebras of certain modules over polynomial rings and exploits some technique
from \cite{ls},\,\cite{jac}.

Let $\Mor(\g,V^*)$ be the $\bbk[\g]$-module of {\sl all\/} polynomial
morphisms $F:\g\to V^*$.
Consider the homomorphism
${\hat\tau}:\Mor(\g,V^*){\to} \Mor(\g,V^*)$
defined by $\hat\tau(F)(x)=x{\cdot}F(x)$.
(Here ``$\cdot$'' refers to the $\g$-module structure on $V^*$.)

\begin{s}{Theorem}   \label{ker-tau}
$\ker\hat\tau$ is a free $\bbk[\g]$-module of rank $m$. 
More precisely, $(F_1,\ldots, F_m)$ is a basis for\/ $\ker\hat\tau$.
\end{s}\begin{proof}
The proof is based on the same idea as the proof of Theorem~1.9 in \cite{jac}.
\\[.6ex]
Clearly, $\ker\hat\tau$ is a torsion-free $\bbk[\g]$-module and 
the rank $\rk(\ker\hat\tau):=\dim(\ker\hat\tau\otimes_{\bkk[\g]} \bbk(\g))$ is 
well-defined. 
%%a finite-dimensional $\bbk(\q)$-vector space. 
An easy argument shows that the rank of $\hat\tau$ over $\bbk(\g)$
equals $\dim V-\max_{x\in \g}\dim (V^*)^x=\dim V-m$. Hence 
$\rk(\ker\hat\tau)=m$.
Obviously, $F_i\in\ker\hat\tau$ and $\displaystyle\bigoplus_{i=1}^m \bbk[\g] F_i$
is a free submodule of $\ker\hat\tau$ of rank $m$. It follows that, for
any $F\in \ker\hat\tau$, there exist $\hat p,p_1,\ldots,p_m\in \bbk[\g]$ such that 
\[
     \hat pF=\sum_i p_i F_i \ . 
\]
Assume $\hat p\not\in \bbk^*$. Let $p$ be a prime factor
of $\hat p$ and $D$  the divisor of zeros of $p$.
Then $\sum_i p_i(v) F_i(v)=0$ for any $v\in D$. Since $\g^{reg}$ is big,
$\g^{reg}\cap D$ is dense in $D$. Because 
$\{F_i(v)\}$ are linearly independent for $v\in\g^{reg}$, we obtain
$p_i\vert_D\equiv 0$. Hence $p_i/p\in \bbk[\g]$ for each $i$, and we are
done.
\end{proof}%
Let $E$ denote the $\bbk[\g]$-module $\Ima\hat\tau$. In view of the previous 
theorem, we have the exact sequence
\begin{equation}  \label{exact}
 0\to \bigoplus_{i=1}^m \bbk[\g]\,F_i \stackrel{\hat\beta}
\to\Mor(\g,V^*)\stackrel{\hat\tau}{\to} E\to 0 \ .
\end{equation}
Choose a basis $\xi_1,\ldots,\xi_n$ for $V^*$. Using this basis, we identify
$\Mor(\g,V^*)=\bbk[\g]\otimes V^*$ with $\bbk[\g]^n$.
Then we can write $F_j(x)=\sum_{i=1}^n F_{ij}(x)\xi_i$, where 
$F_{ij}\in \bbk[\g]$.
If we regard sequence~\re{exact} as a sequence
\[
    0\to\bbk[\g]^m\stackrel{\hat\beta}{\to} \bbk[\g]^n\stackrel{\hat\tau}{\to} E \to 0 \ ,
\]
then $\hat\beta$ becomes an $n\times m$-matrix with entries $F_{ij}$.
Let $I_t(\hat\beta)$ be the ideal generated by $t\times t$ minors of $\hat\beta$.
For $d\in \BN$, consider the following condition \\[1ex]
\hbox to \textwidth{\quad 
$(\mathcal F_d)$\hfil $\hot I_t(\hat\beta)\ge m-t+1+d$ \ \ \ \text{for } $1\le t\le m$. \hfil}

\vskip.7ex\noindent
The ideals $I_t(\hat\beta)$ are independent of the presentation of $E$. These
are  Fitting ideals of $E$, see e.g. \cite[1.1]{vasc}.
Let $\textsl{Sym}_{\bkk[\g]}(E)$ denote the symmetric algebra of 
the $\bbk[\g]$-module $E$.
 
\begin{s}{Theorem}   \label{factorial}
Suppose\/ ${\goth N}(\q)$ is irreducible. Then 
\begin{itemize}
\item[\sf (i)] \ The condition $(\mathcal F_2)$ is satisfied by $E$.
\item[\sf (ii)] \ $\text{Sym}_{\bkk[\g]}(E)$ is a factorial domain of Krull dimension
$\dim\g+n-m$.
\item[\sf (iii)] \ ${\goth N}^u(\q)$ is an irreducible
 factorial complete intersection, and\/
$\bbk[{\goth N}^u(\q)]=\text{Sym}_{\bkk[\g]}E$.
\item[\sf (iv)] \ ${\goth N}^u(\q)=\ov{\Ima(\kappa)}$, where 
$\kappa: \q\to \q$ is defined by $\kappa(x,v)=(x,x{\cdot}v)$, $x\in\g,\,v\in V$.
\end{itemize}
\end{s}\begin{proof}
(i) \ It is easily seen that $\fX_i$ is the zero locus of $I_{i+1}(\hat\beta)$.
Therefore condition $(\mathcal F_2)$ is satisfied in view of
Proposition~\ref{F2}.

(ii) \ The exact sequence \re{exact} shows that $E$ has projective dimension at most
one. Therefore part (ii) follows from (i) combined with 
\cite[Prop. 3 \& 6]{luch}.

(iii) \ The universal property of symmetric algebras implies that
$\textsl{Sym}_{\bkk[\g]}(E)$ is the quotient of $\textsl{Sym}_{\bkk[\g]}(
\bbk[\g]\otimes V^*)=
\bbk[\g\times V]$ by the ideal generated by the image of $\hat\beta$.
It follows from the construction that $\hat\beta(F_i)=\wF_i$.
Hence $\textsl{Sym}_{\bkk[\g]}(E)=\bbk[{\goth N}^u(\q)]$, and the other assertions
follow from (ii).

(iv) \ Clearly, $\ov{\Ima(\kappa)}$ is an irreducible subvariety of $\q$.
Taking the (surjective) projection to $\g$ and looking at the dimension of the 
generic fibre, one finds that $\dim \ov{\Ima(\kappa)}=\dim\g+n-m$.
Thus, $\ov{\Ima(\kappa)}\subset {\goth N}^u(\q)$, both have the same dimension and are 
irreducible. Hence they are equal.
\end{proof}%
{\bf Remark.}  For $V=\g$, i.e., for
the Takiff algebra $\g\ltimes\g$, condition $(\mathcal F_2)$ can be proved directly, without
referring to the irreducibility of ${\goth N}(\q)$, see \cite[Prop.\,2.1]{ls}.
In this special case, the above results for ${\goth N}^u(\q)$ are 
already obtained in \cite[Prop.\,2.4]{ls}.
Actually, ${\goth N}(\q)$ is irreducible if $V=\g$. But this fact, as well as
``Takiff'' terminology, was not used in loc.\,cit.
In Section~\ref{z2}, we give new examples of semi-direct products $\q=\g\ltimes V$
such that ${\goth N}(\q)$ is irreducible
and thereby new instances, where Theorem~\ref{factorial} applies.
 
Now, we proceed to case {\bf B)}.
\\
Recall that $\hat G=G\ltimes \g$ and $\hat V=V\ltimes V$ is a 
$\hat G$-module. To a great extent, our results in this case are similar
to those in case {\bf A)}. A notable distinction is, however, that whereas the
adjoint representation of $G$ has some good properties for granted,
we have to require these properties for $(G:V)$.

We will assume below that $(G:V)$ satisfies properties (1)--(3) of
Theorem~\ref{main-twice}, with $G$ in place of $Q$, and use the respective notation.
In particular, $\hat G^u=1\ltimes V$,
$m=\dim V\md G$, $\bbk[V]^G=\bbk[f_1,\ldots, f_m]$,
and $\wF_i\in \bbk[\hat V]^{\hat G}$ is the invariant
associated with $\textsl{d}f_i$. As in case A), we define a stratification
of $V$ by 
\[
  {\fY}_i =\{ x\in V\mid \dim \text{span}\{
(\textsl{d}f_1)_x,\ldots, (\textsl{d}f_m)_x \}=i\}=
\{ x\in V\mid \rk (\textsl{d}\pi_G)_x=i\} \ .
\]
Then $\ov{{\fY}}_i\subset \ov{{\fY}}_{i+1}$ and $\ov{{\fY}}_m=V$.
Notice that $\fY_0=\{0\}$.
The induced stratification of the null-cone  ${\goth N}^G(V)=\nv$ is
$\fY_i(\nv):=\fY_i\cap\nv$. Since $\pi_G: V\to V\md G$ is onto, $\dim{\goth N}(V)
\ge\dim V -\dim V\md G$. But,
unlike the case of $(G:\g)$, it may
happen that the last inequality is strict and $\fY_m(\nv)=\varnothing$.

\begin{s}{Lemma} 
Suppose $\pi_{\hat G}: \hat V \to \hat V\md \hat G$ is equidimensional.
Then so is $\pi_G: V\to V\md G$ and $\fY_m(\nv)\ne\varnothing$.
\end{s}\begin{proof}
Consider the projection $p: {\goth N}(\hat V)\to {\goth N}(V)$.
If $j$ is the maximal index such that $\fY_j(\nv)\ne\varnothing$,
then $\dim  {\goth N}(\hat V)=\dim \nv+\dim V-j$. Since 
$\dim\nv\ge \dim V-m$, the result follows.
\end{proof}%
Thus, if we are searching for equidimensional quotient morphisms $\pi_{\hat G}$,
then we must assume that \\
\hbox to \textwidth{ \ $(\ast)$ \hfil $\dim\nv=\dim V-\dim V\md G$ \ and \ 
$\fY_m(\nv)\ne\varnothing$. \hfil }

\noindent 
In this setting, analogs of results 
\re{cond-EQ}--\re{factorial} are proved in a quite similar fashion.
Let ${\goth N}^u(\hat V)$ denote the zero-fibre of the morphism
$\bar\pi:\hat V \to \bbk^m$ defined by

$\bar\pi(v_1,v_2)=(\wF_1(v_1,v_2),\ldots,\wF_m(v_1,v_2))=
\bigl(\langle (\textsl{d}f_1)_{v_1}, v_2\rangle,\ldots,
\langle (\textsl{d}f_m)_{v_1}, v_2\rangle\bigr)$.

\begin{s}{Theorem}  \label{cond-EQ-b} Under the 
%%previous assumptions on $G$ and $V$, 
assumptions (1)--(3) of Theorem~\ref{main-twice} and $(\ast)$,
we have
\begin{itemize}
\item[\sf 1.] \ The morphism 
$\pi_{\hat G}:  \hat V \to \hat V\md \hat G$ (resp. 
$\bar\pi:  \hat V \to \bbk^m$)
is equidimensional if and only if\/ 
 $\codim_{\nv} \fY_i(\nv) \ge m-i$ (resp. $\codim_{ V} \fY_i \ge m-i$)
for all\/ $i$.
\item[\sf 2.] \ If\/ ${\goth N}(\hat V)$ (resp. ${\goth N}^u(\hat V)$)
is irreducible, then $\pi_{\hat G}$ (resp. $\bar\pi$)
is equidimensional;
\item[\sf 3.] \  ${\goth N}(\hat V)$ (resp. ${\goth N}^u(\hat V)$)
is irreducible if and only if\/ 
$\codim_{\nv} \fY_i(\nv)\ge m-i+1$ 
(resp. $\codim_{V} \fY_i \ge m-i+1$) for\/ $i<m$.
\item[\sf 4.] \ If\/ $\codim_{\nv} \fY_i(\nv)\ge a$,
%% for\/ $i<m$, 
then $\codim_{V} \fY_i \ge m-i+a$.
\end{itemize}
\end{s}\begin{proof} The proof of parts {\sf 1--3} is similar to the proof of
Theorem~\ref{cond-EQ}. For the last part, we notice that
$\dim \ov{\fY}_i\md G \le i$. Therefore $\dim\fY_i\le i+\dim \fY_i(\nv)\le
i+\dim\nv-a=\dim V-(m-i+a)$. (Cf. the proof of Prop.\,2.1 in \cite{ls}.)
\end{proof}%
%
%%After some preparations, we state below an analogue of Theorem~\ref{factorial}.
%%\\[.6ex]
Consider the homomorphism of $\bbk[V]$-modules
\[
{\hat\mu}: \Mor(V,V^*){\to} \Mor(V,\g^*) 
\]  
defined by $\langle \hat\mu(F)(v),s\rangle:=\langle F(v), s{\cdot}v\rangle$ 
for $v\in V,s\in\g$. 
Here ``$\cdot$'' refers to the $\g$-module structure on $V$ and 
the first (resp. second) $\langle\ ,\ \rangle$
stands for the pairing of $\g$ and $\g^*$ (resp. $V$ and $V^*$).
By \cite[theorem\,1.9]{jac}, $\ker\hat\mu$ is a free $\bbk[V]$-module of
rank $m$ generated by $\textsl{d}f_i$, $i=1,\ldots,m$. Let $\hat E$ denote the
$\bbk[V]$-module $\Ima\hat\mu$.

\begin{s}{Theorem}  \label{factorialB}
Suppose\/ ${\goth N}(\hat V)$ is irreducible. Then 
\begin{itemize}
\item[\sf (i)] \ The condition $(\mathcal F_2)$ is satisfied by $\hat E$.
\item[\sf (ii)] \ $\text{Sym}_{\bkk[V]}(\hat E)$ is a factorial domain of Krull dimension
$2\dim V-m$.
\item[\sf (iii)] \ ${\goth N}^u(\hat V)$ is an irreducible
 factorial complete intersection, and\/
$\bbk[{\goth N}^u(\hat V)]=\text{Sym}_{\bkk[V]}\hat E$.
\item[\sf (iv)] \ ${\goth N}^u(\hat V)=\ov{\Ima(\varkappa)}$, where 
$\varkappa: V\times \g\to V\times V$ is defined by 
$\varkappa(v,x)=(v,x{\cdot}v)$, $x\in\g,\,v\in V$.
\end{itemize}
\end{s}%\begin{proof}
%\end{proof}%
%
The proof of Theorem~%%\ref{cond-EQ-b} and 
\ref{factorialB} is omitted, since it
is similar to the proof of Theorem~%%\ref{cond-EQ} and 
\ref{factorial}.
%%respectively.

                %%%%
%%%%%%%%%%%%%%%%
%%%%%%%%%%%%%%%%   Section 9
%%%%%%%%%%%%%%%%
                %%%%

\section{Isotropy contractions and $\mathbb Z_2$-contractions of semisimple Lie algebras}  
\label{z2}
\setcounter{equation}{0}

\noindent
Let $\h$ be a subalgebra of $\q$ such that $\q=\h\oplus\mathfrak m$ 
for some $\ad\h$-stable subspace $\mathfrak m\subset \q$.
(Such an $\h$ is said to be {\it reductive in\/} $\q$.) 
For instance, if $\vartheta$ is an involutory automorphism of $\q$, then $+1$ and 
$-1$-eigenspaces of $\vartheta$ yield such a decomposition. 
The fixed-point subalgebra of an involutory automorphism is called a symmetric 
subalgebra.

\begin{rem}{Definition}   \label{def:degener}
If $\h$ is reductive in $\q$, then
the semi-direct product 
$\h\ltimes\mathfrak m$ is called an {\it isotropy contraction\/} 
of $\q$. If $\h$ is symmetric,
so the decomposition $\q=\h\oplus\mathfrak m$ is a $\mathbb Z_2$-grading,
then $\h\ltimes\mathfrak m$ is also called a $\mathbb Z_2$-{\it contraction\/}
of $\q$.
\end{rem}%
Notice that $\h\ltimes\mathfrak m$ is a contraction of $\q$
in the  sense of the deformation theory of Lie algebras, see e.g.
\cite[Chapter 7, \S\,2]{t41}.  
More precisely, consider the invertible
linear map $c_t: \q\to \q$, $t\in \bbk\setminus\{0\}$, 
such that $c_t(h+m)=h+t^{-1}m$.  
Define the new Lie algebra multiplication $[\ ,\ ]_{(t)}$
on the vector space $\q$ by the rule
\[
     %%[h_1+m_1, h_2+m_2]_{(t)}=[h_1,h_2]+[h_1,m_2]+[m_1,h_2]+t[m_1,m_2] \ ,
     [x,y]_{(t)}:= c_t\bigl( [ c_t^{-1}(x), c_t^{-1}(y)]\bigr), \quad x,y\in\q \ .
\]
Then, for all $t\ne 0$,  the algebras $\q_{(t)}$ are isomorphic, and
$\lim_{t\to 0}\q_{(t)}=\h\ltimes\mathfrak m$.

\begin{s}{Lemma}  \label{takiff-contr}
Any Takiff Lie algebra is a $\mathbb Z_2$-contraction.
\end{s}\begin{proof}
Consider the direct sum of Lie algebras $\q\dotplus\q$ and the involution
$\vartheta$ permuting the summands. Then the corresponding $\mathbb Z_2$-contraction
is isomorphic to $\q\ltimes \q$.
\end{proof}%
In the rest of the section, we only consider isotropy contractions such that
the initial ambient Lie algebra is semisimple and the subalgebra is reductive.
Let $\ka=\h\ltimes \mathfrak m$ be an isotropy contraction of a
semisimple Lie algebra $\g$. For $\g$, one has equalities 
\[
    \rk\g=\ind\g=\dim\g\md G  \ .
\]
The first natural question is: 

{\it To which extent this remains true for isotropy contractions?}
\\[.7ex]
Recall that the {\it complexity\/} of a homogeneous space $G/H$, denoted
$c(G/H)$, equals $\trdeg\bbk(G/H)^B$, where $B$ is a Borel subgroup of $G$,
and $G/H$ is said to be {\it spherical\/} if $c(G/H)=0$. We refer to \cite{these}
for basic facts on complexity.

\begin{s}{Proposition}  \label{facts}

(1)  \ We have $\ind\ka = \ind \g+ 2c(G/H)$. In particular, $\ind\ka=\ind\g$ if and 
only if $H$ is a spherical subgroup of $G$.

(2) \  $\dim\ka\md K=\dim \z_\g(x)$, where $x\in \h$ is an $\h$-regular semisimple
element.
\end{s}\begin{proof}
(1)  By \cite{rais},  $\ind\ka=\trdeg\bbk(\ma^*)^H+ \ind \es$,
where $\es$ is a generic stabiliser for $(H:\ma^*)$. Since $\ma$ is an orthogonal
$\h$-module, there is no difference between $\ma$ and $\ma^*$, the action $(H:\ma)$
is stable \cite{lu72} 
and therefore $\es$ is reductive. Hence $\ind\ka=\dim \ma\md H+\rk\es$.
On the other hand,
there is a formula for $c(G/H)$ in terms of the isotropy representation
$(H:\ma)$. Namely, $2c(G/H)=\dim\ma\md H-\rk\g+\rk\es$ \cite[Cor.\,2.2.9]{these}.
Hence the conclusion.

(2) By Theorem~\ref{main-reduct}, $\dim\ka\md K=\rk\h+ \dim \ma^{\te_\h}$.
The latter equals $\dim\z_\h(x)+\dim\ma^x$ for 
a regular semisimple element $x\in\te_\h\subset \h$.
\end{proof}%
{\bf Remark.} It is a general fact that the index of a Lie algebra
cannot decrease under contraction. The previous result gives a precise
meaning for this in case of isotropy contractions.

\begin{s}{Corollary}   \label{good-z2}
If\/ $\g=\g_0\oplus\g_1$ is a $\BZ_2$-grading and
$\ka=\g_0\ltimes\g_1$ is the respective $\BZ_2$-contraction, then\/ 
$\ind\ka=\dim\ka\md K=\rk\g$.
\end{s}\begin{proof}
As is well known, any symmetric subgroup $G_0\subset G$ is spherical, and
$\g_0$ contains a regular semisimple element of $\g$.
\end{proof}%
Thus, for $\BZ_2$-contractions one obtains two, usually
different, decompositions of
the rank of $\g$:
\[
\rk\g=\left\{ \begin{array}{l}  \ind\ka=\rk\es + \dim \g_1\md G_0 ; \\
                           \dim\ka\md K=\rk\g_0+\dim(\g_1)^{\te_0} ,
\end{array}\right.
\]
where $\te_0$ is a Cartan subalgebra of $\g_0$. 

\noindent
If $\h$ contains a $\g$-regular semisimple element, then
$\bbk[\g]^G$ and $\bbk[\ka]^K$ are graded polynomial algebras of the same
Krull dimension. The second natural question is: 

{\it Is there a relationship between the degrees of free homogeneous
 generators (basic invariants) ?}
\\[.6ex]
Let $\textsl{Deg}(A)$ denote the multiset of degrees of free
generators of a graded polynomial algebra $A$. The elements of $\textsl{Deg}(A)$
are assumed to be increasingly ordered.

\begin{s}{Theorem}  \label{deform1}
(1) \ If $\h$ contains a $\g$-regular semisimple element, then
$\text{Deg}(\bbk[\ka]^K) \le \text{Deg}(\bbk[\g]^G)$ (componentwise inequalities).

(2) \ Suppose a regular nilpotent element of\/ $\h$ is also regular in $\g$. Then  
\\ 
$\textsl{g-exp}_\g(\g)=\textsl{g-exp}_\h(\h)\uplus \textsl{g-exp}_\h(\ma)$ (the union 
of multisets). Equivalently, $\text{Deg}(\bbk[\ka]^K) = \text{Deg}(\bbk[\g]^G)$.
\end{s}\begin{proof}
(1) \ Recall that $\ka=\lim_{t\to 0}\g_{(t)}$. It is easily seen that this contraction
gives rise to "a curve in the space of algebras of invariants" and 
to an embedding $\lim_{t\to 0}\bbk[\g_{(t)}]^{G_{(t)}}\subset \bbk[\ka]^K$.
The limit exists, because $\bbk[\g_{(t)}]^{G_{(t)}}$ is graded and the (finite)
dimension the homogeneous component of a given degree
does not depend on $t$; so that the limit
is taken in a suitable Grassmannian.

(2) \ Let $\{e,h,f\}$ be a principal $\tri$-triple in $\h$
(see \cite[Ch.\,6,\S\,2.3]{t41}). By the assumption,
it is also a principal $\tri$-triple in $\g$.
By a result of R.\,Brylinski \cite{ranee}, the generalised exponents of a $G$-module $V$ 
are obtained as follows.
Take the subspace $V^T$ and its ``$e$-limit" $\lim_e(V^T)\subset V$, see 
\cite[\S\,2]{ranee} 
for the precise definition. Then $\textsl{g-exp}_\g(V)$ is the multiset of 
$h$-eigenvalues on $\lim_e(V^T)$. It is important that this ``$e$-limit" depends
only on the $\{e,h\}$-module structure on $V$.
In our setting, $\g$ and $\ka$ are isomorphic as $\h$-modules, and $\ka=\h\oplus\ma$ 
as $\h$-module. Therefore 
\[
 \textsl{g-exp}_\g(\g)=\textsl{g-exp}_\h(\g)=
 \textsl{g-exp}_\h(\ka)=\textsl{g-exp}_\h(\h)\uplus \textsl{g-exp}_\h(\ma)  \ .
\]
The second assertion follows from Theorem~\ref{main-reduct}, because
$\textsl{Deg}(\bbk[\ka]^K)=\textsl{g-exp}_\h(\ka)+1$ (componentwise) and likewise for
$\bbk[\g]^G$. 
\end{proof}%
Part (2) of this theorem can be used for finding generalised exponents of certain
representations.

\begin{rem}{Example}   \label{D4G2}
Let $\g$ be ${\frak so}_8$ and $\h$ the exceptional Lie algebra of type
$\GR{G}{2}$ ($\dim\g=28$, $\dim\h=14$). The restriction of the defining representation
of $\g$ to $\h$ is the sum of $V(7)$, the 7-dimensional simple $\h$-module,
and a 1-dimensional trivial module. Let $e\in\h$ be a regular nilpotent element.
It is known that $V(7)$ is a cyclic $e$-module. Therefore, as element of
${\frak so}_8$, $e$ has the Jordan form with blocks of size 7 and 1.
Hence $e$ is also regular in ${\frak so}_8$.
Here $\ma=V(7)\oplus V(7)$.  
Since $\textsl{g-exp}_\g(\g)=\{1,3,3,5\}$ and $\textsl{g-exp}_\h(\h)=\{1,5\}$,
we conclude that $\textsl{g-exp}_\h(V(7))=\{3\}$. That is, the $\bbk[\h]^H$-module
$\Mor_\h(\h, V(7))$ is generated by the covariant of degree 3.

This is also an instructive illustration to Theorem~\ref{cond-EQ} and
Corollary~\ref{m12}. Here $m=\rk\g-\rk\h=2$, hence $\pi_K$ is equidimensional. 
The basic covariant in $\Mor_\h(\h, V(7))$ vanishes on the subregular 
nilpotent orbit in $\N(\h)$. This follows from a result of Broer on the ideal 
defining the closure of the subregular nilpotent orbit \cite[\S\,4]{broer}.
Therefore $\codim_{\N(\h)}\fX_0(\N(\h))=2$ and ${\goth N}(\ka)$ appears to be
reducible.
\end{rem}%
From now on, we assume that $\ka$ is a $\BZ_2$-contraction of $\g$.

\begin{s}{Theorem}  \label{main-z2}
Let\/ $\g=\g_0\oplus\g_1$ be a $\BZ_2$-graded semisimple Lie algebra and\/  
$\ka=\g_0\ltimes\g_1$ its $\mathbb Z_2$-contraction. 
%%its adjoint representation. 
Then ${\goth N}(\ka)$ is irreducible.
\end{s}\begin{proof}
Let $\vartheta$ be the involution of $\g$ determining the $\mathbb Z_2$-grading.
It suffices to handle the case in which $\g$ is not a sum of $\vartheta$-stable ideals.
This means that either $\g$ is simple or $\g=\es\dotplus\es$, where $\es$ is simple and
$\vartheta$ permutes the factors. In the second case, $\ka=\es\ltimes\es$ is a 
Takiff Lie algebra, and the required result is proved in \cite[Theorem\,2.4]{geof1}.
Therefore we concentrate on the first case.

From now on, $\g$ is simple. Write $\N_0$ for the null-cone in $\g_0$ and
$K$ for the Takiff group $G_0\ltimes\g_1$. Since $\g_1$ is an orthogonal $G_0$-module,
we do not distinguish $\g_1$ and $(\g_1)^*$.
%%Our proof is based on 

1) Suppose $\vartheta$ is inner. Then $\rk\g=\rk\g_0$ and therefore
the $\g_0$-module $\g_1$ has no zero weight space. As is noted in Section~\ref{null-cone},
the null-cone ${\goth N}(\ka)$ is then isomorphic to $\N_0\times \g_1$.

2) Suppose $\vartheta$ is outer. This is the difficult part of the proof, which relies
on the classification of the involutions of simple Lie algebras.
Recall that $m=\dim(\g_1)^{\te_0}=\rk\g-\rk\g_0$.

$({a_1})$ \ $\rk\g_0=\rk\g-1$ and $m=1$. Here the assertion follows from
Corollary~\ref{m12}.  This happens if $\g=\sone$ and $\g_0={\frak so}_{2k+1}\times
{\frak so}_{2l+1}$ with $k+l=n-1$.

$(a_2)$ \ $\rk\g_0=\rk\g-2$ and $m=2$. By 
Corollary~\ref{m12}, $\pi_K$ is equidimensional. Still, ${\goth N}(\ka)$
can be reducible {\sl a priori}. To prove that this is not the case,
consider the hierarchy 
$\fX_0(\N_0)\subset \ov{\fX_1(\N_0)}\subset\ov{\fX_2(\N_0)}=\N_0$
determined by the basic covariants of type $\g_1$.
Invoking  the criterion of irreducibility
(Theorem~\ref{cond-EQ}(iii)) with $m=2$ shows that
only the condition with $i=0$ has to be satisfied. That is, we must have
$\codim_{\N_0}\fX_0(\N_0)\ge 3$. This means that each nilpotent orbit in
$\N_0$ of codimension 2 does not belong to $\fX_0(\N_0)$, i.e., there should exist a 
covariant $F\in \Mor_{G_0}(\g_0,\g_1)$ that does not vanish on such an orbit.

There are two involutions with $m=2$ in the exceptional
algebras. In both cases, $\g$ is of type 
$\GR{E}{6}$ and $\g_0$ is either $\GR{F}{4}$ or $\GR{C}{4}$.
Furthermore, the degrees of basic covariants of type $\g_1$ are $4,\,8$ in both cases.
Since $\g_0$ is simple here, $\N_0$ has a unique orbit of codimension~2,
the so-called {\it subregular\/} nilpotent orbit $\co_{sub}$. The closure of $\co_{sub}$
is normal and the equations of $\ov{\co}_{sub}$ in $\bbk[\N_0]$
are explicitly described,
%%algebra of functions on it has an explicit description,
see \cite[\S\,4]{broer}. Therefore, it is not hard to verify that the covariant of degree
4 survives on $\co_{sub}$.

$({a_3})$ \ It remains to handle two series of $(\g,\g_0)$:
$({\frak sl}_{2n}, \spn)$ and $(\sln, \son)$.
In  these cases, we explicitly describe the covariants of type $\g_1$
and verify that the condition of Theorem~\ref{cond-EQ}(iii) is satisfied.
Actually, we show that, for all $\mathbb Z_2$-contractions of simple Lie
algebras, a stronger inequality holds,
see Eq.~\re{brilliant} below.
\end{proof}%
Let us adapt Theorem~\ref{cond-EQ} to our setting. We consider the 
stratification of $\N_0$ determined by covariants of type $\g_1$.
Since $\N_0$ consists of finitely many $G_0$-orbits,
condition~\ref{cond-EQ}(iii) can be verified for each orbit separately.
%%$\dim \fX_i(\N_0)=\max \dim \co$, where $\co$ runs over the
%%$G_0$-orbits in $\fX_i(\N_0)$.
Therefore, it can be written as
\begin{equation}   \label{old-brilliant}
  \dim\z_{\g_0}(x) - \rk\g_0 >m-\dim \text{span}\{F_1(x),\ldots,F_m(x)\}  
\ \ \text{ if } \ x\in \N_0\setminus \fX_m(\N_0) \ ,
\end{equation}
cf. Remark~\ref{remarks-8}(1). What we are going to prove is:
\begin{equation}   \label{brilliant}
 \dim\z_{\g_0}(x) - \rk\g_0 \ge 2 \bigl(m-\dim \text{span}\{F_1(x),\ldots,F_m(x)\} \bigr)
\ \ \text{ for any } \ x\in \N_0 \ .
\end{equation}
Clearly, the last version is stronger and has an advantage of being stated more
uniformly. 

\begin{s}{Theorem}  \label{thm-brill}
Inequality~\re{brilliant} holds for any $\BZ_2$-grading of a simple Lie algebra $\g$.
\end{s}\begin{proof}
Since the difference in the left-hand side of \re{brilliant} is always even,
there is no distinction between inequalities~\re{old-brilliant} and
\re{brilliant} for ${m\le 2}$.
Therefore the proof of Theorem~\ref{main-z2} shows that
it remains to verify Eq.~\re{brilliant} for the following series
of $\mathbb Z_2$-gradings: 
%%and respective $\mathbb Z_2$-contractions:
\begin{itemize}
\item[$\bullet$]  \ $\g_0={\frak sp}(V)$, $\g_1=\wedge^2_0(V)$, $\dim V=2n$.
\item[$\bullet$]  \ $\g_0={\frak so}(V)$, $\g_1=\mathcal S^2_0(V)$. Here one actually 
has two series, depending on the parity of $\dim V$.
\end{itemize}
We use familiar matrix models of classical Lie algebras and their representations.
In the following computations, we need the fact that the nilpotent 
$G_0$-orbits are classified by certain partitions of $\dim V$, see \cite[IV.2.15]{ss},
\cite[Ch.\,6 \S 2.2]{t41}.
A minor unpleasant phenomenon related to $\sone$ is that there are
two isomorphic $SO_{2n}$-orbits corresponding to a
``very even partition". This does not affect, however, our computations.
For $x\in \N_0$, let
$\boldsymbol{\eta}=(\eta_1,\eta_2,\ldots)$ denote the corresponding partition.  
Write $(\hat\eta_1,\hat\eta_2,\ldots,\hat\eta_s)$ for the dual partition. 
This means in particular that $s=\eta_1$. What we need from these partitions
is an explicit formula for $\dim\z_{\g_0}(x)$ and a way to determine
$i$ such that $x\in \fX_i(\N_0)$.

Let us begin with the symplectic case. Let $J$\/ be a skew-symmetric non-degenerate
bilinear form on $V$, which is identified with its matrix in a certain basis for
$V$. Then 

$-$ \ $\spn=\spv={\frak sp}(V,J)$ is the space
of matrices $\{ x \in\glv\mid xJ \text{ is symmetric}\}$;

$-$ \ the representation space $\wedge^2_0(V)$ can be regarded as the 
space of skew-symmetric matrices modulo one-dimensional subspace generated by 
$J$. The $\spn$-action on the space of skew-symmetric matrices is given by
$(x,A)\mapsto xJA+A(xJ)^t$. 
\\[.6ex]
In this case $m=n-1$, i.e., there are $n-1$ basic covariants of type
$\g_1$. Since any regular nilpotent element in $\spn$ is also regular in 
${\frak sl}_{2n}$, the generalised exponents of the $\g_0$-module $\g_1$
can be found using Theorem~\ref{deform1}(2).
These are $2,4,\ldots, 2n-2$.
The key observation is that the corresponding covariants have a very simple expression.
Namely, consider the maps $(x\in\spn)\mapsto F_i(x)=x^{2i}J$, $i=1,\ldots, n-1$. 
It is easily seen that $x^{2i}J$ is skew-symmetric 
and each $F_i$ is $Sp_{2n}$-equivariant. Because the $F_i$'s are linearly independent
over $\bbk[\g_0]^{G_0}$, these are precisely the basic covariants.

\begin{s}{Proposition}  \label{Sp-case}
Inequality~\re{brilliant} is satisfied for $(Sp(V), \wedge^2_0(V))$.
\end{s}\begin{proof}
By \cite[Corollary\,3.8(a)]{hess},
the dimension of the centraliser of $x$ in $\g_0=\spn$
is given by the formula \ 
$\displaystyle 
  \dim\z_{\g_0}(x)=\frac{1}{2}(\sum_i\hat\eta_i^2 + \#\{j\mid \eta_j \text{ is odd}\})$.
\\[.7ex]
The maximal nonzero power of $x$ is determined by the size of the maximal Jordan 
block, i.e., $\eta_1$. Therefore $x\in \fX_i(\N_0)$ if and only if $x^{2i}\ne 0$
and $x^{2i+2}=0$ 
if and only if $[\frac{\eta_1-1}{2}]=i$.
Hence inequality~\re{brilliant}, which we wish to prove, can be written as
\[
  2\left[\frac{\eta_1-1}{2}\right]+
\frac{1}{2}\bigl(\sum_{i=1}^s\hat\eta_i^2 +\#\{j\mid \eta_j \text{ is odd}\}\bigr)-n-2(n-1)\ge 0\ .
\]
Using the relations $\sum \hat\eta_i=2n$ and
$\eta_1=s$, the left-hand 
side is transformed as follows:
\begin{multline*}
2\left[\frac{\eta_1-1}{2}\right]+
\frac{1}{2}\bigl(\sum_{i=1}^s\hat\eta_i^2 + \#\{j\mid \eta_j \text{ is odd}\}\bigr)-
\frac{3}{2}\sum\hat\eta_i+2 = 
\\
2\left[\frac{s-1}{2}\right]+
\frac{1}{2}\bigl(\sum_{i=1}^s(\hat\eta_i^2-3\hat\eta_i) + 
\#\{j\mid \eta_j \text{ is odd}\}\bigr)+2= 
\\
\frac{1}{2}\bigl(\sum_{i=1}^s(\hat\eta_i-1)(\hat\eta_i-2)+\#\{j\mid \eta_j \text{ is odd}\}\bigr)
+2\left[\frac{s+1}{2}\right]-s \ .
\end{multline*}
The first group of summands is non-negative, and so is the last group.
Thus, inequality~\re{brilliant} holds for any nilpotent orbit in $\spn$.
\end{proof}%
We continue with the orthogonal case, with $\dim V=N$.
Here $\g_0$ is the space of skew-symmetric $N{\times} N$-matrices and $\g_1=\mathcal S^2_0(V)$
is the space of traceless symmetric $N{\times} N$-matrices.
\\[.6ex]
If $N=2n+1$, then $m=n$. In this case, a regular
nilpotent element of $\sono$ is also regular in ${\frak sl}_{2n+1}$,
so that  Theorem~\ref{deform1}(2)
applies, and $\textsl{g-exp}_{G_0}(\g_1)=\{2,4,\ldots,2n\}$.
Similarly to the symplectic case, we find that
$x\mapsto F_i(x)=x^{2i}$, $i=1,2,\ldots,n$, are the basic covariants.
If $N=2n$, then $m=n-1$. A regular nilpotent element of $\sone$ is not
regular in ${\frak sl}_{2n}$, but 
%%an easy {\sl ad hoc\/} argument shows that
$F_1,\ldots,F_{n-1}$ are still the basic covariants.
For, the $F_i$'s are linearly independent over $\bbk[\g_1]^{G_0}$ and neither
of them vanishes on the regular nilpotent orbit in $\sone$.

\begin{s}{Proposition}  \label{SO-case}
Inequality~\re{brilliant} is satisfied for $(SO(V), \mathcal S^2_0(V))$.
\end{s}\begin{proof}
By \cite[Corollary\,3.8(a)]{hess},
the dimension of the centraliser of $x$ in $\g_0={\frak so}_N$
is given by the formula \ 
$\displaystyle
  \dim\z_{\g_0}(x)=\frac{1}{2}(\sum_i\hat\eta_i^2 -\#\{j\mid \eta_j \text{ is odd}\})$.
The constraints imposed on partitions in the orthogonal case imply that $\hat\eta_1
\equiv \dim V\pmod 2$.
\\[.6ex] 
The maximal nonzero power of $x$ is determined by 
the size of the maximal Jordan 
block. Therefore $x\in \fX_i(\N_0)$ if and only if $x^{2i}\ne 0$ and $x^{2i+2}=0$
if and only if $[\frac{\eta_1-1}{2}]=i$.
The following computations are slightly different for $\sono$ and $\sone$.

{\sf 1.} \ $N=2n+1$. Here 
inequality~\re{brilliant} can be written as
\[
  2\left[\frac{\eta_1-1}{2}\right]+
\frac{1}{2}\bigl(\sum_i\hat\eta_i^2 - \#\{j\mid \eta_j \text{ is odd}\}\bigr)-3n \ge 0 \ .
\]  
Using the relations $\sum \hat\eta_i=2n+1$ and $\eta_1=s$, 
the left-hand side is transformed as follows:  
%%\newpage
\begin{multline*}
2\left[\frac{s-1}{2}\right]+
\frac{1}{2}\bigl(\sum_{i=1}^s\hat\eta_i^2 - \#\{j\mid \eta_j \text{ is odd}\}\bigr)-
\frac{3}{2}(\sum_{i=1}^s\hat\eta_i - 1) = 
\\
2\left[\frac{s-1}{2}\right]+
\frac{1}{2}\bigl(\sum_{i=1}^s(\hat\eta_i^2-3\hat\eta_i+2) - 2s +3 -
\#\{j\mid \eta_j \text{ is odd}\}\bigr)= 
\\
\frac{1}{2}\bigl(\sum_{i=1}^s(\hat\eta_i-1)(\hat\eta_i-2)-\#\{j\mid \eta_j \text{ is odd}\}
+4\left[\frac{s+1}{2}\right]-2s-1\bigr)=:\mathsf L \ .
\end{multline*}
To see that  $\mathsf L$ is nonnegative, consider several cases.

(a) \ $\hat\eta_1=1$ and hence all $\hat\eta_i=1$. Then  $s=2n+1$ and
$\mathsf L=0$.

(b) \ $\hat\eta_1=3$ and therefore $\boldsymbol{\eta}=(\eta_1,\eta_2,\eta_3)$.
Then $\sum_{i=1}^s(\hat\eta_i-1)(\hat\eta_i-2)=2\eta_3$. Hence
\[
 \mathsf L=\eta_3+2\left[\frac{\eta_1+1}{2}\right]-\eta_1-
\frac{1}{2}(1+\#\{j\mid \eta_j \text{ is odd}\}) \ .
\]
Taking into account that the even parts in $(\eta_1,\eta_2,\eta_3)$ occur
pairwise and $\eta_1+\eta_2+\eta_3$ is odd, 
one quickly verifies that $\mathsf L$ is always nonnegative.

(c) \ $\hat\eta_1\ge 5$. Then $\sum_{i=1}^s(\hat\eta_i-1)(\hat\eta_i-2) \ge \hat\eta_1+2
\ge \#\{j\mid \eta_j \text{ is odd}\}+2$. 
Next, \\ 
$4\left[\frac{s+1}{2}\right]-2s-1\ge -1$. Hence $\mathsf L$ is positive.
\\[.6ex]
Thus, inequality~\re{brilliant} holds for any nilpotent orbit in $\sono$.

{\sf 2.} \ $N=2n$. 
Here the inequality we need to prove reads
\[
    2\left[\frac{\eta_1-1}{2}\right]+
\frac{1}{2}\bigl(\sum_i\hat\eta_i^2 - \#\{j\mid \eta_j \text{ is odd}\}\bigr)-n-2(n-1) \ge 0 \ .
\]
Using the relations $\sum \hat\eta_i=2n$ and $\eta_1=s$, 
the left-hand side is being transformed to
\[
 \frac{1}{2}\bigl(\sum_{i=1}^s(\hat\eta_i-1)(\hat\eta_i-2)-\#\{j\mid \eta_j \text{ is odd}\}\bigr)
+2\left[\frac{s+1}{2}\right]-s=:\mathsf L \ .   
\]
Again, consider some cases. 
%%Here the constraints imposed on partitions in the orthogonal case imply that $\hat\eta_1$ 
%%(the number of Jordan blocks) is even. 

(a) $\hat\eta_1=2$, i.e., $x$ has only two Jordan blocks $(\eta_1,\eta_2)$.
Then $\eta_1,\eta_2$ have the same parity, and in both cases $\mathsf L=0$. 

(b) $\hat\eta_1\ge 4$. Then $\sum_{i=1}^s(\hat\eta_i-1)(\hat\eta_i-2) \ge \hat\eta_1+2
> \#\{j\mid \eta_j \text{ is odd}\}$. Since $2\left[\frac{s+1}{2}\right]-s\ge 0$, the total 
expression is  positive.
\\[.6ex]
Thus, inequality~\re{brilliant} holds for any nilpotent orbit in $\sone$.
\end{proof}%
This completes the proof of Theorem~\ref{thm-brill}. 
\end{proof}%
Thus, all verifications needed to complete the proof of Theorem~\ref{main-z2}
are done. Below, we gather our results on $\mathbb Z_2$-contractions of
semisimple Lie algebras.

\begin{s}{Theorem}   \label{main-free}
Let $\ka=\g_0\ltimes\g_1$ be a $\BZ_2$-contraction of a semisimple Lie algebra
$\g$. Then 
\begin{itemize}
\item[\sf (1)] \  $\bbk[\ka]^K$ is a polynomial algebra of Krull dimension $\rk\g$;
\item[\sf (2)] \  ${\goth N}(\ka)$ is an irreducible complete intersection. 
If\/ $\bbk[\ka]^K=\bbk[f_1,\ldots, f_l]$, $l=\rk\g$, then 
the ideal of ${\goth N}(\ka)$ in $\bbk[\ka]$ is generated by $f_1,\ldots, f_l$.
\item[\sf (3)] \  the quotient morphism $\pi_K:\ka\to \ka\md K$ is equidimensional;
\item[\sf (4)] \ $\bbk[\ka]$ is a free $\bbk[\ka]^K$-module.
\item[\sf (5)] \ if $\kappa: \g_0\oplus\g_1 \to \g_0\oplus\g_1$ is defined
by $\kappa(x_0,x_1)=(x_0,[x_0,x_1])$, then $\ov{\Ima\kappa}={\goth N}^u(\q)$
and it is a factorial complete intersection of codimension $\rk\g-\rk\g_0$.
\item[\sf (6)] \ the coadjoint representation of $\ka$ has a generic stabiliser.
\end{itemize}
\end{s}\begin{proof}
Part (1) follows from Theorem~\ref{main-reduct}. 
The irreducibility in Part~(2) is just Theorem~\ref{main-z2}. 
Let $x\in\N_0$ be a regular nilpotent element. Then $\tilde x=(x,0)\in {\goth N}(\ka)$, and
the description of basic invariants $f_1,\ldots, f_l$ in Theorem~\ref{main-reduct}
shows that  
%%$x\in {\goth N}(\ka)$ such that 
$(\textsl{d}f_1)_{\tilde x},\ldots,(\textsl{d}f_l)_{\tilde x}$ are linearly independent.
Then a standard argument shows that the ideal of 
${\goth N}(\ka)$ is generated by $f_1,\ldots, f_l$ (cf. \cite[Prop.\,6]{ko63}.)
Part~(3) follows from (2) and Theorem~\ref{cond-EQ}(2).
Part~(4) is a formal consequence of Parts (1) and (3).
Part~(5) follows from Theorem~\ref{factorial} and the irreducibility of 
${\goth N}(\ka)$. Since the isotropy representation of
any symmetric subalgebra of $\g$ is polar, part~(6) follows from Theorem~\ref{main-Q2}.
\end{proof}%
To prove the irreducibility of ${\goth N}(\ka)$,  inequality~\re{old-brilliant}
is sufficient. However, our efforts in proving stronger
inequality~\re{brilliant} are not in vain, because that result 
also has a geometric meaning.

\begin{s}{Theorem}  \label{not-vain}
Let\/ $\g=\g_0\oplus\g_1$ be a $\BZ_2$-grading. Consider
the semi-direct product $\tilde\ka=\g_0\ltimes (\g_1\oplus\g_1)$
and the corresponding adjoint representation $(\tilde K:\tilde\ka)$.
Then the quotient morphism $\pi_{\tilde K}$ is equidimensional.
\end{s}\begin{proof}
The criterion for equidimensionality of $\pi_{\tilde K}$, 
Theorem~\ref{cond-EQ}(i), written out in
this case yields precisely inequality~\re{brilliant}.
\end{proof}%
Main efforts in this section were devoted to 
$\BZ_2$-contractions of $\g$. However, there are interesting examples of 
other isotropy contractions with full bunch of good properties.

\begin{rem}{Examples}  \label{B3G2}
1. Suppose $\g={\frak so}_7$ and $\h$ is a simple subalgebra of type 
$\GR{G}{2}$. It is a "truncation" of Example~\ref{D4G2}. Here $\ma=V(7)$,
and one easily verifies that all conclusions of 
Theorem~\ref{main-free} hold for $\ka=\h\ltimes\ma$.

2. $\g={\frak sl}_{2n+1}$ and $\h=\spn=\spv$. Here the $\spv$-module
$\ma$ equals $\wedge^2(V)\oplus V\oplus V$. Since the $\spv$-module $V$
has no zero-weight space, the structure of ${\goth N}(\h\ltimes\ma)$
is essentially the same as for the $\BZ_2$-contraction of the symmetric
pair $({\frak sl}_{2n},\spn)$.
\end{rem}%
{\bf Remark.} Our proofs of Theorems~\ref{main-z2} and
\ref{thm-brill} use classification of involutory automorphisms
and explicit considerations of cases. It would be extremely interesting
to find a case-free proof for the irreducibility of ${\goth N}(\ka)$. 
Especially, because the corresponding irreducibility
result for the Takiff algebra $\g\ltimes\g$ can derived without checking cases. 
We discuss this topic in the following section.

                %%%%
%%%%%%%%%%%%%%%%
%%%%%%%%%%%%%%%%   Section 10
%%%%%%%%%%%%%%%%
                %%%%

\section{Reductive Takiff Lie algebras and their representations}  
\label{tt}
\setcounter{equation}{0}

\noindent
The attentive reader may have noticed that we stated and proved the stronger
inequality~\re{brilliant} only for the $\BZ_2$-gradings of {\sf simple} Lie algebras,
leaving aside the permutation of two factors in $\g\times \g$ and the corresponding
Takiff algebra $\hat\g$.

The situation here is as follows.
By Theorem~\ref{cond-EQ}, the counterpart of inequality~\re{old-brilliant}
for $\hat\g$ is equivalent to the irreducibility of $\ntg$, and this 
was already proved by Geoffriau \cite{geof1}. His proof
%%of \re{old-brilliant} for Takiff algebras 
consists of explicit
verifications for all simple types. It was noticed by M.\,Brion \cite{brion2}
that a classification-free proof of \re{old-brilliant} for $\hat\g$, and hence
the irreducibility of $\ntg$, can be derived
from the fact that $\N$ is a complete intersection having only
rational singularities, see below. The advantage of the Takiff algebra case is that
the rather mysterious term
$\dim \text{span}\{F_1(x),\ldots,F_m(x)\}$ is being interpreted as
the rank of the differential of the quotient
map $\pi_G:\g \to \g\md G$ at $x$.

On the other hand, we will prove here the counterpart of \re{brilliant}
for $\hat\g$, using the classification. Brion's idea
cannot be applied directly to obtain a case-free proof of that stronger result.
The reason for being interested in
proving a counterpart of ~\re{brilliant} for $\hat\g$ is that we deduce from this the 
equidimensionality of some other quotient morphisms, see 
Theorems~\ref{takif-max-rang},\ref{not-vain2}.

We work in the setting of case {\bf B)} from Section~\ref{null-cone}.
\begin{rem}{Definition}  \label{extrem-good}
Let $\rho: G\to GL(V)$ be a representation of a connected
reductive group $G$. Then $V$ or $\rho$ is said to be {\it extremely good\/}
if 
\begin{itemize}
\item[ (1)] \ $\bbk[V]^G$ is a polynomial algebra;

\item[ (2)] \ $\max\dim_{x\in V} G{\cdot}x=\dim V-\dim  V\md G$;

\item[ (3)] \ If\/ $\pi_G: V\to  V\md G$ is the quotient morphism, then
$\Omega:=\{x\in V \mid (\text{d}\pi_G)_x\ \text{ is onto }\}$ is a big open subset of
$V$;
\item[ (4)] \ $\nv:=\pi_G^{-1}(\pi_G(0))$ consists of finitely many $G$-orbits. 
\item[ (5)] \ ${\nv}$ is irreducible and has only rational
singularities; 
\end{itemize}
\end{rem}%
Note that properties (1)--(3) are those appearing in Theorem~\ref{main-twice}.
Recall from Section~\ref{takif-2} that if $G$ is semisimple, then 
(2) and (3) are always satisfied.

\begin{s}{Theorem} \label{10.2}
Let\/ $V$ be an extremely good $G$-module and
$\hat V=V\ltimes V$ the corresponding $\hat G$-module. 
Then\/ 
\begin{itemize}
\item[\sf (i)] \ ${\goth N}^{\hat G}(\hat V)=\ntv$ 
is an irreducible complete intersection; 
\item[\sf (ii)] \ the ideal
of\/ $\ntv$ in $\bbk[\hat V]$ is generated by the basic invariants
in\/ $\bbk[\hat V]^{\hat G}$;
\item[\sf (iii)] \ $\pi_{\hat G}: \hat V\to \hat V\md \hat G$ is equidimensional
and\/ $\bbk[\hat V]$ is a free $\bbk[\hat V]^{\hat G}$-module.
\end{itemize}
\end{s}\begin{proof}
Let $f_1,\ldots,f_m$ be algebraically independent generators of $\bbk[V]^G$.
By Theorem~\ref{main-twice}, $\bbk[\hat V]^{\hat G}$ is freely generated
by the polynomials $f_1,\ldots,f_m,\wF_{f_1},\ldots,\wF_{f_m}$.
Recall from Section~\ref{null-cone} the stratification of the null-cone:
\[
   \fY_i(\nv)=\{v\in\nv\mid \rk(\textsl{d}\pi_G)_v=i\}, \quad i=0,1,\ldots,m \  .
\]
Since $\nv$ contains finitely many $G$-orbits, $\pi_G$ is equidimensional.
If $G{\cdot}x$ is the dense $G$-orbit in $\nv$, then $\dim G{\cdot}x=\dim V-m$
and therefore $x\in \fY_m(\nv)$
\cite[Korollar\,2]{knop}. (Corollary~2 is stated in Knop's article under the 
assumption that $G$ is semisimple. However, that proof works also for reductive
groups as long as conditions~(2) and (3) are satisfied.)
Since $\nv$ is irreducible and $\fY_m(\nv)\ne \varnothing$, it is a complete intersection.
The condition of the
irreducibility of $\ntv$ (Theorem~\ref{cond-EQ-b}(iii)) can be written as
\begin{equation}   \label{old-brilV}
  \dim V-\dim G{\cdot}v+ \rk(\textsl{d}\pi_G)_v > 2m  
\ \ \text{ if } \ v\not\in \fY_m(\nv) \ .
\end{equation}
We derive this inequality from a property of the local ring of (the
closure of) the orbit $G{\cdot}v \subset \nv$. Let ${\goth O}$ be this local
ring. Then $\dim{\goth O}=\dim\nv-\dim G{\cdot}v$ and
$\ed {\goth O}= \dim T_v\nv -\dim G{\cdot}v=\dim \g_v-\rk(\textsl{d}\pi_G)_v$.
Here $\ed{\goth O}$ is the embedding dimension of ${\goth O}$ 
and $T_v \nv$ is the tangent space of $\nv$ at $v$.
Since $\nv$ has only rational singularities, so has ${\goth O}$.
By a result of Goto-Watanabe (see \cite[Theorem\,2']{nw}), 
if a local ring ${\goth O}$ is a complete intersection with only rational
singularities and $\dim{\goth O}>0$, then 
$\ed{\goth O} < 2\dim {\goth O}$. Using the above expressions for
$\ed{\goth O}$ and $\dim{\goth O}$, one obtains  
inequality~\re{old-brilV}, and thereby the irreducibility of $\ntv$.

All other statements of the theorem are consequences of the fact that
$\ntv$ is irreducible.
By Theorem~\ref{cond-EQ-b}(ii), $\pi_{\hat G}$ is equidimensional.
If $v\in \fY_m(\nv)$, then the differentials of
the generators $f_1,\ldots,f_m,\wF_{f_1},\ldots,\wF_{f_m}$
are linearly independent
at $(v,0)\in \ntv \subset\hat V$. This fact and the irreducibility of 
$\ntv$ imply that $\ntv$ is a complete intersection whose ideal is
generated by the polynomials $f_1,\ldots,f_m,\wF_{f_1},\ldots,\wF_{f_m}$
(cf. \cite[Prop.\,6]{ko63}).
\end{proof}%
{\bf Remark.} The most subtle point in the definition of extremely good representations
is the rationality of singularities of $\nv$.
For the adjoint representations, this result is due to W.\,Hesselink~\cite{hess2}.
The idea to exploit the fact that $\N={\goth N}(\g)$ is a complete
intersection with only rational singularities, and
to use the Goto-Watanabe inequality for local rings is due to M.\,Brion \cite{brion2}.
Since $(G, \Ad)$ is extremely good, this approach yields
a conceptual proof of \cite[Theorem\,2.4]{geof1}.

\begin{s}{Corollary}
If\/ $V$ is extremely good, then the closure of
the image of the map 
\[
\varkappa: V\times \g\to V\times V, \quad (v,x)\mapsto (v, x{\cdot}v), 
\]
is a factorial complete intersection of codimension $m=\dim V\md G$ and the
ideal of\/ $\ov{\Ima\varkappa}$ is generated by $\wF_{f_1},\ldots,\wF_{f_m}$.
\end{s}\begin{proof}
This follows from the irreducibility of $\ntv$ and Theorem~\ref{factorialB}.
\end{proof}%
Since conditions~(4) and (5) are quite restrictive,
there are only a few extremely good representations. 
Below is a list of such irreducible representations known to this author
such that $G$ is simple and $\bbk[V]^G\ne\bbk$, except the adjoint ones:

$(\GR{B}{n} \text{ or } \GR{D}{n}, \vp_1),\,(\GR{B}{3}, \vp_3),\,(\GR{B}{4}, 
\vp_4),\,(\GR{G}{2}, \vp_1),\,(\GR{A}{n}, 2\vp_1),\,(\GR{A}{2n-1}, \vp_2)
,\,(\GR{E}{6}, \vp_1),$

$(\GR{C}{3}, \vp_3),\,(\GR{A}{5}, \vp_3),\,(\GR{D}{6}, \vp_6),\,(\GR{E}{7}, \vp_1)
,\,(\GR{B}{5}, \vp_5),\,(\GR{F}{4}, \vp_1),\,(\GR{C}{n}, \vp_2)$.
\\[.6ex]
The representations are given by their highest weights, and $\{\vp_i\}$ are fundamental
weights of $G$ with numbering from \cite{t41}. 
For all representations in the list but the last one, the algebra of covariants, 
$\bbk[V]^U$, is polynomial \cite{br83} (here $U$ is a maximal unipotent
subgroup of $G$). Therefore the same is true for $\bbk[\nv]^U$.
Then a result of Kraft (see \cite[1.5-6]{br-these}) shows that $\nv$ has 
rational singularities.
\\
I conjecture that if $G$ is simple and $V$ is a simple $G$-module, 
then $V$ is extremely good if and only if $\dim V\le \dim G$. Practically, this means that 
one has to only verify that $\nv$ has rational singularities for the 
following representations:\\
$(\GR{A}{6}, \vp_3),\,(\GR{A}{7}, \vp_3),\,(\GR{B}{6}, \vp_6),\,(\GR{D}{7},\vp_7)$.
\\[.7ex]
For $V=\g$, inequality~\re{old-brilV} reads
\begin{equation}  \label{old-bril-t}
    \dim\z_\g(x)+\rk(\textsl{d}\pi_G)_x > 2\rk\g=2m
\ \ \text{ if } \ x\not\in  \fX_m(\N)=\N^{reg} \ .
\end{equation} 
This inequality was proved in \cite[2.6-2.15]{geof1} in a case-by-case fashion.
Below, we prove a stronger result, which is the counterpart of
inequality~\re{brilliant} in the context of Takiff algebras.
By the Morozov-Jacobson theorem \cite[Ch.\,3, Theorem\,1.3]{t41}, any $x\in\N\setminus \{0\}$ 
can be embedded
in an $\tri$-triple $\{x,h,y\}$, where $h$ is semisimple; 
$x$ is said to be {\it even\/} if the $\ad h$-eigenvalues in $\g$ are even.
Following E.B.\,Dynkin, $h$ is called a {\it characteristic\/} of $x$.

\begin{s}{Theorem}                \label{thm-bril-takif}
Let $\g$ be a simple Lie algebra and $x\in\N$. Then
\begin{equation}   \label{bril-takif}
   \mathsf L:=\dim\z_\g(x)+2\rk(\text{d}\pi_G)_x- 3\rk\g\ge 0 \ .
\end{equation}
If $\g=\sln$, then\/ $\mathsf L=0$ if and only if
the matrix $x$ has at most two Jordan blocks.
Furthermore, if $\g\ne {\frak sl}_{2n+1}$, then $\mathsf L=0$ if and only if 
$x$ is even and $[\z_\g(h),\z_\g(h)]$
is a sum of several copies of $\tri$. \textrm{\rm (Here $h$ is a characteristic 
of $x$)}. 
\end{s}\begin{proof}
The proof is case-by-case.  However, the computations themselves are much 
shorter and more transparent than those in \cite{geof1}, because
our inequality is stronger, and we use formulae for 
$\dim\z_\g(x)$ in terms of dual partitions
(already used for $Sp$ and $SO$ in the proof of
Propositions~\ref{Sp-case} and \ref{SO-case}).

For the classical series, we work with the partition of $x$; while for the exceptional 
algebras the explicit classification of nilpotent orbits is used. 
If $\g=\g(\BV)$ is classical and $x\in \g(\BV)$ is nilpotent, then
$\boldsymbol{\eta}=(\eta_1,\eta_2,\ldots)$ is the partition of $\dim \BV$  corresponding to 
$x$ and $(\hat\eta_1,\ldots,\hat\eta_s)$ is the dual partition.
Here $s=\eta_1$. For $Sp$ and $SO$, 
%%the symplectic and orthogonal series, 
our analysis is quite similar to those in Propositions~\ref{Sp-case} and \ref{SO-case}.

(A) \ $\g={\frak sl}(\BV)$, $\dim\BV=n+1$. \\[.7ex]
Here $\dim\z_\g(x)=\sum_{i=1}^s\hat\eta_i^2 - 1$ and 
$\rk(\textsl{d}\pi_G)_x=\eta_1-1=s-1$ \cite[Theorem\,4.2.1]{ri2}. Then 
%%\begin{multline*}
\[
\mathsf L=\sum_{i=1}^s\hat\eta_i^2 - 1+2(s-1)-3n= 
\sum_{i=1}^s\hat\eta_i^2-3\sum_{i=1}^s\hat\eta_i+2s=
\sum_{i=1}^s(\hat\eta_i-1)(\hat\eta_i-2)\ge 0 .
%%\end{multline*}
\]
This expression equals zero if and only if all $\hat\eta_i \le 2$, i.e.,
$x$ has at most two Jordan blocks.

(B) \ $\g={\frak so}(\BV)$, $\dim\BV=2n+1$. \\[.7ex]
Here $\hat\eta_1$ is odd, $\dim\z_\g(x)= 
\frac{1}{2}\bigl(\sum_i\hat\eta_i^2 -\#\{j\mid \eta_j \text{ is odd}\}\bigr)$, and 
$\rk(\textsl{d}\pi_G)_x=[s/2]$ \cite[Theorem\,4.3.3]{ri2}. Then
\begin{multline*}
 \mathsf L=\frac{1}{2}\bigl(\sum_{i=1}^s\hat\eta_i^2 -\#\{j\mid \eta_j \text{ is odd}\}\bigr)
+2[s/2]-3n= \\
\frac{1}{2}\bigl(\sum_{i=1}^s (\hat\eta_i-1)(\hat\eta_i-2)+3-
\#\{j\mid \eta_j \text{ is odd}\}-2(s-2[s/2])\bigr) \ .
\end{multline*}
If $\hat\eta_1=1$, then $\mathsf L=0$. This is the case of regular nilpotent
elements.
\\[.6ex]
If $\hat\eta_1=3$, then $\sum_{i=1}^s (\hat\eta_i-1)(\hat\eta_i-2)=
2\eta_3 \ge 2$. Therefore 
$2\mathsf L=2\eta_2+3-\#\{j\mid \eta_j \text{ is odd}\}-2(s-2[s/2])$. Since
$\#\{j\mid \eta_j \text{ is odd}\}\le 3$ and $2(s-2[s/2])\le 2$, $2\mathsf L$ is
nonnegative.
Furthermore, $\mathsf L=0$ if and only if $\eta_3=1$ and all $\eta_i$'s are odd.
\\[.6ex]
If $\hat\eta_1\ge 5$, then $\sum_{i=1}^s(\hat\eta_i-1)(\hat\eta_i-2) \ge \hat\eta_1+2
\ge \#\{j\mid \eta_j \text{ is odd}\}+2$. 
Next, \\ 
$3-2(s-2[s/2]) \ge 0$. Hence $\mathsf L$ is positive.

(C) \ $\g={\frak sp}(\BV)$, $\dim\BV=2n$. \\[.7ex]
Here $\dim\z_\g(x)= 
\frac{1}{2}\bigl(\sum_i\hat\eta_i^2 +\#\{j\mid \eta_j \text{ is odd}\}\bigr)$ and 
$\rk(\textsl{d}\pi_G)_x=[s/2]$ \cite[Theorem\,4.3.3]{ri2}. Then 
\begin{multline*}
  \mathsf L=\frac{1}{2}\bigl(\sum_{i=1}^s\hat\eta_i^2 +\#\{j\mid \eta_j \text{ is odd}\}\bigr)
+2[s/2]-3n= \\
\frac{1}{2}\bigl(\sum_{i=1}^s (\hat\eta_i-1)(\hat\eta_i-2)+
\#\{j\mid \eta_j \text{ is odd}\}-2(s-2[s/2])\bigr) \ .
\end{multline*}
It is easily seen that $\mathsf L=0$ if and only if $\hat\eta_1 \le 2$.
Otherwise  it is positive.

(D) \ $\g={\frak so}(\BV)$, $\dim\BV=2n$. \\[.7ex]
Here  $\hat\eta_1$ is even and
$\dim\z_\g(x)$ is as in (B).
%%= \frac{1}{2}\bigl(\sum_i\hat\eta_i^2 -
%%\#\{j\mid \eta_j \text{ is odd}\}\bigr)$ and 
For the rank of $\textsl{d}\pi_G$, we have \cite[Theorem\,4.4.2]{ri2}\\
$
\rk(\textsl{d}\pi_G)_x=\left\{\begin{array}{rl} [s/2], & \text{ if } \hat\eta_1\ge 4; \\
(2n-i+1)/2, & \text{ if } \boldsymbol{\eta}=(2n-i,i) \text{ with $i$ odd};\\
 l,  & \text{ if } \boldsymbol{\eta}=(n,n) \text{ and } n=2l \ .
\end{array}\right.
$
\\[.6ex]
Then
\begin{multline*}
\mathsf L= \frac{1}{2}\bigl(\sum_{i=1}^s\hat\eta_i^2 -\#\{j\mid \eta_j \text{ is odd}\}\bigr)
+ 2\rk(\textsl{d}\pi_G)_x- 3n= \\
\frac{1}{2}\bigl(\sum_{i=1}^s (\hat\eta_i-1)(\hat\eta_i-2)+4\rk(\textsl{d}\pi_G)_x
-2s-\#\{j\mid \eta_j \text{ is odd}\}\bigr) \ .
\end{multline*}
Now, a consideration of cases shows that $\mathsf L=0$ if 
$\hat\eta_1=2$. If $\hat\eta_1\ge 4$, then $\mathsf L> 0$ unless
$\boldsymbol{\eta}=(\eta_1,\eta_2,1,1)$, where $\eta_1,\eta_2$ are both odd.

(EFG) \ $\g$ is exceptional. \\[.7ex]
It is enough to check inequality~\re{bril-takif}
for sufficiently large orbits (with $\dim\z_\g(x) \le 3\rk\g$).
To this end, one can consult the tables in \cite[Ch.\,8]{CM} for dimensions
of orbits and \cite[Appendix]{ri2} for the values of $\rk(\textsl{d}\pi_G)_x$.
Below we list all non-regular nilpotent orbits with $\mathsf L=0$.
The orbits 
are represented by their Dynkin-Bala-Carter labels.

\begin{itemize} 
\item[$\GR{G}{2}$:] \quad $\GR{G}{2}(a_1)$;
\item[$\GR{F}{4}$:] \quad $\GR{F}{4}(a_1)$, $\GR{F}{4}(a_2)$;
\item[$\GR{E}{6}$:] \quad $\GR{E}{6}(a_1)$, $\GR{D}{5}$, $\GR{E}{6}(a_3)$;
\item[$\GR{E}{7}$:] \quad $\GR{E}{7}(a_1)$, $\GR{E}{7}(a_2)$, $\GR{E}{6}$, $\GR{E}{6}(a_1)$;
\item[$\GR{E}{8}$:] \quad $\GR{E}{8}(a_1)$, $\GR{E}{8}(a_2)$, $\GR{E}{8}(a_3)$, $\GR{E}{8}(a_4)$. 
%%$\GR{E}{6}$, $\GR{E}{6}(a_1)$;
\end{itemize}
Inspecting the tables in \cite[Ch.\,8]{CM}
shows that these are precisely the even nilpotent orbits whose weighted
Dynkin diagrams have no adjacent zeros, which exactly means that $[\z_\g(h),\z_\g(h)]$ is
a sum of several $\tri$'s.
\\[.6ex]
For $\g$ classical, there is a rule for writing out the characteristic $h$
in terms of $\boldsymbol{\eta}$ \cite[Ch.\,IV]{ss}.
Hence the Levi subalgebra $\z_\g(h)$ can be computed.
This yields the last assertion of the theorem.
\end{proof}%
A geometric meaning of \re{bril-takif} will be made clear in the following result.
Let $\g=\g_0\oplus\g_1$ be a $\BZ_2$-grading and $\vartheta$ the corresponding 
involutory automorphism of $\g$.
Then $\vartheta$ can be  extended to an involution of the Takiff algebra 
$\hat\g$ by letting $\vartheta(x+\eps y)=\vartheta(x)+\eps\vartheta(y)$.
The corresponding eigenspaces are $\hat\g_0=\g_0\ltimes\g_0$ and
$\hat\g_1=\g_1\ltimes\g_1$.
Here $\hat\g_1$ is a $\hat\g_0$-module just in the sense of definition
given in Section~\ref{takif-2}.
The $G_0$-module $\g_1$ is not {\it extremely good},
so that Theorem~\ref{10.2} cannot be applied. But
it is `good enough' in the sense that it satisfies properties (1),\,(2),\,(4)
of Definition~\ref{extrem-good}.

%Let $\g=\g_0\oplus\g_1$ be a $\BZ_2$-grading of {\it maximal rank}.
%This means that $\g_1$ contains a Cartan subalgebra of $\g$. 
%(See \cite{an} about involutions of maximal rank.)

\begin{s}{Theorem}    \label{takif-max-rang}
Suppose\/ $\g=\g_0\oplus\g_1$ is a $\BZ_2$-grading of maximal rank, i.e.,
$\g_1$ contains a Cartan subalgebra of\/ $\g$.
Then the quotient morphism\/ $\hat\pi :\hat\g_1 \to \hat\g_1\md \hat G_0$
is equidimensional.
\end{s}\begin{proof}
%%The null-cone ${\goth N}(\g_1)$ consists of finitely many $G_0$-orbits.
Recall the relationship between orbits and null-cones for
the actions $(G:\g)$ and $(G_0:\g_1)$.
The null-cones are $\N$ and ${\goth N}(\g_1)$, respectively.
\begin{itemize} 
\item[$\bullet$] \ ${\goth N}(\g_1)=\N\cap \g_1$;
\item[$\bullet$] \ $G{\cdot}x\cap \g_1$ is a union of finitely many $G_0$-orbits;
\item[$\bullet$] \ If $x\in\g_1$, then $\dim G_0{\cdot}x=\frac{1}{2}\dim G{\cdot}x$; 
\item[$\bullet$] \ For any $x\in\g$, we have $G{\cdot}x\cap\g_1\ne \varnothing$;
\item[$\bullet$] \ $\bbk[\g]^G\simeq \bbk[\g_1]^{G_0}$.
\end{itemize}
The first three properties hold for all $\BZ_2$-gradings, whereas the last two are
characteristic for the involutions of maximal rank, see~\cite{an}.

Let us see what the equidimensionality criterion (Theorem~\ref{cond-EQ-b}(i)) means here.
We have $V=\g_1$, $G=G_0$, and $m=\dim\g_1\md G_0$.
Since ${\goth N}(\g_1)$ consists of finitely many $G_0$-orbits, that criterion
reads
\[
   \dim{\goth N}(\g_1)-\dim G_0{\cdot}x\ge \dim\g_1\md G_0- \rk(\textsl{d}\pi_{G_0})_x
\]
for any $x\in {\goth N}(\g_1)$. Here $\pi_{G_0}: \g_1\to \g_1\md G_0$
is the quotient morphism.
In view of the above properties of such $\BZ_2$-gradings, we have
$\dim{\goth N}(\g_1)=\frac{1}{2}\dim\N=\frac{1}{2}(\dim\g-\rk\g)$,
$\dim\g_1\md G_0=\rk\g$, and $\rk(\textsl{d}\pi_{G_0})_x=\rk(\textsl{d}\pi_{G})_x$.
The latter stems from the isomorphism $\bbk[\g]^G\simeq \bbk[\g_1]^{G_0}$.
Rewriting the previous inequality using this data yields precisely 
inequality~\re{bril-takif}\,! Thus, the fact that $\hat\pi$ is
equidimensional is essentially equivalent to Theorem~\ref{thm-bril-takif}.
\end{proof}%
Yet another geometric application of Eq.~\re{bril-takif} is 
the following (cf. Theorem~\ref{not-vain}):

\begin{s}{Theorem}    \label{not-vain2}
Set ${\g}^{[n]}=\g\ltimes\g^{\oplus n}$, where $\g^{\oplus n}$ 
(the sum of $n$ copies) regarded
as commutative Lie algebra and $n\ge 1$.
Consider the adjoint action $({G}^{[n]}:{\g}^{[n]})$. Then $\pi_{{G}^{[n]}}$ is
equidimensional if and only if $n\le 2$.
\end{s}\begin{proof} For $n=1$, the assertion is already proved.
Next, $\dim(\g^{\oplus n})^T=n\rk\g$ and for $x\in \N$ the 
equidimensionality condition
of Theorem~\ref{cond-EQ}(i) reads $\dim\z_\g(x)-\rk\g\ge n(\rk\g-
\rk (\textsl{d}\pi_G)_x)$, which is exactly \re{bril-takif} for $n=2$.
Conversely, if $n\ge 3$, then this condition is not satisfied for the subregular 
nilpotent orbit.
\end{proof}%
{\bf Remark.} In the last theorem, the null-cone ${\goth N}(\g^{[2]})$
is always reducible. Indeed, each nilpotent $G$-orbit such that
$\mathsf{L}=0$ in \re{bril-takif} gives rise to an irreducible component
of ${\goth N}(\g^{[2]})$, see Remark~\ref{remarks-8}(1). 
The proof of Theorem~\ref{thm-bril-takif}
shows that, for any $\g$, there are at least two orbits with 
$\mathsf{L}=0$.
\\[.6ex]
There are several equivalent ways to present inequality~\re{bril-takif}. 
Let $\mathcal B$ denote the variety of Borel subgroups of $G$. 
For any $x\in\N$, set
$\mathcal B_x=\{B'\in \mathcal B\mid x\in \Lie B'\}$.
Recall that $\fX_i=\fX_{i,\g}=\{x\in\g\mid \rk(\textsl{d}\pi_G)_x=i\}$ and $\fX_{i,\g}(\N)=
\fX_{i,\g}\cap\N$. This stratification is determined by the covariants of type $\g$.

\begin{s}{Proposition}  \label{equiv-ways}
Let $\g$ be a simple Lie algebra and
%%For any $x\in {\goth N}(\g_1)$, 
$m=\rk\g$. Then the following holds:
%% following conditions are equivalent:

(1) \ $\codim_\N\fX_{i,\g}(\N) \ge 2(m-i)$ for any $i=0,1,\ldots,m$;

(2) \ $\dim\mathcal B_x +\rk(\textsl{d}\pi_G)_x \ge \rk\g$ for any $x\in\N$;

%%(4) \ $\hat\pi :\g_1\ltimes\g_1 \to \g_1\ltimes\g_1\md \hat G_0$
%%is equidimensional;

(3) \ If $\co$ is the local ring of any $G$-orbit in $\N$, then
$\ed\co\le \frac{3}{2}\dim\co$;

(4) \ If $\g=\g_0\oplus\g_1$ is a $\BZ_2$-grading of maximal rank and 
$x\in{\goth N}(\g_1)$, %%and $\pi_0: \g_1\to \g_1\md G_0$, then
then $\dim (G_0)_x+ \rk(\textsl{d}\pi_{G_0})_x \ge \rk\g$.

(5) \ If $\g=\g_0\oplus\g_1$ is a $\BZ_2$-grading of maximal rank and 
$\co'$ is the local ring of a $G_0$-orbit in ${\goth N}(\g_1)$, then
$\ed\co'\le {2}\dim\co'$.
\end{s}\begin{proof} 
In fact, all these conditions are equivalent to
inequality~\re{bril-takif}.
Since $\N$ contains finitely many $G$-orbits,
(1) can be written as $\codim_\N G{\cdot}x \ge 2(m-\rk(\text{d}\pi_G)_x)$
for any $x\in\N$, which makes it clear that (1) is equivalent to \re{bril-takif}.
For (2), one should use the fact that 
$\dim\mathcal B_x=\frac{1}{2}(\dim\z_\g(x)-\rk\g)$, see e.g.
\cite[4.3.10, 4.5]{spr}.
For (3), one have to use formulae for $\dim\co$ and $\ed\co$
written out in the proof of Theorem~\ref{10.2}.
For (4), we notice that since $\dim\g_1-\dim\g_0=\rk\g$, the equality 
$\dim G_0{\cdot}x=\frac{1}{2}\dim G{\cdot}x$ is equivalent to that 
$\dim (G_0)_x=\frac{1}{2}(\dim\z_\g(x)-\rk\g)=\dim\mathcal B_x$.
Finally, the inequalities in (4) and (5) are obtained from each other
via simple transformations.
\end{proof}%
{\bf Remark.} Concerning (5), we note that this inequality is weaker
than the Goto-Watanabe inequality from the proof of
Theorem~\ref{10.2}, but ${\goth N}(\g_1)$
is not normal and can be reducible.

\begin{s}{Corollary}  \label{3mal}
$\codim_\g\fX_{i,\g} \ge 3(m-i)$ for any $i=0,1,\ldots,m$.
\end{s}\begin{proof}
It follows from the definition of $\fX_{i,\g}$ that $\dim \ov{\fX}_{i,\g}\md G=i$. 
Since $\ov{\fX}_{i,\g}$ is conical, the fibre of the origin of the 
morphism $\ov{\fX}_{i,\g} \to \ov{\fX}_{i,\g}\md G$ has the maximal
dimension, i.e., $\dim\fX_{i,\g} \le i +\dim \fX_{i,\g}(\N)$, which is exactly what we want,
in view of Proposition~\ref{equiv-ways}(1).
\end{proof}%
There are many open problems and observations related to our results on 
reductive Takiff algebras and $\BZ_2$-contractions. Here are some of them.

$1^o$. \ It seems that if $H$ is a spherical reductive subgroup of $G$ and 
$\ka=\h\ltimes\ma$ is the corresponding isotropy contraction of $\g$, then
$\pi_K$ is always equidimensional. At least, I have verified this in case
$G$ is simple. In fact, Examples~\ref{D4G2} and \ref{B3G2} present several instances
of this verification.

$2^o$. \ It would be quite interesting to have a case-free proof for 
Theorem~\ref{thm-bril-takif} or, equivalently, \ref{takif-max-rang}.
Various equivalent forms of that result presented in Proposition~\ref{equiv-ways}
suggest that there might be different approaches to proving it.
From the geometric point of view, the equidimensionality of
$\hat\pi$ means that there exists a transversal subspace to 
${\goth N}(\hat\g_1)$, i.e., a subspace $U$ such that $\dim U=\dim\hat\g_1\md\hat G_0$
and $U\cap {\goth N}(\hat\g_1)=\{0\}$.

$3^o$. \ Whenever some quotient morphism is equidimensional,
it is interesting to find a natural transversal subspace to the null-cone.
One may ask for such a subspace in the setting  of Theorems~\ref{main-free},
\ref{10.2}, \ref{takif-max-rang}. 
Even for the adjoint representation of $\hat\g=\g\ltimes\g$
it is not known how to naturally construct a transversal space to ${\goth N}(\hat\g)$.
If $\Delta_\te\subset\g\ltimes\g$ is the diagonally embedded Cartan subalgebra,
then $\Delta_\te \cap {\goth N}(\hat\g)=\{0\}$, so that one has
a "one-half" of a transversal space. The problem is to construct
the second half.
Similarly, if $\ka$ is a $\BZ_2$-contraction of a simple Lie algebra,
it is not known how to construct a transversal space to ${\goth N}^u(\ka)$.

$4^o$. \ If one knows that some null-cone ${\goth N}$ is irreducible, then
it is tempting to find a resolution of singularities for ${\goth N}$.

%$5^o$. \ We did not say much on regular invariants of the
%{\sl coadjoint\/} representation of 
%a $\BZ_2$-contraction $\ka=\g_0\ltimes\g_1$. For $\BZ_2$-gradings of maximal rank,
%the action $(G_0:\g_1)$ has trivial generic stabiliser. Therefore,
%in view of Theorem~\ref{main-dual},
%$\bbk[\ka^*]^K\simeq \bbk[\g_1]^{G_0}$ is polynomial. But we do not know what happens
%for the other $\BZ_2$-gradings.

$5^o$. \ A case-by-case verification shows that
$\fX_{1,\g}(\N)$ is irreducible for any simple $\g$,
and the dense $G$-orbit in it is Richardson.

                %%%%
%%%%%%%%%%%%%%%%
%%%%%%%%%%%%%%%%   Section 11
%%%%%%%%%%%%%%%%
                %%%%

\section{On invariants and null-cones for generalised Takiff Lie algebras}  
\label{generalised}
\setcounter{equation}{0}

\noindent
Following \cite{rt}, we recall the definition of a generalised Takiff Lie algebra.
The infinite-dimensional $\bbk$-vector space $\q_\infty:=\q\otimes \bbk[\mathsf T]$ has a natural
structure of a Lie algebra such that $[x\otimes \mathsf T^l, y\otimes \mathsf T^k]=
[x,y]\otimes \mathsf T^{l+k}$.
Then $\q_{\ge (n+1)}=\displaystyle\bigoplus_{j\ge n+1} \q\otimes \mathsf T^j$ is an ideal of
$\q_\infty$, and the 
respective quotient
is a {\it generalised Takiff Lie algebra}, denoted $\q\langle n\rangle$.
We also say that $\q\langle n\rangle$ is the $n$-{\it th Takiff algebra}.
Write $Q\langle n\rangle$ for the corresponding connected group.
Clearly, $\dim\q\langle n\rangle=(n+1)\dim\q$ 
and  $\q\langle 1\rangle\simeq \q\ltimes \q$.
The main results of \cite{rt} are the following:

(i) \ $\ind\q\langle n\rangle=(n+1)\ind\q$, 

(ii) if $\q=\g$ is semisimple, then $\bbk[\g\langle n\rangle]^{G\langle n\rangle}$ is
a polynomial algebra whose set of basic invariants is explicitly described.
\\[.7ex]
Actually, the authors of \cite{rt} work with invariants of 
the coadjoint representation of $G\langle n\rangle$, but this makes no difference,
since $\g\langle n\rangle$ is quadratic.

In this section, we generalise the results from (ii) in the spirit of
Section~\ref{takif-2}. Let $\q\langle n\rangle^u$ denote the image of
$\q_{{\ge} 1}$ in $\q\langle n\rangle$. It is a nilpotent Lie algebra, 
which is noncommutative for $n\ge 2$,
and $\q\langle n\rangle\simeq \q \ltimes \q\langle n\rangle^u$.
Accordingly, one obtains the semi-direct product structure of the group:
$Q\langle n\rangle=Q\ltimes Q\langle n\rangle^u$.

\begin{s}{Theorem}  \label{gen-takiff}
Suppose $\q$ satisfies conditions %%(1)--(3) of Theorem~\ref{main-twice}.

(1) \ $\bbk[\q]^Q$ is a polynomial algebra;

(2) \ $\max\dim_{x\in \q} Q{\cdot}x=\dim \q-\dim \q\md Q$;

(3) \ If\/ $\pi_Q: \q\to \q\md Q$ is the quotient morphism and\/
$\Omega:=\{x\in \q \mid (\text{d}\pi_Q)_x\ \text{ is onto }\}$, then\/ 
$\q\setminus\Omega$ contains no divisors.
\\
Then 
\begin{itemize}
\item[\sf (i)] \  $\bbk[\q\langle n\rangle]^{Q\langle n\rangle^u}$ is a polynomial 
algebra of Krull dimension\/ $\dim\q+n\dim\q\md Q$ whose algebraically independent 
generators can explicitly be described;
\item[\sf (ii)] \ $\bbk[\q\langle n\rangle]^{Q\langle n\rangle}$ is a polynomial 
algebra of Krull dimension\/ $(n+1)\dim\q\md Q$ whose algebraically independent 
generators can explicitly be described;
%%\item[\sf (iii)] \ The algebra $\hat\q$ satisfies conditions (1)--(3), too.
\end{itemize}
\end{s}\begin{proof*}
Let $\boldsymbol{x}=x_0{+}\eps x_1{+}\ldots {+}\eps^n x_n$ denote the image of 
$\sum_{i=0}^n x_i\otimes \mathsf T^i$ in $\q\langle n\rangle$. Here each $x_i\in\q$ and
$\eps$ is regarded as the image of $\mathsf T$ in $\bkk[\mathsf T]/(\mathsf T^{n+1})$
Set $m=\dim\q\md Q$, and let $f_1,\ldots,f_m$
be a set of basic invariants in $\bbk[\q]^Q$.
Expand the polynomial $f_i(x_0+\eps x_1+\ldots +\eps^n x_n)$ using the relation
$\eps^{n+1}=0$. We obtain
\[
 f_i(x_0+\eps x_1+\ldots +\eps^n x_n)=
\sum_{j=0}^n \eps^j \wF_{i}^{(j)}(x_0,x_1,\ldots, x_n) \ .
\]
Following the argument in \cite[Sect.\,III]{rt}, one proves that
\begin{gather} 
\text{$\wF_{i}^{(j)}$ depends only on $x_0,\ldots,x_j$ 
and} \notag \\   \label{gen-Fij}
\wF_{i}^{(j)}(x_0,\ldots,x_j)=\langle (\textsl{d}f_i)_{x_0},x_j\rangle+
p_{ij}(x_0,\ldots,x_{j-1}) .
\end{gather}
It follows from the construction that 
all $\wF_{i}^{(j)}$ belong to $\bbk[\q\langle n\rangle]^{Q\langle n\rangle}$.

(i) Making use of Lemma~\ref{igusa} and Eq.~\re{gen-Fij}, 
we prove that the polynomials $\wF_{i}^{(j)}$,
$i=1,\ldots,m$, $j=1,\ldots,n$, and the coordinates on the first factor in 
$\q\langle n\rangle$ freely generate $\bbk[\q\langle n\rangle]^{Q\langle n\rangle^u}$.

Consider the mapping 
\[
   \psi: \q\langle n\rangle \to \q\times \bbk^{nm}\ ,
\]
given by $\psi(\boldsymbol{x})=
(x_0, \wF_{1}^{(1)}(\boldsymbol{x}),\ldots,\wF_{m}^{(n)}(\boldsymbol{x}))$.
Here we regard $\q$ as $\q\langle n\rangle/\q\langle n\rangle^u$, so that
$\q\times \bbk^{nm}$ is a variety with trivial $Q\langle n\rangle^u$-action.
%%We identify $\bbk^{n+m}$ with $\q\times \bbk^{m}$.
If $x_0\in\Omega$, then $(\textsl{d}f_i)_{x_0}$ are linearly independent. Therefore
Eq.~\re{gen-Fij} shows that
the system $\wF_i^{(j)}(x_0{+}\eps y_1{+}\ldots{+}\eps^n y_n)=\ap_i^{(j)}$, $i=1,\ldots,m$, $j=1,\ldots,n$
has a solution, say $(y_1,\ldots,y_n)$, for any $nm$-tuple
$\boldsymbol{\ap}=(\ap_1^{(1)},\ldots,\ap_m^{(n)})$.
Indeed, $(y_1,\ldots,y_n)$ can be computed consecutively: First $y_1$, then $y_2$,
and so on.
Hence $\Ima\psi\supset \Omega\times \bbk^{nm}$, i.e.,
$\Ima\psi$ contains a big open subset of $\q\times\bbk^{nm}$.
This also implies that the coordinates on $\q$ and the polynomials $\wF_i^{(j)}$
are algebraically independent. It follows that
%%$\displaystyle
\[
\max_{\boldsymbol{x}\in \q\langle n\rangle}
\dim Q\langle n\rangle^u{\cdot}\boldsymbol{x}\le
\dim \q\langle n\rangle-\dim\q-mn=n(\dim\q-m) \ .
\]
Next, consider $\Omega'=\Omega\cap\{z\in \q\mid \dim Q{\cdot}z=\dim\q-m\}$. In view of
condition~(2), it is still a non-empty open 
$Q$-stable subset of $\q$. Fix $x_0\in\Omega'$,  and let 
$(\bar y_1,\ldots,\bar y_n)$ be a solution to the system
$\wF_i^{(j)}(x_0{+}\eps y_1{+}\ldots{+}\eps^n y_n)=
\ap_i^{(j)}$, $i=1,\ldots,m$, $j=1,\ldots,n$. 
%%That is, $\psi^{-1}(x,\boldsymbol{\ap})\ni x{+}\eps y_1{+}\ldots{+}\eps^n y_n$. 
Then  
$\psi^{-1}(x_0,\boldsymbol{\ap})\supset 
Q\langle n\rangle^u{\cdot}(x_0{+}\sum_{i=1}^n\eps^i \bar y_i)$.
Since $x_0\in\Omega$, we have $\dim \psi^{-1}(x_0,\boldsymbol{\ap})=
n(\dim\q-m)$. On the other hand, the following holds
\\[.8ex]
{\bf Claim.} \ {\sl If $x\in\q^{reg}$, then 
$\dim Q\langle n\rangle^u{\cdot}(x+\eps y_1{+}\ldots{+}\eps^n y_n)=
n\dim Q{\cdot}x=n(\dim\q-m)$ for any $(y_1,\ldots,y_n)\in \q^n$.}
\\[.8ex]
{\sl Proof of the claim.}
We argue by induction on $n$. For $n=1$, the assertion is obvious.
Assume that $n\ge 2$. Consider the $Q\langle n\rangle^u$-equivariant projection
\[
(x+\sum_{1}^n \eps^i y_i\in \q\langle n\rangle) \stackrel{p}{\mapsto}
(x+\sum_{1}^{n-1} \eps^i y_i\in \q\langle n{-}1\rangle)\ .
\]
Let $\co_n$ denote the $ Q\langle n\rangle^u$-orbit of
$x+\sum_{1}^{n} \eps^i y_i$. Then $p(\co_n)=\co_{n-1}$. By the induction
hypothesis, $\dim\co_{n-1}=(n-1)(\dim\q-m)$.
It is easily seen that
\[
   p^{-1}(x+\sum_{1}^{n-1} \eps^i y_i)\cap \co_n \supset
x+\sum_{1}^{n-1} \eps^i y_i+\eps^n(y_n +[\q,x]) \ .
\]
For, the right hand side is precisely an orbit of the subgroup
$\exp(\eps^n\q) \subset Q\langle n\rangle^u$.
Hence $\dim\co_n \ge n(\dim\q-m)$. But it is already proved 
that the dimension of every $Q\langle n\rangle^u$-orbit is at most
$n(\dim\q-m)$.
\qus%
\noindent
Hence 
$\psi^{-1}(x_0,\boldsymbol{\ap})=
Q\langle n\rangle^u{\cdot}(x_0{+}\eps \bar y_1{+}\ldots{+}\eps^n \bar y_n)$ 
for dimension reason.
Thus, a generic fibre of $\psi$ is an $Q\langle n\rangle^u$-orbit, and 
Lemma~\ref{igusa} applies here.

(ii) Follows from (i) and the description of $Q\langle n\rangle^u$-invariants.
\end{proof*}%
\begin{rem}{Remark} It was noticed in Section~\ref{z2} that any Takiff algebra $\q\ltimes\q$
is a $\BZ_2$-contraction of $\q\dotplus\q$. Similar phenomenon holds for the
generalised Takiff algebras: $\q\langle n\rangle$ is a contraction of
${\q\dotplus\ldots\dotplus\q}=(n+1)\q$. The starting point for constructing
such a contraction is to consider the action $\BZ_{n+1}$ on $(n+1)\q$ that 
cyclically permutes the summands. 
On the other hand, given $\q\langle n\rangle$, it can further be contracted
to $\q\ltimes \q^{\oplus n}$, the "usual" semi-direct product, where $\q^{\oplus n}$ 
is regarded as commutative Lie algebra. The details are left to the reader.
Thus, 
\[
 {\q\dotplus\ldots\dotplus\q}=(n+1)\q \ \leadsto \ \q\langle n\rangle 
\ \leadsto \ \q\ltimes \q^{\oplus n}
\]
is a chain of contractions.
\end{rem}%
Using Eq.~\re{bril-takif} and Eq.~\re{gen-Fij} one can easily prove that 
if $\g$ is semisimple and $\g\langle 2\rangle$ is the second Takiff
Lie algebra, then the quotient morphism 
$\pi_{G\langle 2\rangle}: \g\langle 2\rangle \to \g\langle 2\rangle\md
G\langle 2\rangle$ is equidimensional.
This is a particular case of the theorem of Eisenbud-Frenkel mentioned in
the introduction.

\begin{comment}
%%%%%%%%%%%%%%%%%%%%%%%%%%%%%%%%%%%%%%%%%%%%%%%%%%%
Since the Lie algebra $\g\langle 2\rangle$ is quadratic, the following
is a standard consequence of Theorems~\ref{gen-takiff}(ii) and 
\ref{eq-gen-takiff} (cf. \cite{ko63, geof1}):

\begin{s}{Corollary}
The enveloping algebra\/ ${\goth U}(\g\langle 2\rangle)$ is a free module over its
centre.
\end{s}%
%
%%%%%%%%%%%%%%%%%%%%%%%%%%%%%%%%%%%
\end{comment}

\end{document}